\documentclass[a4paper,12pt]{amsart}

\usepackage{hyperref}

\usepackage[top=3cm,left=2.5cm,right=2.5cm,bottom=3cm]{geometry}
\usepackage{relsize}
\usepackage{lineno}

\usepackage{amsmath,amsthm,pb-diagram,amssymb,comment}
\usepackage{amsfonts,graphicx,color}
\usepackage{enumitem, fancyhdr, dsfont, multicol}

\usepackage{thmtools,cleveref}
\usepackage[normalem]{ulem}
\usepackage{stmaryrd}
\usepackage{textcomp}
\usepackage{mathtools}
\usepackage{fontawesome}
\usepackage{textcomp}
\usepackage{tikz,float}
\usepackage{mleftright}
\usepackage{mathrsfs}
\usepackage{tikz-cd}
\usepackage{dsfont}
\usepackage{stmaryrd}

\hypersetup{
    colorlinks=true,
    linkcolor=magenta,
    filecolor=magenta,
    urlcolor=magenta,
    citecolor=magenta,
}
\definecolor{caribbeangreen}{rgb}{0.0, 0.8, 0.6}

\setlist{topsep=0ex,itemsep=1ex}

    \DeclareMathOperator{\dom}{{\rm dom}}

    \DeclareMathOperator{\ran}{{\rm ran}}

    \newcommand{\menos}{\smallsetminus}
    
    \DeclareMathOperator{\pts}{\mathcal{P}}
    
    \newcommand{\frestr}{{\upharpoonright}}

     \newcommand{\M}{\mathcal{M}}
    \newcommand{\Mcal}{\mathcal{M}}

    \newcommand{\N}{\mathbb{N}}
    \newcommand{\Z}{\mathbb{Z}}
    
    \newcommand{\R}{\mathbb{R}}

    \newcommand{\la}{\langle}
    \newcommand{\ra}{\rangle}

\DeclareMathOperator{\Lb}{\mathrm{Lb}}

\newcommand{\id}{\mathrm{id}}

   \DeclareMathOperator{\supsum}{\overline{\rm{S}}}
    \DeclareMathOperator{\infsum}{\underline{\rm{S}}}

\newcommand{\set}[2]{\left\{#1 \colon #2\right\}}
\newcommand{\Seq}[2]{\la #1 \colon #2\ra}
\newcommand{\lseq}[2]{\left\la #1 \colon #2\right\ra}
\newcommand{\largeset}[2]{\left\{#1 \colon #2\right\}}

\newcommand{\Fn}{\mathrm{Fn}}

\newcommand{\cantor}{{}^\omega2}

\newcommand{\Pbf}{\mathbf{P}}

\newcommand{\redq}[1]{{\color{red} #1}}

\newcommand{\ZFC}{\mathrm{ZFC}}

\newcommand{\At}{\mathrm{At}}

\newcommand{\bfC}{\mathbf{C}}

\newcommand{\bfP}{\mathbf{P}}

\newcommand{\bbQ}{\mathbb{Q}}
\newcommand{\bbR}{\mathbb{R}}

\newcommand{\cA}{\mathscr{A}}
\newcommand{\cB}{\mathscr{B}}
\newcommand{\cC}{\mathscr{C}}

\newcommand{\cI}{\mathscr{I}}

\newcommand{\cM}{\mathscr{M}}

\newcommand{\calB}{\mathcal{B}}
\newcommand{\calP}{\mathcal{P}}

\newcommand{\calX}{\mathcal{X}}

\newcommand{\varp}{\varepsilon}

\newcommand{\rest}{{\restriction}}

\newenvironment{PROOF}[2][\proofname.]
   {\begin{proof}[#1]\renewcommand{\qedsymbol}{\textsquare$_{\mathrm{#2}}$}}
   {\end{proof}}

\renewcommand{\setminus}{\smallsetminus}

\DeclareMathOperator{\stone}{\mathrm{St}}
\DeclareMathOperator{\osc}{\mathrm{osc}}
\newcommand{\Sm}{\mathrm{S}}
\newcommand{\cJ}{\mathscr{J}}

\DeclareMathOperator{\cp}{\mathrm{cp}}

\newcommand{\Andres}[1]{{\color{purple}Andres says: #1}}

\usepackage{accents}
\MakeRobust{\underaccent} 
\newcommand{\unbar}[1]{\underaccent{\bar}{#1}}
\definecolor{sub0}{RGB}{29,32,137}
\definecolor{sub1}{RGB}{1,71,157}
\definecolor{sub2}{RGB}{1,104,183}
\definecolor{sub3}{RGB}{0,160,234}
\definecolor{sug}{RGB}{0,154,68}
\definecolor{suy}{RGB}{208,219,1}

\newcommand{\subo}[1]{{\color{sub0}#1}}

\definecolor{redun}{rgb}{0.65, 0.11, 0.19}
\definecolor{greenun}{rgb}{0.58, 0.71, 0.23}
\definecolor{dodger}{rgb}{0.0,0.5,1.0}
\definecolor{carrotorange}{rgb}{0.93, 0.57, 0.13}

\definecolor{bluet}{rgb}{0.0, 0.4, 0.6}

\hypersetup{
    colorlinks=true,
    linkcolor=bluet,
    filecolor=bluet,
    urlcolor=bluet,
    citecolor=bluet,
}
    
    \newcommand{\DM}[1]{{\color{blue}Diego says: #1}}

    \newcommand{\MC}[1]{{\color{carrotorange}Miguel says: #1}}

\makeindex

\title{Finitely additive measures on Boolean algebras}

\author{Miguel A. Cardona}
\address{Faculty of Engineering, Instituci\'on Universitaria Pascual Bravo. Calle 73 No.~73A -- 226, Medell\'in, Colombia;
and\newline 
Einstein Institute of Mathematics, 
Edmond J. Safra Campus, Givat Ram\\
The Hebrew University of Jerusalem\\
Jerusalem, 91904, Israel}
\email{\href{mailto:miguel.cardona@pascualbravo.edu.co}{miguel.cardona@pascualbravo.edu.co}}
\urladdr{\url{https://sites.google.com/view/miacardonamo}}

\author{Diego A.~Mej\'ia}
\address{Graduate School of System Informatics, Kobe University. 1-1 Rokkodai-cho, Nada-ku, Kobe, Hyogo 657-8501 Japan}
\email{\href{mailto:damejiag@people.kobe-u.ac.jp}{damejiag@people.kobe-u.ac.jp}}
\urladdr{\url{https://researchmap.jp/mejia?lang=en}}

\author{Andr\'es F. Uribe-Zapata}
\address{TU Wien, Faculty of Mathematics and Geoinformation, Institute of Discrete Mathematics and Geometry, Wiedner Hauptstrasse 8--10, A--1040 Vienna, Austria }
\email{\href{mailto:andres.zapata@tuwien.ac.at}{andres.zapata@tuwien.ac.at}}
\urladdr{\url{https://sites.google.com/view/andres-uribe-afuz}}

\date{\today}

\begin{document}

\makeatletter
\def\@roman#1{\romannumeral #1}
\makeatother

\newcounter{enuAlph}
\renewcommand{\theenuAlph}{\Alph{enuAlph}}

% \numberwithin{equation}{section}
% \renewcommand{\theequation}{\thesection.\arabic{equation}}

\theoremstyle{plain}
  \newtheorem{theorem}{Theorem}[section]
  \newtheorem{corollary}[theorem]{Corollary}
  \newtheorem{lemma}[theorem]{Lemma}
  \newtheorem{mainlemma}[theorem]{Main Lemma}
  \newtheorem{mainproblem}[theorem]{Main Problem}
  \newtheorem{construction}[theorem]{Construction}
  \newtheorem{prop}[theorem]{Proposition}
  \newtheorem{claim}[theorem]{Claim}
  \newtheorem{fact}[theorem]{Fact}
  \newtheorem{exer}[theorem]{Exercise}
  \newtheorem{question}[theorem]{Question}
  \newtheorem{problem}[theorem]{Problem}
  \newtheorem{cruciallem}[theorem]{Crucial Lemma}
  \newtheorem{conjecture}[theorem]{Conjecture}
  \newtheorem{assumption}[theorem]{Assumption}
  \newtheorem*{thm}{Theorem}
    \newtheorem*{thma*}{Main Lemma A}
    \newtheorem*{thmb*}{Theorem B}
    \newtheorem*{thmc*}{Theorem C}
    \newtheorem*{thmd*}{Theorem D}
  \newtheorem{teorema}[enuAlph]{Theorem}
  \newtheorem{pregunta}[enuAlph]{Question}
  \newtheorem*{corolario}{Corollary}
\theoremstyle{definition}
  \newtheorem{definition}[theorem]{Definition}
  \newtheorem{example}[theorem]{Example}
  \newtheorem{remark}[theorem]{Remark}
    \newtheorem{hremark}[theorem]{Historical Remark}
    \newtheorem{observation}[theorem]{Observation}
  \newtheorem{notation}[theorem]{Notation}
  \newtheorem{context}[theorem]{Context}

  \newtheorem*{defi}{Definition}
  \newtheorem*{acknowledgements}{Acknowledgements}

\numberwithin{equation}{theorem}
\renewcommand{\theequation}{\thetheorem.\arabic{equation}}

\def\sectionautorefname{Section}
\def\subsectionautorefname{Subsection}

\begin{abstract}
    In this article, we conduct a detailed study of \emph{finitely additive measures} (fams) in the context of Boolean algebras, focusing on three specific topics: freeness and approximation, existence and extension criteria, and integration theory. In the first topic, we present a classification of \emph{free} finitely additive measures, that is, those for which the measure of finite sets is zero, in terms of approximation to uniform probability measures. This inspires a weaker version of this notion, which we call the \emph{uniform approximation property}, characterized in terms of freeness and another well-determined type of fams we call \emph{uniformly supported}. In the second topic, we study new criteria for the existence and extension of finitely additive measures. In particular, we provide an extension of the so-called \emph{compatibility theorem}---which characterizes when two finitely additive measures can be extended---yielding a precise and compact characterization of when three finitely additive measures can be simultaneously extended, under the assumption that one of them is an ultrafilter.
    
    Finally, we study a Riemann-type integration theory on fields of sets with respect to finitely additive measures, allowing us to extend and generalize some classical concepts and results from real analysis, such as Riemann integration over rectangles in $\mathbb{R}^{n}$ and the Jordan measure. We also generalize the extension criteria for fams allowing desired values of integrals of a given set of functions. At the end, we explore the connection between integration in fields of sets and the Lebesgue integration in the Stone space of the corresponding field, where we establish a characterization of integrability in the sense of the Lebesgue-Vitali theorem, which follows as a consequence of our results.
\end{abstract}

%\keywords{Finitely additive measure, Boolean algebra, extension criteria, integration theory.}

\thanks{The first author was partially supported by the Slovak Research and Development Agency under Contract No.~APVV-20-0045 and by Pavol Jozef \v{S}af\'arik University in Ko\v{s}ice at a postdoctoral position; 
The second author was partially supported by the Grant-in-Aid for Scientific Research (C)  23K03198, Japan Society for the Promotion of Science; and the third author was supported by the Austrian Science Fund (FWF): project number P33895.
} 

\thanks{The authors thank Professor Piotr Borodulin-Nadzieja for hosting the visits of the second and third authors at the University of Wroc{\l}aw and for providing valuable references that enriched this work.}



% 28A60 Measures on Boolean rings, measure algebras

% 03G05 Logical aspects of Boolean algebras

% 06E05 Structure theory of Boolean algebras

% 06E15 Stone spaces (Boolean spaces) and related structures

% 06E99 None of the above, but in this section (Boolean algebra and Boolean rings)

%28A25 Integration with respect to measures and other set functions

%54H99 None of the above, but in this section (Connections of general topology with other structures, applications)

\date{\today}

\maketitle

{
\small
\hypersetup{linkcolor=black}
\tableofcontents
}

\section{Introduction}

The notion of \emph{measure} is a central concept in mathematics that has given rise to a complete field of study: Measure Theory. %It is dedicated to its formalization and subsequent study. 
This theory not only formalizes, extends, and generalizes some natural and intuitive concepts such as length, area, and volume, but also allows for the addressing of problems ranging from Archimedes' geometric methods---such as his \emph{method of exhaustions}, which laid the foundations for the origins of the current calculus---to more modern and complex issues in several areas of mathematics, such as mathematical analysis, functional analysis, probability theory, stochastic processes, and integration theory, among many other disciplines. In these fields, the notion of measure provides the foundation for the formalization and rigorous study of structures where notions of measurability exist.

While standard measure theory---that is, the theory based on $\sigma$-additive measures---is extremely broad and powerful, there are abstract contexts in which it becomes overly restrictive or even inappropriate to use them. One of such examples arises in \emph{Game Theory}. For instance, in~\cite{CoreGames}, a case is presented where, in games with an infinite number of decisions or players, using weak forms of measures is not only useful but, in fact, necessary. Another example that highlights this can be found in the context of topology. As mentioned in~\cite{pivato}, it is often not possible to construct measures on $\mathfrak{R}(S)$, the collection of regular open subsets of a topological space $S$. This difficulty arises because $\mathfrak{R}(S)$ generally lacks the necessary structural properties to support a well-defined measure (see also~\cite{Fremlin74}). For this reason, it is necessary to consider weaker forms of measures or alternative mathematical frameworks when dealing with such spaces. The impossibility of using measures in these and other contexts motivates the introduction of a concept that extends the traditional framework of measure theory. This concept---known as \emph{finitely additive measure}, abbreviated \emph{fam}---provides a natural way to weaken the standard notion of measure by weakening the requirement of $\sigma$-additivity, enforcing it just for a finite number of disjoint sets. This generalization not only broadens the scope of measures, allowing their use in situations where $\sigma$-additive measures are not applicable, but has also led some mathematicians to regard them as even more interesting and natural objects than standard measures. For example, in the words of Dorothy Maharam: ``\emph{Some years ago S.~Bochner remarked to the author that finitely additive measures are more interesting, and perhaps more important, than countably additive ones. Certainly there has been increasing interest in them shown by mathematicians and statisticians}'' (see~\cite{Maharam1976}). 

Beyond their presence in certain contexts of Game Theory, finitely additive measures play a significant role in various areas of mathematics. In Functional Analysis, these measures arise naturally in the duality theory of Banach spaces, where they provide fundamental tools for the study of linear operators and functionals (see~\cite{VectorMeasures}). In Real Analysis, they can be used to characterize Riemann integration over rectangles in $\mathbb{R}^n$ and to prove that this integral is absolute for transitive models of $\ZFC$ (see \autoref{t16} and~\cite{PU2025}). In the field of Mathematical Logic, finitely additive measures have multiple applications (see e.g.~\cite{Piotr1},~\cite{Piotr2},~\cite{VANDOUWEN1992223} and~\cite{Dzamonja2006}), in particular because they generalize the notion of \emph{ultrafilter} (which is, essentially, a $0$-$1$-finitely additive measure). Ultrafilters are particularly useful in model theory for constructing structures such as ultraproducts and ultrapowers, and on the other hand, it allows the characterization of Boolean algebras through the well-known \emph{Stone Representation Theorem}. 
This theorem establishes deep connections between Boolean algebras and topology (see \autoref{stone}). In the study of cardinal invariants of the continuum, ultrafilters have proven to be a powerful tool for solving fundamental problems. A notable example is the work of Martin Goldstern, the second author, and Saharon Shelah, who in~\cite{GMS} used ultrafilters to separate the left side of Cichoń's diagram. This achievement represented a crucial first step toward the proof of \emph{Cichoń's Maximum}, a result later proven by Goldstern, Kellner, the second author, and Shelah (see~\cite{GKS} and~\cite{GKMS}), which is considered one of the most significant advances in the combinatorics of real numbers to date. 

In the context of fields of sets,~\cite{BhaskaraRa} develops an extensive and well-known study of finitely additive measures---which they call \emph{charges}---providing a detailed overview of their properties and their significance in various mathematical settings. In particular, it is shown how these measures can be used to extend classical results in analysis that would not be possible within the framework of $\sigma$-additive measures. 

The main objective of this work is to develop a detailed study of finitely additive measures within the framework of Boolean algebras and fields of sets, focusing on three specific topics: freeness and approximation (\autoref{2}), existence and extension criteria (\autoref{3}), and integration theory (\autoref{3.5}--\ref{stone}). %To this end, in some cases, we develop new tools such as the content presented in \autoref{2.3}. In others, we adapt existing results from other areas and make the necessary modifications for them to be applicable in our Boolean setting, as in the case of \autoref{m48}, whose version for field of sets can be found in~\cite[Ch.~3]{BhaskaraRa}.
The motivation for this arises from the study of \emph{cardinal invariants} and \emph{Forcing Theory}. In 2000, Saharon Shelah, in~\cite{Sh00}, introduced a powerful forcing technique using finitely additive measures on $\calP(\omega)$, which solved a problem of cardinal invariants that had remained open for at least 30 years (see~\cite[Sec.~5.1]{uribethesis}). Later, Jakob Kellner, Saharon Shelah, and Anda T\u{a}nasie, in~\cite{KST}, revisited Shelah’s method and applied it to other problems related to cardinal invariants. In 2023, the third author’s master’s thesis expanded some of the ideas developed in the previous works and---since in \cite[Sec.~1]{Sh00} Shelah concisely presents the properties of finitely additive measures on $\calP(\omega)$ that he will use throughout his work, but without providing proofs or additional details---included an entire chapter dedicated to the study of finitely additive measures on Boolean algebras, which served as the initial inspiration and motivation for the development of this article (see~\cite[Ch.~3]{uribethesis}). Finally, in 2024, the authors in~\cite{CMU}, building on the aforementioned works, introduced a general theory of iterated forcing using finitely additive measures on Boolean algebras. In this work, it was necessary to use many (particular cases) of tools presented here, ultimately consolidating the content of this article, which is fundamentally divided into three parts as we next describe in detail. 

In the first part, corresponding to \autoref{2}, we introduce the fundamental notions of Boolean algebras, define the concept of \emph{finitely additive measure} in this context, and study some of its fundamental properties. In general terms, a finitely additive measure on a Boolean algebra $\cB$ is a function $\Xi \colon \cB \to [0, \infty]$ that satisfies the property of additivity for a finite number of sets, that is, $
\Xi(a \vee b) = \Xi(a) + \Xi(b) $ 
whenever $a, b \in \cB$ and $a \wedge b = 0_{\cB}$. In many contexts, when $\cB$ is a field of sets, finitely additive measures that assign positive values to finite sets exhibit pathological behavior. For this reason, it is common to assume---or even incorporate into the definition itself---that the measure of any finite set in $\cB$ is zero, but we develop our work without this restriction. In connection with the notion of \emph{free ultrafilter} (i.e.\ ultrafilters not containing finite sets), the concept of \emph{free finitely additive measure} formalizes the previous idea (see \autoref{m4}). In~\cite{Sh00} (see also~\cite[Lem.~1.2]{KST}) it is shown that free finitely additive measures can be locally approximated to the uniform measure. In the following result, we generalize and reformulate this fact (corresponding to \autoref{m420} and~\ref{m421}).
%We present a classification of this kind of measure through the following result, corresponding to \autoref{m420}. 

\begin{teorema}\label{in0}
    Let $\cB$ be a Boolean subalgebra of $\pts(X)$ and let $\Xi$ be a finite fam on $\cB$ with $\delta\coloneqq \Xi(X)>0$. Then, for any 
    any $\varp>0$ and any finite partition $P$ of $X$ by sets in $\cB$, there is some non-empty finite $u\subseteq X$ of size ${\leq}\min\left\{|P|,\left\lceil\frac{\delta}{\varp}\right\rceil\right\}$ and a probability measure $\mu$ on $\pts(u)$ such that $|\delta\mu(b) - \Xi(b)|<\varp$ for all $b\in P$.
    
    Moreover, the following statements are equivalent.
    \begin{enumerate}[label = \rm (\roman*)]
        \item\label{m420i.in} For any $\varp>0$, any finite $F\subseteq X$ and any finite partition $P$ of $X$ by sets in $\cB$, there is some non-empty finite $u\subseteq X\menos F$  such that, for all $b\in P$,
        \[\left|\delta\frac{|u\cap b|}{|u|} - \Xi(b)\right| < \varp.\]
        
        \item\label{m420ii.in} All finite sets in $\cB$ have measure zero.
    \end{enumerate}
    In~\ref{m420i.in}, the size of $u$ can be bounded by $\left\lceil\frac{\delta}{\varp}\right\rceil$.    
\end{teorema}

If we remove the finite set $F$ in \autoref{in0}~\ref{m420i.in}, we obtain a weaker form of freeness, which we call the \emph{uniform approximation property} (uap). Our main result in this section provides a classification of the finitely additive measures that satisfy the uap, also using a type of finitely additive measure that we call \emph{uniformly supported}. 

\begin{teorema}\label{in1}
    Let $\Xi$ be a finite fam on a Boolean subalgebra $\cB$ of $\pts(X)$ and $\delta\coloneqq \Xi(X)$. Then the following statements are equivalent.
    \begin{enumerate}[label = \rm (\roman*)]
        \item $\Xi$ has the uap.
        
        \item\label{in1.2} Either all finite sets in $\cB$ have measure zero, or $\Xi$ is \emph{uniformly supported}, i.e.\ there is some natural number $d>0$ such that, for any $B\in\cB$, $\Xi(B)= \frac{\ell}{d}\delta$ for some natural $0\leq \ell\leq d$, and $\ell\leq |B|$ when $B$ is finite.
    \end{enumerate}
\end{teorema}

In general, without assuming that $\cB$ is a subalgebra of $\calP(X)$, we can represent the uniformly supported finitely additive measures as a linear combination of finitely additive measures obtained from ultrafilters (see \autoref{m410}, this is the main reason for the name ``uniformly supported"). 

In the second part, which is covered in \autoref{3}, we explore existence and extension criteria for finitely additive measures. More precisely, we address problems of the following type: given a Boolean algebra $\cB$ and a function $f$ whose domain is a subset of $\cB$ containing the top element $1_{\cB}$ of $\cB$, can $f$ be extended to a finitely additive measure on $\cB$ with certain desired properties? A classical instance of this problem arises when $f$ is itself a finitely additive measure $\Xi$ defined on a subalgebra of $\cB$ (see \autoref{f700}). A well-known example that is closely related is the so-called \emph{Compatibility Theorem} for fams (see e.g. \cite{Guy1961}, \cite{Kindler88} and \cite{BhaskaraRa}). This theorem characterizes when \emph{two} finitely additive measures $\Xi_{0}$ and $\Xi_{1}$, defined on subalgebras $\cB_{0}$ and $\cB_{1}$ of $\cB$, respectively, can be extended by a single measure. While this result appears in \cite[Ch.~3]{BhaskaraRa} for fields of sets, we provide a self-contained and more direct proof for Boolean algebras in \autoref{m48}. One particular case is the extension of two compatible finitely additive measures when one of them corresponds to a filter (i.e.\ a finitely additive measure onto $\{0,1\}$), presented in Shelah's~\cite[Prop.~1.3~(3)]{Sh00} and later in  \cite[Fact.~1.3~(a)]{KST}, both without proof), which we prove in \autoref{m69}. This stands out as one of the most significant results in the development of~\cite{CMU}, specifically, it serves as the main tool in the proof of the central theoretical result presented there (see~\cite[Sec.~7.3]{CMU}).

It is well known that---as an alternative to the proof given in \cite{BhaskaraRa} or to our own proof of the compatibility theorem, which relies on linear algebraic arguments---this result can also be obtained using the so-called Sandwich Theorems (see, e.g., \cite[Prop.~2.2]{Hansel86} and \cite{Kindler88}). Moreover, inspired by this type of argument, it is possible to generalize the compatibility theorem to any number (even infinite) of fams, as we present in \autoref{f701}. The compatibility theorem then appears as a particular case of this more general statement, since in the case of only two fams the conditions in \autoref{f701} admit a very clear and elegant characterization. However, when dealing with three or more fams, these conditions can no longer be stated in a similarly transparent and compact form.
Such situations arise, for example, in bi- or tri-dimensional forcing constructions, which play a significant role in applications to infinite combinatorics. 

Thus, despite the general formulation given in \autoref{f701}, the question remains whether cleaner and more workable characterizations can be obtained. In the particular case of filters (i.e., finitely additive measures taking values in ${0,1}$), the answer is well known and is provided by the \emph{finite intersection property} (see, e.g., \cite{BellM}). However, in the context of more general finitely additive measures, the problem becomes significantly more complex. One of the main results of this work---presented in \autoref{m80}---provides such a characterization for three finitely additive measures when one of them is a filter. The result establishes:

% However, in none of these cases is it clear how to generalize it 
% in a clean and well-defined manner to situations involving three or more measures. 

%Si bien en \autoref{f701}, ofrecemos una caracterización para extender una cantidad finita de fams, este criterio 

% Hence, the following question arises:

% \begin{pregunta}\label{in20}
%     Is there a clear and general compatibility theorem for three or more finitely additive measures?
% \end{pregunta}

\begin{teorema}\label{in25}
    For $e\in\{0,1\}$, let $\Xi_e$ be a finitely additive measure on a Boolean subalgebra $\cB_e$ of $\cB$, and let $\langle  b_i \colon i\in I\rangle\subseteq\cB$. Assume that $\Xi_0(1_\cB)= \Xi_1(1_\cB)=\delta\in(0,\infty)$. Then, the following statements are equivalent:
    \begin{enumerate}[label=\normalfont(\roman*)]
        \item\label{in.m80i} There is a finitely additive measure $\Xi$ extending $\Xi_0\cup\Xi_1$ such that $\Xi(b_i)=\delta$ for all $i\in I$.

        \item\label{in.m80ii} The following conditions are satisfied:
          \begin{enumerate}[label = $(\bullet_{\arabic*})$]
              \item\label{in.m80b1} For any $e\in\{0,1\}$, $b\in\cB_e$ and $J\subseteq I$ finite, if $\Xi_e(b)>0$ then $b\wedge\bigwedge_{i\in J}b_i\neq 0_\cB$.
              \item\label{in.m80b2} For any $a\in\cB_0$, $b\in \cB_1$ and $J\subseteq I$ finite, if $\Xi_0(a)>0$ and $a\wedge\bigwedge_{i\in J}b_i \leq b\wedge\bigwedge_{i\in J}b_i$ then $\Xi_0(a)\leq \Xi_1(b)$.
          \end{enumerate}
    \end{enumerate}
\end{teorema}

Finally, in the third part, corresponding to \autoref{3.5}--\ref{stone}, we study integration theory on fields of sets with respect to finitely additive measures. It is well known that in this context, integration can be developed in two distinct ways. The first approach, presented in detail in \autoref{3.5}, follows ideas from the Riemann integral on the real line. Given a field of sets $\cB \subseteq \calP(X)$ (that is, a subalgebra under the usual set-theoretic operations), a bounded function $f \colon X \to \bbR$, and a finitely additive measure $\Xi \colon \cB \to [0, \infty)$, we define the \emph{integral of $f$ on $X$ with respect to $\Xi$}, denoted by $\int_{X} f d\Xi$, whenever $f$ is \emph{$\Xi$-integrable}. These notions appear in~\cite[Sec.~4.5]{BhaskaraRa} under the name \emph{$S$-integrability} (see also \cite{Luxemburg1991} and \cite{Toland20}). This integral preserves many of the basic properties of the classical Riemann integral, as shown in \autoref{3.5} (see, for instance, \autoref{t47},~\ref{t55},~\ref{t58}, and~\ref{t29}).  The second approach, developed in \autoref{stone}, is based on the Lebesgue measure and the Stone representation theorem (see also \cite[Sec .~8]{pivato}).

Several results in the literature connect these two types of integration in slightly different contexts. For instance, Pivato and Vergopoulos~\cite{pivato} study the integration theory on the Boolean algebra of regular open sets of a compact space endowed with a finitely additive measure. However, this algebra is not a field of sets, since the join $\vee$ and complement $\sim$ operations do not coincide with the set-theoretic union and complement. Therefore, their approach differs from the framework developed in this work. One of the main results presented here shows that it is not only possible to relate the two integration approaches discussed above, but also to characterize $\Xi$-integrability in terms of Lebesgue-type integration properties on the corresponding Stone space. Recall that any Boolean algebra $\cB$ is isomorphic to the Boolean algebra of clopen subsets of its Stone space $\stone(\cB)$ via the map $a \mapsto [a] \coloneqq \{ u \in \stone(\cB) \colon a \in u \}$. The finitely additive measure $\Xi$ can then be naturally transferred to the algebra of clopen subsets of $\stone(\cB)$ by setting $\mu^\Xi_-([a]) \coloneqq \Xi(a)$ for $a \in \cB$. Since $\stone(\cB)$ is compact, $\mu^\Xi_-$ is $\sigma$-additive on the algebra of clopen sets. Hence, there exists a unique $\sigma$-additive measure $\mu^\Xi$ on the $\sigma$-algebra of $\stone(\cB)$ generated by the clopen sets such that $\mu^\Xi$ extends $\mu^\Xi_-$. Denote by $\cM^\Xi$ the completion of this $\sigma$-algebra with respect to $\mu^\Xi$.
If $\cB$ is a field of sets over $X$ and $f \colon X \to \bbR$ is a bounded function, we define a function $\hat{f}$ on the set of all ultrafilters $u \in \stone(\cB)$ for which $f$ is $\Xi_u$-integrable by $\hat{f}(u) \coloneqq \int_{X} f\, d\Xi_{u}$, where $\Xi_{u} \colon \cB \to \{0,1\}$ is given by $\Xi_{u}(b) = 1$ if $b \in u$ and $\Xi_{u}(b) = 0$ otherwise (see \autoref{2.2}). The function $\hat{f}$ is continuous on its domain. Inspired by the classical Lebesgue–Vitali theorem---which characterizes Riemann integrability in terms of the measure of the set of discontinuities---our main result in \autoref{stone} (see \autoref{u10}) extends this idea to the setting of $\Xi$-integration.

\begin{teorema}\label{in16}
    Let $f\colon X\to \bbR$ be bounded. Then $f$ is $\Xi$-integrable iff  the complement of $\dom\hat f$ has $\mu^\Xi$-measure zero.
\end{teorema}

This result has consequences in the integration theory of any compact zero-dimensional space $S$. If $\mu$ is a measure on (the completion of) the $\sigma$-algebra generated by the clopen sets, and $\mu_-$ is the restriction of $\mu$ to the algebra of clopen sets, then any bounded real-value function on $S$ is $\mu_-$-integrable iff it is continuous $\mu$-almost everywhere (the Lebesgue-Vitali Theorem for compact zero-dimensional spaces). As a consequence, after analyzing the integration theory on the Cantor space and its relation to Riemann integration on $[0, 1]$, we obtain a new proof of the well-known Lebesgue-Vitali theorem: a bounded function $g \colon [0, 1] \to \bbR$ is Riemann integrable if and only if it is continuous almost everywhere (see \autoref{u50}).  

%we have the actual Lebesgue-Vitali theorem for $\R$  as a corollary after analyzing the integration theory on the Cantor space and its relation to Riemann integration on $[0, 1]$ %, we obtain a new proof of the well-known Lebesgue-Vitali theorem: a bounded function $g \colon [0, 1] \to \bbR$ is integrable if and only if it is continuous almost everywhere 

The reader might wonder why the authors take the trouble to develop \autoref{3.5} in such detail, given that the traditional approach is to lift to the Lebesgue measure in the Stone space. However, several results established there have proven to be highly useful. For instance, the integration approach developed there is precisely the one used in~\cite{Sh00}, \cite{KST}, and~\cite{CMU} to construct a general theory of iterated forcing using finitely additive measures, where relying on the Lebesgue measure and the Stone space would be unnecessarily difficult and technically cumbersome. More recently, the third author, together with Parra-Londo\~no, used the characterization given in \autoref{t16}, along with several other results (see \autoref{s5}, \autoref{s2}, and \autoref{da1}), to prove that the Riemann integral on $\bbR^{n}$ is absolute for transitive models of $\ZFC$ (see \cite{PU2025}). This framework is also useful for reformulating various approximation and extension criteria for finitely additive measures involving integrability conditions, which have proven to be fundamental for \cite{CMU} and which we develop further in \autoref{7}.

In \autoref{5}, we study theorems on limits and convergence within the integration theory developed in \autoref{3.5}, extending well-known results for $\sigma$-additive measures to the setting of finitely additive measures. While some results of this type appear in various contexts---particularly in functional analysis and related areas---it is often unclear which statements apply directly to finitely additive measures, and the literature does not provide a unified treatment. As such, the theorems presented here, which appear to be new in comparison with the references cited in this work, are valuable: they provide a clean and direct formulation, and they proves that $\sigma$-additivity is not strictly necessary to obtain many standard results on limits and convergence (see \autoref{t66} and \autoref{t67}). Finally, we use this theorem as the essential technical tool to generalize a classical result from integration theory over $\bbR^{n}$ about Jordan measurable sets: \autoref{in14}. For this, we study integration over subsets of $X$ using characteristic functions to be able to introduce the notion of \emph{Jordan measurable sets} for such fields and define the associated \emph{Jordan measure} $\hat{\Xi}$, which extends the original finitely additive measure $\Xi$. Based on this, we define the collection $\cJ^{\Xi}$ of Jordan measurable sets (see \autoref{t311}), which forms a Boolean subalgebra of $\calP(X)$ naturally extending $\cB$. The main result of this section shows that integration with respect to $\Xi$ is equivalent to integration with respect to $\hat{\Xi}$, where as said the theorems on limits and convergence play a fundamental role.

\begin{teorema}\label{in14}
    Let $f\colon X\to\bbR$ be bounded. Then $f$ is $\Xi$-integrable if and only if it is $\hat\Xi$-integrable, in which case $\int_X f \, d \Xi = \int_X f d \hat\Xi$.
\end{teorema}

\section{Basic notions of finitely additive measures}\label{2}

In this section, we review the notion of  \emph{finitely additive measure} on a Boolean algebra and present some of its fundamental properties. In particular, we define the concept of \emph{freeness}, which characterizes those finitely additive measures that assign zero to finite sets, as well as the \emph{uniform approximation property}, which will allow us to establish approximation criteria for certain finitely additive measures. For this, we start by reviewing some basic concepts and notation about Boolean algebras generally following~\cite{Halmos} and~\cite{BellM}.

\subsection{Boolean algebras}\label{2.1}

Let $\cB$ be a Boolean algebra. We denote $\cB^{+} \coloneqq \cB \setminus \{ 0_{\cB} \}.$ For $B \subseteq \cB,$ the \emph{Boolean subalgebra generated by $B$}, denoted by $\langle B \rangle_{\cB}$, is the $\subseteq$-smallest Boolean subalgebra of $\cB$ containing $B.$ In this case, $B$ is called the \emph{generator set} of $\langle B \rangle_{\cB}.$ If $B$ is a finite set, we say that $\langle B \rangle_{\cB}$ is a \emph{finitely generated Boolean algebra}. It is well-known that $\cB$ can be endowed with a partial order: for any $a, b \in \cB,$ $ a \leq_{\cB} b $ iff $ a \wedge b = a.$  Recall that, if $I$ is an ideal on $\cB,$ then we can define the quotient between $\cB$ and $I,$ which we denote by $\cB / I.$   An \emph{atom of $\cB$} is a $\leq_{\cB}$-minimal non-zero element of $\cB.$  Denote by $\At_{\cB}$ the set of all atoms of $\cB.$  We are particularly interested in the atomic structure of finitely generated Boolean algebras, which we describe below. For this, we introduce the following definition but recall before that for sets $I$ and $J$, we denote
\[\Fn(I,J)\coloneqq\set{\sigma}{\text{$\sigma$ is a finite partial function from $I$ into $J$}}.\]

\begin{definition}\label{b41}
    Let %$\cB$ be a Boolean algebra and 
    $B \subseteq \cB$. %Then, 

    \begin{enumerate}[label=\normalfont(\arabic*)]
        \item For any $b \in \cB$ and $d \in \{ 0, 1 \}$ we define:  
        
        $$b^{d} \coloneqq 
        \left\{ \begin{array}{ll}
             b &   \text{if $d = 0$,} \\[1ex]
             {\sim }b &  \text{if  $d = 1$.}
            \end{array}
        \right.$$ \index{$b^{d}$}

        \item For $\sigma \in \Fn(B, 2)$, we define $ a_{\sigma} \coloneqq \bigwedge_{b \in \dom(\sigma)} b^{\sigma(b)}$, which clearly is al element of $\langle B \rangle_{\cB}$.  
    \end{enumerate}
\end{definition}

In the case when $\cB$ is finitely generated by some finite subset $B$, we have that: $$\At_{\cB} = \{ a_{\sigma} \colon \sigma \in {}^{B}2 \wedge a_{\sigma} \neq 0_\cB \}.$$

Moreover, $\At_{\cB}$ is finite, $\vert \At_{\cB} \vert \leq 2^{\vert B \vert}$, and $|\cB|=2^{|\At_{\cB}|}$ (so $\cB$ is finite). The latter is connected to the fact that $\calP(\At_{\cB})$ and $\cB$ are isomorphic via the map $A\mapsto \bigvee A$.

Recall that any element of $\cB$ is in the subalgebra generated by a set $B$ iff it can be written as a finite join of finite meets of
elements and complements of elements from $B$, that is, 
$$ \langle B \rangle = \left \{ \  \bigvee_{\sigma \in C}a_{\sigma} \colon C \in [\Fn(B, 2)]^{< \omega} \right \}. 
$$ 
In particular, if $B$ is finite, then 
$$ \langle B \rangle = \left   \{ \ \bigvee_{\sigma \in C} a_{\sigma} \colon C \subseteq {}^{B}2  \right \}. 
$$   

Recall that a \emph{filter on $\cB$} is a non-empty set $F \subseteq \cB$ such that:  

\begin{enumerate}[label=\normalfont(\roman*)]
    \item If $x, y \in F,$ then $x \wedge y \in F,$

    \item If $x \in F$ and $x\leq y,$ then $y \in F,$

    \item $0_{\cB} \notin F.$  

    \index{Boolean algebra! filter}
\end{enumerate} 

An \emph{ultrafilter on $\cB$} is a filter $F \subseteq \cB$ such that, for any $b \in \cB, $ either $b \in F$ or $\sim \! b \in F$. Equivalently, it is a maximal filter on $\cB$. 

For any $A\subseteq\cB$ denote $A^\sim\coloneqq \set{{\sim}a}{a\in F}$.
If $F \subseteq \cB$ is a filter, then $\langle F \rangle = F \cup F^{\sim}$. Therefore, $F$ is an ultrafilter on $\cB$ iff $\langle F \rangle = \cB.$ 

Recall that a \emph{Boolean homomorhpism} is a function $h \colon \cB \to \cC$---where $\cB$ and $\cC$ are Boolean algebras---that preserves the Boolean operations, that is,  $h(a \wedge b) = h(a) \wedge h(b),$ $h(a \vee b) = h(a) \vee h(b)$, and   $h({\sim} a) =  {\sim} h(a)$ for all $a, b \in \cB$, which implies that $h(0_{\cB}) = 0_{\cC}$ and $h(1_{\cB}) = 1_{\cC}$.\footnote{Because $h(a\vee {\sim} a) = h(a)\vee {\sim}h(a)$ and $h(a\wedge {\sim} a) = h(a)\wedge {\sim}h(a)$.} A \emph{Boolean isomorphism} from $\cB$ into $\cC$ is a bijective Boolean homomorphism from $\cB$ onto $\cC.$ 

\begin{example}\label{s16}
    Let $X$, $Y$ be non-empty sets and $h \colon X \to Y$ a function. 
    
    \begin{enumerate}[label=\normalfont(\arabic*)]
        \item\label{s16.1} The map $f_{h} \colon \calP(Y) \to \calP(X)$ defined by $f_{h}(A) \coloneqq h^{-1}[A]$ for every $A \subseteq Y$ is an homomorphism. Furthermore, $f_{h}$ is an isomorphism if, and only if, $h$ is a bijection. 

        \item More generally, if $\cB$ is a Boolean subalgebra of $\calP(X)$, then $\cC \coloneqq h^{\to}(\cB)$ is a Boolean subalgebra of $\calP(Y)$, where $ h^{\to}(\cB) \coloneqq \{ A \subseteq Y \colon h^{-1}[A] \in \cB \}$, and the map $f_{h}\frestr\cC$ %\colon \cC \to \cB$ defined by $f_{h}(c) \coloneqq h^{-1}[c]$ for every $c \in \cC$ 
        is a homomorphism into $\cB$. Furthermore: 
        \begin{enumerate}[label = \normalfont (\alph*)]
            \item $f_{h}$ is one-to-one iff $h$ is onto: If $h$ is not onto then $Y\menos \ran h\neq\emptyset$ and $h^{-1}[\emptyset]=h^{-1}[Y\menos \ran h]=\emptyset\in\cB$, so $Y\menos \ran h\in\cC$ and $f_h$ is not one-to-one. The converse ($h$ is onto implies that $f_h$ is one-to-one) is clear. 
            \item $\cC\subseteq F[\cB]$ iff $h$ is onto, where $F \colon \calP(X) \to \calP(Y)$ is defined by $F(A) \coloneqq h[A]$ for every $A \subseteq X$: if $\cC\subseteq F[\cB]$ then $Y\in F[\cB]$, which implies that $h$ is onto. Conversely, if $h$ is onto and $A\in\cC$, i.e.\ $A\subseteq Y$ and $h^{-1}[A]\in\cB$, then $A=h[h^{-1}[A]]\in F[\cB]$. 
        \end{enumerate}
    \end{enumerate}
\end{example}

By the Stone Representation Theorem, we know that every Boolean algebra is isomporhic to a Boolean subalgebra of $\calP(X)$ for some set $X$. More precisely, this $X$ is 
\[\stone(\cB)\coloneqq \set{u\subseteq \cB }{u \text{ is an ultrafilter on }\cB},\]
which is known as the \emph{Stone space of $\cB$}. This is a compact Hausdorff space endowed with the topology with basic clopen sets $[a]\coloneqq\set{u\in\stone(\cB)}{a\in u}$, in fact, these are precisely all the clopen subsets of $\stone(\cB)$. Therefore, the map $a\mapsto[a]$ is an isomorphism from $\cB$ onto the collection of clopen subsets of $\stone(\cB)$.

We define the following type of filters for Boolean algebras of the form $\calP(X)$.

\begin{definition}\label{b67}
    Let $X$ be a non-empty set. We say that $F \subseteq \calP(X)$ is a \emph{free filter} if it is a filter containing all the cofinite subsets of $X$.
\end{definition}

\subsection{Basic notions of finitely additive measures}\label{2.2}

\begin{definition}\label{m4}    
    %Let $\cB$ be a Boolean algebra. 
    A \emph{finitely additive measure on $\cB$} is a function $\Xi \colon \cB \to [0,\infty]$ satisfying:
    
    \begin{enumerate}[label=\normalfont(\roman*)]
        \item $\Xi (0_{\cB})=0$,
        
        \item $\Xi(a\vee b)=\Xi(a)+\Xi(b)$ whenever $a,b\in\cB$ and $a \wedge b= 0_{\cB}$.
    \end{enumerate}

    If $\Xi$ satisfies the following additional property, we say that it is a \emph{$\sigma$-additive measure}:

    \begin{enumerate}[resume*]
        \item If $\{ b_{n} \colon n < \omega \} \subseteq \cB$ is such that $\bigvee_{n < \omega} b_{n} \in \cB$, and $b_{i} \wedge b_{j} = 0_{\cB}$ whenever $i\neq j$, then 
        $$ \Xi \! \left( \bigvee_{n < \omega} b_{n} \right) = \sum_{n < \omega} \Xi(b_{n}).$$
    \end{enumerate}
\end{definition}

Most of the time we exclude the \emph{trivial measure}, that is,  when talking about finitely additive measures, we will always assume $\Xi(1_\cB) > 0.$ Also, we will occasionally use the acronym ``fam'' or ``FAM'' to refer to finitely additive measures.

\begin{notation}\label{n364}
    If $b \in  \cB$, then $\cB |_{b} \coloneqq \{ b \wedge b' \colon b' \in \cB \}$. This is a Boolean algebra iff $b\neq 0_\cB$ under the same operations, however, although $0_{\cB|_{b}} = 0_{\cB}$, $1_{\cB|_{b}} = b$. For this reason, $\cB|_{b}$ is not a Boolean subalgebra on $\cB$, unless $b=1_{\cB}$ (in which case $\cB|_{b} = \cB$). 
    
    If $\Xi$ is a finitely additive measure on $\cB$, we define $\Xi |_{b} \coloneqq \Xi \rest \cB |_{b}$, which is clearly a fam on $\cB |_{b}$. 
\end{notation}

There are several types of finitely additive measures. In the following definition, we introduce those that will be most relevant to us.

\begin{definition}\label{m9}
    %Let $\cB$ be a Boolean algebra and 
    Let $\Xi$ be a finitely additive measure on $\cB$. %Then,

    \begin{enumerate}[label=\normalfont(\arabic*)]
        \item   We say that $\Xi$ is \emph{finite}, if $\Xi(1_{\cB}) < \infty.$
    
        \item When $\Xi(1_\cB) = 1$ we say that $\Xi$ is a \emph{finitely additive probability measure} (abbreviated \emph{probability fam}).
        
        \item If $\Xi(b)>0$ for any $b \in \cB^{+}$, we say that $\Xi$ is \emph{strictly positive}.
    \end{enumerate}
\end{definition}

Examples of probability fams are naturally obtained from filters. Note that, whenever $\cC$ is a subalgebra on $\cB$ and $\Xi$ is a finite finitely additive measure on $\cC$ with $\Xi(1_\cB)=\delta>0$, $F_\Xi\coloneqq \set{b\in\cC}{\Xi(b)=\delta}$ is a filter on $\cC$, which can be extended to a filter on $\cB$ by closing upwards.  Conversely, if $F$ is a filter on $\cB$, then $\Xi_{F} \colon \langle F \rangle \to \{0, 1\}$ defined by 
    $$\Xi_{F}(b)= \left\{ \begin{array}{ll}
             1 &   \text{if  $b \in F$,}  \\[1ex]
             
             0 &  \text{if   $b \in F^{\sim}$,} 
             \end{array}
   \right.$$ 
is a probability fam. Furthermore if $G$ is another filter on $\cB,$ then  $$F \subseteq G \Leftrightarrow \Xi_{F} \leq \Xi_{G}.$$
This actually shows an equivalence between $\{0,1\}$-valued probability fams and ultrafilters. 

Next, we present an example of a probability measure that will appear in several places throughout this work.

\begin{example}\label{m37}
    Let $X$ be a non-empty set and fix a finite non-empty set $u \in \calP(X).$ We define $\Xi^{u} \colon \calP(X) \to [0, 1]$ such that, for any $x \in \calP(X), \, \Xi^{u}(x) \coloneqq \frac{\vert x \cap u \vert}{ \vert u \vert}$ and we call it  the \emph{uniform measure with support $u$}. It is clear that it is a probability measure ($\sigma$-additive as well). Notice that, in general, $\Xi^{u}$ is not strictly positive.
\end{example}

\begin{example}\label{m36}
   Let $\cB$, $\cC$ be Boolean algebras and $f \colon \cC \to \cB$ a Boolean homomorphism. If $\Xi$ is a finitely additive measure on $\cB,$  define the \emph{finitely additive measure induced by $f$ and $\Xi$ on $\cC$}, denoted $\Xi_{f}$, by $\Xi_{f}(c) \coloneqq \Xi(f(c))$ for all $c \in \cC$. In the particular case when $\cB$ is a Boolean subalgebra of $\calP(X)$ for some set $X$, $\Xi$ a finitely additive measures on $\cB$, and $\cC \coloneqq h^{\to}(\cB)$, we define $\Xi_{h} \coloneqq \Xi_{f_{h}}$, where $f_{h}$ is as in \autoref{s16}~\ref{s16.1}. Notice that $\Xi_{h}$ is a finitely additive measure and it is a probability fam iff $\Xi$ is, too. 
\end{example}

The following result shows one very interesting example of a strictly positive probability fam. Recall that a Boolean algebra $\cB$ is \emph{$\sigma$-centered} if $\cB^+$ can be written as a countable union of filters (or, equivalently, as a countable union of ultrafilters).

\begin{lemma}\label{m38}
    Every $\sigma$-centered Boolean algebra admits a strictly positive probability fam.
\end{lemma}

\begin{PROOF}{\ref{m38}}
    Let $\cB$ be a $\sigma$-centered Boolean algebra, so there is a countable family $\{ F_{n} \colon n < \omega \}$ of ultrafilters on $\cB$ such that $\cB^{+} = \bigcup_{n < \omega} F_{n}$. For any $b \in \cB,$ we define $\omega_{b} \coloneqq \{ n < \omega \colon b \in  F_{n} \}$. The function $\Xi \colon  \cB \to [0, 1]$ defined by
    \[\Xi(b) \coloneqq  \sum_{n \in \omega_{b}} \frac{1}{2^{n+1}}\]
    is the required fam.
\end{PROOF}

In general, to prove the existence of interesting finitely additive measures, we require the axiom of choice, that is, non-constructivist methods (see \cite{Lauwers}), like in the case of constructing interesting ultrafilters.  

We close this subsection with the following (well-known) interesting fact about fams.

\begin{lemma}\label{m50}
    Let $\Xi$ be a finite fam on a Boolean algebra $\cB$ and let $A\subseteq \cB$ such that each element of $A$ has positive measure and $a\wedge b=0_\cB$ whenever $a\neq b$ in $A$. Then, $A$ must be countable.
\end{lemma}
\begin{PROOF}{\ref{m50}}
    Pick some decreasing sequence $\Seq{\varp_n}{n<\omega}$ of positive real numbers converging to $0$. For each $n<\omega$, define $A_n\coloneqq \set{a\in A}{\Xi(a)\geq\varp_n}$. Hence $A=\bigcup_{n<\omega}A_n$.

    It is enough to show that each $A_n$ is finite. Towards a contradiction, assume that $A_n$ is infinite for some $n<\omega$. Pick some finite $F\subseteq A_n$ of size ${>}\frac{\Xi(X)}{\varp_n}$. Then, $\Xi\left(\bigvee F\right) = \sum_{a\in F}\Xi(a) \geq |F|\,\varp_n>\Xi(X)$, a contradiction.
\end{PROOF}

\subsection{Freeness and approximation}\label{2.3}

People working with finitely additive measures on Boolean subalgebras of $\calP(X)$ for some set $X$ (i.e.\ \emph{fields of sets}) are not interested, typically, in the pathological case when the measure of finite sets is positive. Although we consider the most general cases, we separate this situation as follows.
%For this reason, we introduce the following notion of \emph{freeness}. 

\begin{definition}\label{m12}
    \index{finitely additive measure!free finitely additive measure}

    If $X$ is a non-empty set, $\cB$ is a Boolean subalgebra of $\mathcal{P}(X)$ and $\Xi$ is a finitely additive measure on $\cB$, we say that $\Xi$ is a \emph{free finitely additive measure} if, for any $x\in X$, $\{x\}\in\cB$ and $\Xi(\{x\}) = 0$. 
\end{definition}

Notice that this implies that $[X]^{<\aleph_0}\subseteq\cB$ and, effectively, $\Xi(A)=0$ for any finite $A\subseteq X$. 
We adopt the name ``free finitely additive measure'' in connection with ``free filter'' (see \autoref{b67}). Note that a filter $F$ on $\calP(X)$ is free iff the measure $\Xi_F$ on $\la F\ra$ is free.

%In order to get a characterization of free finitely additive measures, we prove the following result.  

Fix the following notation about partitions.

\begin{notation}
    Let $\cB$ be a Boolean algebra. Denote by $\Pbf^\cB$ the set of finite partitions of the unity, that is, the collection of $P\subseteq \cB$ such that $\bigvee P =1_\cB$ and $b\wedge b'=0_\cB$ for all $b,b'\in P$. When $\Xi$ is a fam on $\cB$, we also write $\Pbf^\Xi\coloneqq \Pbf^\cB$,  
\end{notation}

In general, fams on fields of sets enjoy approximation properties to measures on finite sets. The following general result is inspired in~\cite[Lem.~1.2]{KST}.

\begin{theorem}\label{m390}
    Let $\cB$ be a Boolean subalgebra of $\pts(X)$ and let $\Xi$ be a finite fam on $\cB$ with $\delta\coloneqq \Xi(X)>0$. Then, for any $\varp>0$, any finite $F\subseteq X$ and any $P\in\Pbf^\Xi$, setting $c\coloneqq \left\lceil \frac{\delta}{\varp} \right\rceil$, there are some non-empty finite $u\subseteq X$ and some probability measure $\mu$ on $\pts(u)$ such that, for all $b\in P$:
    \begin{enumerate}[label = \rm (\alph*)]
        \item $\left|\delta\mu(u\cap b) - \Xi(b)\right| < \varp$,
        
        \item whenever $b$ is infinite, $u\cap b\cap F = \emptyset$ and $\mu(\{i\})=\frac{1}{c}$  for any $i\in u\cap b$,
        
        \item $b\cap u =\emptyset$ whenever $\Xi(b) = 0$, and 

        \item\label{m390.d} $\vert u \vert \leq c$ and the equality holds when there are no finite $b\in P$ with positive measure. In the latter case, for any $b\in P$, $\mu(u \cap b) = \frac{\vert u \cap b \vert}{\vert u \vert}$. 
    \end{enumerate}
\end{theorem}
\begin{PROOF}{\ref{m390}}
Let $\varp>0$, $F\subseteq X$ finite, and let $P=\set{b_m}{m<m^*}\in\Pbf^\Xi$. 
Find a sequence $\la r_m\colon m<m^*\ra$ of nonnegative rational numbers with denominator $c$ such that
\begin{itemize}
    \item $r_m =0$ iff $\Xi(b_m) = 0$,
    
    \item $|\Xi(b_m) - \delta r_m| < \varp$ for all $m<m^*$, and
    
    \item $\sum_{m<m^*}r_m =1$.
\end{itemize}
We can write $r_m=\frac{k_m}{c}$ for some natural number $0 \leq k_m \leq c$. The construction proceeds by recursion on $m$, also ensuring that $\left\lfloor \sum_{\ell<m}c\frac{\Xi(b_\ell)}{\delta} \right\rfloor = \sum_{\ell<m}k_\ell$ along the way, which is clear for $m=0$. Assume we have constructed $\Seq{k_\ell}{\ell<m}$. If $\Xi(b_m)=0$ just let $k_m\coloneqq 0$; otherwise, compare $\left\lfloor \sum_{\ell\leq m}c\frac{\Xi(b_\ell)}{\delta} \right\rfloor$ with $\sum_{\ell<m}k_\ell$. In case they are equal set $k_m\coloneqq \left\lfloor c\frac{\Xi(b_m)}{\delta} \right\rfloor$, else $k_m\coloneqq \left\lceil c\frac{\Xi(b_m)}{\delta} \right\rceil$. At the end of the construction, $\sum_{m<m^*}k_m = \left\lfloor \sum_{m<m^*} c \frac{\Xi(b_m)}{\delta}\right\rfloor = c$.

When $b_m$ is either infinite or of measure zero, we can find $u_m\subseteq b_m\menos F$ such that $|u_m|=k_m$; when $b_m$ is finite with positive measure, let $u_m = \{i_m\}$ for some $i_m\in b_m$. Set $u\coloneqq \bigcup_{m<m^*}u_m$, which has size $\leq c$. For $m<m^*$ and $i\in u_m$, define
\[\mu(\{i\})\coloneqq \left\{
\begin{array}{ll}
    \frac{1}{c}, & \text{when $b_m$ is infinite,}\\
    r_m, & \text{when $b_m$ is finite.}
\end{array}
\right.\]
Then, $\mu$ extends to a probability measure on $\pts(u)$ and it is as required.
\end{PROOF}

Concerning free fams, we have the following classification.

\begin{corollary}\label{m420}
    Let $\cB$ be a Boolean subalgebra of $\pts(X)$ and let $\Xi$ be a finite fam on $\cB$ with $\delta\coloneqq \Xi(X)>0$. Then the following statements are equivalent.
    \begin{enumerate}[label = \rm (\roman*)]
        \item\label{m420i} For any $\varp>0$, any finite $F\subseteq X$ and any $P\in\Pbf^\Xi$, there is some non-empty finite $u\subseteq X\menos F$ such that, for all $b\in P$,
        \[\left|\delta\frac{|u\cap b|}{|u|} - \Xi(b)\right| < \varp.\]
        \item\label{m420ii} All finite sets in $\cB$ have measure zero.
    \end{enumerate}
    Even more, the set $u$ in~\ref{m420i} can be found of size $\left\lceil \frac{\delta}{\varp} \right\rceil$ and disjoint with all $b\in P$ of measure zero.
\end{corollary}

\begin{PROOF}{\ref{m420}}
    \ref{m420i}${}\Rightarrow{}$\ref{m420ii}: Assume~\ref{m420i} and that there is some (non-empty) finite $ F\in \cB$ such that $\Xi(F)$ has positive measure. Applying~\ref{m420i} to $\varp < \Xi(F)$, to $F$ and to the partition $\{F, X\menos F\}$, we find $u\subseteq X\menos F$ such that $\left|\delta\frac{|u\cap F|}{|u|} - \Xi(F) \right| <\varp$, which implies that $\Xi(F) <\varp$, a contradiction.

    \ref{m420ii}${}\Rightarrow{}$\ref{m420i}: This implication is immediate  by \autoref{m390}~\ref{m390.d}. 
%
%
% Assume~\ref{m420ii}. Let $\varp>0$, $F\subseteq X$ finite, and let $\la B_m\colon m<m^*\ra$ be a finite partition of $X$ by sets in $\cB$. 
% Find $u\subseteq X$, $c\in\omega$ and $\mu$ as in \autoref{m390}. Then, whenever $B_m$ has measure zero, $B_m\cap u=\emptyset$, otherwise $B_m$ must be infinite, in which case $\mu(\{i\}) = \frac{1}{c}$ for all $i\in B_m\cap u$, and therefore for all $i\in u$. Hence, $|u|=c$ and $\mu(u\cap B_m) = \frac{|u\cap B_m|}{|u|}$, so~\ref{m420i} follows. \Andres{Esta implicación es inmediata usando (d) de \autoref{m390}.} \DM{Pero toca escribir para verificar que $\mu(u\cap B_m) = \frac{|u\cap B_m|}{|u|}$, aunque eso se puede incluir en \autoref{m390}.}
% Find a sequence $\la r_m\colon m<m^*\ra$ of nonnegative rational numbers such that
% \begin{itemize}
%     \item $r_m =0$ iff $\Xi(B_m) = 0$,
%     \item $|\Xi(B_m) - \delta r_m| < \varp$ for all $m<m^*$, and
%     \item $\sum_{m<m^*}r_m =1$.
% \end{itemize}
% We can write $r_m=\frac{k_m}{c}$ for some natural number $0 \leq k_m \leq c$, where $c>0$ is a natural number (not depending on $m$). Since $\Xi$ is free, any set of positive measure must be infinite, so we can find $u_m\subseteq B_m\menos F$, for all $m<m^*$, such that $|u_m|=k_m$. Then $u\coloneqq \bigcup_{m<m^*}u_m$ is as required in~\ref{m420i}.
\end{PROOF}

We also have an alternative version of \autoref{m390} that allows bounding $u$ by $|P|$.

\begin{corollary}\label{m421}
    Under the assumptions of \autoref{m390}, for any $\varp>0$ and $P\in\Pbf^\Xi$, there are some non-empty finite $u\subseteq X$ of size ${\leq}\min\left\{\left\lceil \frac{\delta}{\varp} \right\rceil, |P|\right\}$ and a probability measure $\mu$ on $\pts(u)$ such that $|\delta\mu(u\cap b)-\Xi(b)|<\varp$ for all $b\in P$, even more, $u\cap b=\emptyset$ whenever $b\in P$ has measure zero. 
\end{corollary}
\begin{PROOF}{\ref{m421}}
    In the case $\left\lceil \frac{\delta}{\varp} \right\rceil \leq |P|$ use \autoref{m390}; in the case $|P| \leq \left\lceil \frac{\delta}{\varp} \right\rceil$, pick some $u$ such that $|b\cap u| =1$ for all $b\in P$ of positive measure, and $u\cap b=\emptyset$ whenever $b\in P$ has measure zero. It is clear that $|u|\leq|P|$.

    For $i\in u$ there is a unique $b\in P$ of positive measure such that $i\in b$, so set $\mu(\{i\})\coloneqq \frac{\Xi(b)}{\delta}$. 
\end{PROOF}

For a fixed $\varp$, a finer $P\in\Pbf^\Xi$ ensures a better approximation of $\Xi$. Notice that the set $\set{|P|}{P\in\Pbf^\Xi}$ is bounded iff $\cB$ is finite, so the upper bound $|P|$ is irrelevant when we work with an infinite field of sets (as we can refine a partition to be as large as we want). On the other hand, the bound $|P|$ is relevant when $\cB$ is finite. In this case, the partition in $\Pbf^\Xi$ of maximum size is precisely the partition of $X$ into atoms of $\cB$, and this partition gives the best approximation of $\Xi$ to a finite measure. As it can be observed from the proof, the approximation is quite trivial.

We propose the following notion that results from omitting the finite set $F$ in~\ref{m420i} of the previous result. 

\begin{definition}\label{m392}
    Let $\Xi$ be a finite fam on a Boolean subalgebra $\cB$ of $\pts(X)$ and $\delta\coloneqq \Xi(X)$. We say that $\Xi$ has the \emph{uniform approximation property (uap)} if, for any $\varp>0$ and any $P\in\Pbf^\Xi$, there is some non-empty finite $u\subseteq X$ such that, for all $b\in P$,
        \[\left|\delta\frac{|u\cap b|}{|u|} - \Xi(b)\right| < \varp.\]
\end{definition}

This means that $\Xi$ can be approximated by \emph{uniform probability measures on finite sets.}

We aim to characterize the finite fams with the uap.

\begin{lemma}\label{m395}
    Let $\Xi$ be a finite fam on a Boolean subalgebra of $\pts(X)$ with $\delta\coloneqq \Xi(X)$. Assume that $P\in\Pbf^\Xi$ contains some finite $b_0$ with positive measure. Then, for small enough $\varp>0$, if $u\subseteq X$ is as in \autoref{m392} then $|u| \leq \frac{|b_0|}{\Xi(b_0)/\delta}$ and $\Xi(b) = \delta\frac{|u\cap b|}{|u|}$ for all $b\in P$.
\end{lemma}
\begin{PROOF}{\ref{m395}}
    %For any finite $F\in\cB$, define $L_F \coloneqq  \set{\frac{\ell}{d}\delta}{d>0,\ \ell\leq |F|,\ \ell,d\in\omega}\menos\{\Xi(F)\}$, which is a closed set with unique limit point $0$ when $\Xi(F)>0$.
    For any $b\in\cB$ define $L'_b \coloneqq  \set{\frac{\ell}{n}\delta}{ \ell\leq n\leq  \frac{|b_0|}{\Xi(b_0)/\delta}\text{ in } \omega }\menos\{\Xi(b)\}$, which is a closed subset of $\R$. Let $\varp >0$ smaller than the (euclidean) distance of $\Xi(b)$ with $L'_{b}$ for all $b\in P$,\footnote{We stipulate that the distance between some $x\in\R$ and $\emptyset$ is $\infty$.} and such that the integer parts of $\frac{|b_0|}{\Xi(b_0)/\delta}$ and $\frac{|b_0|}{(\Xi(b_0)-\varp)/\delta}$ are equal.

    Assume that $u\subseteq X$ is non-empty and finite such that, for all $b\in P$,
    \[\left|\delta\frac{|u\cap b|}{|u|} - \Xi(b)\right| < \varp.\]
    %
    %By applying \autoref{m392} to the partition $\set{B_m\cap B_0}{m<m^*}\cup \set{B_m\menos B_0}{m<m^*}$, there is some non-empty finite $u\subseteq X$ such that, for all $m<m^*$,
    %\[\left|\delta\frac{|u\cap B_m\cap B_0|}{|u|} - \Xi(B_m\cap B_0)\right| < \frac{\varp}{2m^*} \text{ and }\left|\delta\frac{|u\cap B_m\menos B_0|}{|u|} - \Xi(B_m\menos B_0)\right| < \frac{\varp}{2m^*}.\]
    Then
    \[%\left|\delta\frac{|u\cap B_0|}{|u|} - \Xi(B_0)\right| < \frac{\varp}{2} 
      %\text{, so } 
      |u| \leq \frac{|b_0|}{\frac{\Xi(b_0)-\varp}{\delta}}.\]
    Therefore, $|u|\leq \frac{|b_0|}{(\Xi(b_0))/\delta}$ (because the latter has the same integer part as $\frac{|b_0|}{(\Xi(b_0)-\varp)/\delta}$).
    
    On the other hand, 
    %\[\left|\delta\frac{|u\cap B_m|}{|u|} - \Xi(B_m)\right| < \frac{\varp}{m^*} \text{ for $m<m^*$,}\]
    %which implies that 
    $\Xi(b) = \delta\frac{|u\cap b|}{|u|}$ for all $m<m^*$ because the distance between $\Xi(b)$ and $L'_{b}$ is larger than $\varp$.
\end{PROOF}   

Next, we provide a characterization of the fams satisfying the uap.

\begin{theorem}\label{m400}
    Let $\Xi$ be a finite fam on a Boolean subalgebra $\cB$ of $\pts(X)$ and $\delta\coloneqq \Xi(X)$. Then the following statements are equivalent.
    \begin{enumerate}[label = \rm (\roman*)]
        \item\label{m400i} $\Xi$ has the uap.
        \item\label{m400ii} Either all finite sets in $\cB$ have measure zero, or there is some natural number $d>0$ such that, for any $b\in\cB$, $\Xi(b)= \frac{\ell}{d}\delta$ for some natural $0\leq \ell\leq d$, and $\ell\leq |b|$ when $b$ is finite.
    \end{enumerate}
    Even more, the second case of~\ref{m400ii} implies that $d\leq \frac{\delta}{\Xi(F)/|F|}$ for any finite $F\in\cB$ of positive measure.
\end{theorem}
\begin{PROOF}{\ref{m400}}
\ref{m400ii}${}\Rightarrow{}$\ref{m400i}: 
In the case when all finite sets in $\cB$ have measure zero,~\ref{m400i} follows by \autoref{m420} (or trivially when $\delta = 0$). So assume the second case of~\ref{m400ii}. Let $\cA_{0}$ be the set of minimal finite sets in $\cB$ with positive measure. Note that the members of $\cA_{0}$ are pairwise disjoint, which implies that $|\cA_{0}|\leq d$ (because $\frac{1}{d}\delta$ is a lower bound of the measure of the sets in $\cB$ with positive measure). Also, for every finite $F\in \cB$, $\Xi(F) = \sum\set{\Xi(E)}{E\subseteq F,\ E\in \cA_{0}}$. 

Let $F^*\coloneqq \bigcup \cA_{0}$, which is finite and in $\cB$. Then, $\Xi(F\menos F^*) =0$ for every finite $F\in\cB$. Let $P\in\Pbf^\Xi$.
For each $b\in P$, $b\menos F^*$ has measure $\frac{\ell_b}{d}\delta$ for some $\ell_b\leq d$, and $\ell_b>0$ implies that $b\menos F^*$ is infinite. So we can find some $v_m\subseteq b\menos F^*$ of size $\ell_b$. 
On the other hand, $\Xi(b\cap F^*)=\delta\frac{\ell'_b}{d}$ for some $\ell'_b\leq |b\cap F^*|$, hence there exists some $v'_b\subseteq b\cap F^*$ of size $\ell'_b$.
Let $u_b\coloneqq v_b\cup v'_b$ and $u\coloneqq \bigcup_{b\in P}u_b$, which is as required, even more, $\Xi(b)=\delta\frac{|u\cap b|}{|u|}$ for all $b\in P$ and $|u|=d$.

On the other hand, if $F\in\cB$ is finite with positive measure, then $\Xi(F)=\delta\frac{|u\cap F|}{d}\leq \delta\frac{|F|}{d}$, so $d\leq \frac{\delta}{\Xi(F)/|F|}$.

\ref{m400i}${}\Rightarrow{}$\ref{m400ii}: Assume~\ref{m400i} and that $\cB$ contains some finite set of positive measure. As before, let $\cA_{0}$ be the set of minimal finite sets in $\cB$ with positive measure (which is non-empty and pairwise disjoint). Note that, for any $F\in\cA_{0}$, the only subsets of $F$ in $\cB$ are $\emptyset$ and itself.
Pick one $F_0\in\cA_{0}$.

As a consequence of \autoref{m395}:
\begin{enumerate}[label = \rm (\arabic*)]
    \item $\Xi(b) \in L'\coloneqq  \set{\frac{\ell}{n}\delta}{ \ell\leq n\leq  \frac{|F_0|}{\Xi(F_0)/\delta}}$ for any $b\in \cB$.
    \item Any pairwise disjoint family of sets in $\cB$ with positive measure is finite. In particular, $\cA_{0}$ is finite.
    \item For any $b\in\cB$ with positive measure, there is some $b'\subseteq b$ in $\cB$ with minimal positive measure, i.e.\ $\Xi(b')> 0$ and any subset of $b'$ in $\cB$ has either measure zero or the same measure as $b'$.
\end{enumerate}

Also, note that the members of $\cA_{0}$ have minimal positive measure.

Let $\cC_0$ be a pairwise disjoint family of sets in $\cB$ with minimal positive measure, extending $\cA_{0}$, and maximal with this property. Then $\cC_0$ is finite and $\Xi\left(\bigcup \cC_0\right) = \Xi(X)$. Let $d\coloneqq |u|$ where $u$ results by applying \autoref{m395} to $F_0$ and $\cC_0\cup\{X\menos \bigcup \cC_0\}$. This $d$ is as required because, for any $b\in\cB$, $\Xi(b)=\sum_{c\in\cC_0}\Xi(b\cap c)$, and each $\Xi(b\cap c)$ is either $0$ or equal to $\Xi(c) = \delta\frac{|u\cap c|}{|u|}$.
\end{PROOF}

As a consequence, the size of the $u$ found from the uap can be bounded as in \autoref{m420}.

\begin{corollary}\label{m401}
    Assume that $\Xi$ has the uap and $\delta\coloneqq \Xi(X)\in(0,\infty)$. Then, for small enough $\varp>0$, the $u$ in \autoref{m392} can be found of size ${\leq}\left\lceil \frac{\delta}{\varp} \right\rceil$.
\end{corollary}
\begin{PROOF}{\ref{m401}}
    We have two cases according to \autoref{m400}. If all finite sets in $\cB$ have measure zero then the result holds by \autoref{m420} (for any $\varp>0$); the other case follows by the proof of \ref{m400ii}${}\Rightarrow{}$\ref{m400i} of \autoref{m400} by noting that $u$ can be found of size $d$, so the result is valid for any $\varp$ below $\min\largeset{\frac{\delta}{\Xi(F)/|F|}}{F\in\cB \text{ is finite with positive measure}}$.
\end{PROOF}

We give a name to the fams satisfying the second case of \autoref{m400}~\ref{m400ii}.

\begin{definition}\label{m403}
    Let $\cB$ be a Boolean subalgebra of $\pts(X)$ and $\Xi$ a finite fam on $\cB$. We say that $\Xi$ is \emph{uniformly supported} if there is some natural number $d>0$ such that, for any $b\in\cB$, $\Xi(b)= \frac{\ell}{d}\Xi(X)$ for some natural $0\leq \ell\leq d$, and $\ell\leq |b|$ when $b$ is finite.
\end{definition}

The name is motivated by the following fact for Boolean algebras in general. 
%The fams of the second case of \autoref{m400}~\ref{m400ii} can be classified as follows, even for Boolean algebras in general.

\begin{lemma}\label{m410}
    Let $\cB$ be a Boolean algebra, $\Xi$ a finite fam on $\cB$ with $\delta\coloneqq \Xi(1_\cB)$, and let $d>0$ in $\omega$. Then the following statements are equivalent:
    \begin{enumerate}[label =\rm (\roman*)]
        \item\label{410i} For any $b\in\cB$, $\Xi(b)= \frac{\ell}{d}\delta$ for some natural number $0\leq \ell\leq d$.
        \item\label{410ii} There is some $P\in\Pbf^\Xi$ such that, for any $b\in P$, $\Xi(b)= \frac{\ell_b}{d}\delta$ for some natural $0\leq \ell_b\leq d$, and any $a\leq b$ in $\cB$ with positive measure has the same measure as $b$.
        %\item\label{410iii} The quotient of $\cB$ modulo the ideal of measure zero sets, with the fam defined in a natural way from $\Xi$, is isomorphic with $\pts(N)$ for some natural number $N$ with a measure $\mu$ on $\pts(N)$ such that $\mu(\{i\}) = \delta\frac{k_i}{d}$ for some natural $0<k_i\leq d$. \Diego{maybe not too interesting}
    \end{enumerate}
\end{lemma}

Condition~\ref{410ii} means that $\Xi$ is the direct sum of fams obtained from ultrafilters on each $b\in P$, i.e.\ for each $b\in P$ there is some ultrafilter $U_b$ on $\cB |_{b}$ (see \autoref{n364}) such that $\Xi(a) = \sum_{b\in P}\frac{\ell_b}{d}\Xi_{U_b}(a\wedge b)$ for all $a\in\cB$.

\begin{PROOF}{\ref{m410}}
    \ref{410i}${}\Rightarrow{}$\ref{410ii}:  
    Define $b_m\in\cB$ by induction on $m<\omega$: if $b^0_m\coloneqq {\sim}\bigvee_{j<m}b_j$ has measure zero, we let $b_m\coloneqq b^0_m$, otherwise find the minimum natural number $\ell_m\in (0,d]$ such that $\frac{\ell_m}{d}\delta$ is the minimum measure of an $a\leq b^0_m$ in $\cB$ with positive measure, and we choose some $b_m\leq b^0_m$ in $\cB$ of measure $\frac{\ell_m}{d}\delta$. Within this induction, we can find a (minimal) $m^*$ such that $b_{m^*}=\emptyset$. Then $b_n=\emptyset$ for all $n\geq m^*$, so $P=\{ b_m\colon m<m^*\}$ is the required partition.

    %\ref{410ii}${}\Rightarrow{}$\ref{410iii}: \redq{To complete}

    %\ref{410iii}${}\Rightarrow{}$\ref{410ii}: \redq{To complete}

    \ref{410ii}${}\Rightarrow{}$\ref{410i}: Same idea of the final sentence of the proof of \autoref{m400}.
\end{PROOF}

%Condition~\ref{410ii} above indicates that $\Xi$ looks like a \emph{sum of ultrafilters from $b$ for $b\in P$} since, whenever $b$ has positive measure, the collection of of points below $b$ in $\cB$ with positive measure forms an ultrafilter on $\cB$ restricted to $b$. Conversely, if $U_b$ is an ultrafilter on $\cB$ restricted to $b$ for each $b\in P$, then the fam on $\cB$ defined by $\Xi(a)\coloneqq \delta\sum_{b\in P}\frac{\ell_b}{d}\Xi_{U_m}(a\cap b)$ for some $\ell_b\in\omega$ ($b\in P$) satisfying $\sum_{b\in P}\ell_m = d$ satisfies condition~\ref{410ii} above.

We can then define the support of a uniformly supported fam.

\begin{definition}\label{m404}
    Let $\Xi$ be a uniformly supported fam on a field of sets $\cB$ over $X$. A \emph{support of $\Xi$} is a $P\in\Pbf^\Xi$ satisfying \autoref{m410}~\ref{410ii} and containing all the minimal finite sets of positive measure.
\end{definition}

By  the proof of \autoref{m400}, any uniformly supported fam has a support.

As a consequence of \autoref{m400}, when $K$ is finite we can characterize fams with the uap using uniform measures (see \autoref{m37}). 

\begin{corollary}\label{t95}
    Let $K$ be a finite set and $\Xi$ a probability fam on $\calP(K)$. Then, $\Xi$ has the uap iff there is some (unique) $K_0\subseteq K$ such that $\Xi = \Xi^{K_0}$. Even more, such a fam is uniformly supported.
\end{corollary}
\begin{PROOF}{\ref{t95}}
    We just deal with the non-trivial implication. According to \autoref{m400}, either all finite sets in $\cB$ have measure zero, or there is some natural number $d>0$ such that, for any $B\in\cB$, $\Xi(B)= \frac{\ell}{d}$ for some natural $0\leq \ell\leq d$, and $\ell\leq |B|$ when $B$ is finite. In the first case, clearly $K_{0} \coloneqq \emptyset$ is the required set. 
    Assume now the conditions of the second case. Define $K_{0} \coloneqq \{ k \in K \colon \Xi(\{ k \}) > 0 \}$. Clearly, $\vert K_{0} \vert = d$ and for any $k \in K_{0}$, $\Xi(\{ k \}) = \frac{1}{d}$. Let $a \subseteq K$ and $m \coloneqq \vert a \cap K_{0} \vert$. Then, 
    $$
    \Xi^{K_{0}}(a) = \frac{\vert a \cap K_{0} \vert}{\vert K_{0} \vert} = \frac{m}{d} = \sum_{k \in a \cap K_{0}} \Xi(\{ k \}) = \Xi(a \cap K_{0}) = \Xi(a). 
    $$
    %Thus, $\Xi^{K_{0}} = \Xi$. 
    
    Thus, $\Xi^{K_{0}} = \Xi$. Finally, to prove the uniqueness, assume that $K_{0}^{\ast} \subseteq K$ is such that $\Xi = \Xi^{K_{0}^{\ast}}$. Clearly, for any $k \in K$, %Clearly, any element in $K_{0}^{\ast}$ has positive measure, hence $K_{0}^{\ast} \subseteq K_{0}$. On the other hand, if $k \in K_{0} \setminus K_{0}^{\ast}$ for some $k \in K$, then $0 = \Xi^{K_{0}^{\ast}}(\{  k \}) = \Xi^{K_{0}}(\{ k \}) > 0$, which is a contradiction. Thus, $K_{0} = K_{0}^{\ast}$, proving uniqueness.
    $k \in K_{0} \Leftrightarrow \Xi^{K_{0}}(\{ k \}) > 0 \Leftrightarrow \Xi^{K_{0}^{\ast}}(\{  k\}) > 0 \Leftrightarrow k \in K_{0}^{\ast}$, hence $K_{0} = K_{0}^{\ast}$. 
\end{PROOF}

Another example of uniformly supported fams are those associated with filters.

\begin{corollary}\label{m402}
    Let $\cB$ be a Boolean subalgebra of $\pts(X)$ and $F$ a filter on $\cB$. Then, $\Xi_F$ is uniformly supported. 
\end{corollary}
\begin{PROOF}{\ref{m402}}
    For $b\in\la F\ra_\cB$, $\Xi_F(b)=1$ whenever $b\in F$, and $\Xi_F(b)=0$ otherwise.
\end{PROOF}

Notice that $\{X\}$ is a support of $\Xi_F$.

\begin{example}\label{m35}
    Let $h\colon X\to Y$, $\cB$ a Boolean subalgebra of $\pts(X)$, $\Xi$ a finite fam on $\cB$ and let $\Xi_h$ be the fam on $h^\to(\cB)$ as in \autoref{m36}. 
    In the case when $h$ is one-to-one, if all finite sets in $\cB$ have $\Xi$-measure zero, then all finite sets in $h^\to(\cB)$ have $\Xi_h$-measure zero. 
    
    However, uap and ``uniformly supported" do not follow this way. For example, let $X= 10 = \{0,1,\ldots,9\}$, $\cB=\pts(X)$, $Y=\omega$, $\Xi$ is the uniform probability measure on $X$ (i.e.\ all singletons have measure $\frac{1}{10}$), and $h\colon X\to Y$ is the function that sends $0$ to $0$ and the rest to $1$. Then $\Xi_h(\{0\})=\frac{1}{10}$, $\Xi_h(\{1\})=\frac{9}{10}$ and $\Xi_h(\{n\})=0$ for $n\geq 2$. Although $\Xi$ is uniformly supported, $\Xi_h$ does not have the uap. 
    %Furthermore, if for any $x \in X$, $\{ x \} \in \cB$, and $h$ is finite-to-one, then $\Xi_{h}$ is free iff $\Xi$ is free. 
\end{example}

We conclude this subsection with the following general fact.

\begin{lemma}\label{m405}
    Let $\cC$ be a Boolean subalgebra of $\pts(X)$ and $\cB$ a Boolean subalgebra of $\cC$. Assume that $\Xi$ is a finite fam on $\cC$.
    \begin{enumerate}[label = (\alph*)]
        \item If $\Xi$ has the uap, then so does $\Xi\frestr\cB$. 
        \item If all the finite sets in $\cC$ have measure zero, then so do all finite sets in $\cB$.
        \item If $\Xi$ is uniformly supported, then so is $\Xi\frestr \cB$.
    \end{enumerate}
\end{lemma}
\begin{PROOF}{\ref{m405}}
    Easy by the definitions.
\end{PROOF}

\section{Extension criteria}\label{3}

In this section, we present several existence and extension criteria for finitely additive measures. More precisely, we consider questions of the following type: if $\cB$ is a Boolean algebra and $\Xi$ is a finitely additive measure defined on a Boolean subalgebra of $\cB$, can $\Xi$ be extended to a finitely additive measure on $\cB$? A further variant of the problem arises when $\Xi$ is replaced by a function $f$ whose domain is an arbitrary subset of $\cB$ containing the top element $1_{\cB}$ of $\cB$ (see \autoref{f700}). We also consider a family of finitely additive measures $\langle \Xi_{i} \colon i\in I \rangle$ defined on subalgebras $\langle \cB_{i} \colon i \in I \rangle$ of $\cB$, respectively, and we establish simple conditions characterizing when it is possible to find another fam $\Xi$ on the subalgebra generated by $\bigcup_{i\in I} \cB_i$ that simultaneously extends all the measures $\Xi_i$ (see \autoref{3.2}). The particular case of two families is known as the \emph{Compatibility Theorem} for fams in \cite{BhaskaraRa} (see \autoref{m48}). 

\subsection{Compactness: the main element for extension criteria}\label{3.3}

We begin by studying the compactness of the set of finitely additive measures in function spaces, which will be the fundamental point for establishing existence and extension criteria of fams. The existence criteria we present in this section are inspired from results (without proof) of~\cite[Sec.~1]{Sh00}.

\begin{theorem}\label{m40}
    Let $\cB$ be a Boolean algebra, $\delta\in[0,\infty]$ and $Z \coloneqq {}^{\cB} [0,\delta]$ with the product topology,\footnote{Where $[0, \delta]$ has the usual topology inherited from $[- \infty, \, \infty].$} which is compact.\footnote{By \emph{Tychonoff's Theorem}.} Then, the set of finitely additive measures on $\cB$ with $\Xi(1_{\cB}) \leq \delta$ is closed in $Z.$ As a consequence, it is compact in $Z.$ 

    The same is obtained if we demand $\Xi(1_\cB)=\delta$ (instead of just $\leq$).
\end{theorem}

\begin{PROOF}{\ref{m40}}
    Notice that the set $\{z\in Z \colon  z(0_\cB)=0 \}$ is closed in $Z$ because we can write it as the following product of closed sets: $$\{0\} \times \left( \prod_{a\in\cB \setminus \{0_\cB \}}[0,\delta] \right) .$$ 
    
    So, it remains to show that $\bigcap\{F_{a,b} \colon a,b\in\cB,\ a\wedge b=0_\cB\}$ is closed where, for any $a, b \in \cB$ with $a \wedge b = 0_{\cB}, \, F_{a,b}\coloneqq \{z\in Z \colon  z(a \vee b)=z(a)+z(b)\}.$ Let $a, b \in \cB$ such that $a \wedge b = 0_{\cB}.$ It is enough to show that $F\coloneqq F_{a,b}$ is closed. For this,  suppose that $z\in \mathrm{cl}(F)$ and 
    distinguish two cases: 

\noindent\textbf{Case} $z(a \vee b)<\infty.$ In this case we must have $z(a), z(b) < \infty$. If this is not true,  without loss of generality we can assume that $z(a) = \infty.$ Consider the open set $U'\coloneqq \prod_{x\in \cB}U'_x$ where:

        \begin{itemize}
            \item $U'_{a\vee b}\coloneqq  \left( z(a \vee b)-\frac{1}{2}, \, z(a \vee b)+\frac{1}{2} \right) \cap[0,\infty),$

            \item $U'_a\coloneqq \left( z(a \vee b)+\frac{1}{2},\infty \right],$

            \item $U'_x\coloneqq [0,\infty],$ if $x  \in \cB \setminus \{ a, a \vee b \}.$
        \end{itemize}
        
        Since $U'$ is an open neighborhood of $z$, we can find a $z'\in U'\cap F_{a,b}.$ So, by the definition of $U',$ we have that $$z'(a \vee b)=z'(a)+z'(b), \ |z'(a \vee b)-z(a \vee b)|<\frac{1}{2} \ \text{ and } \ z'(a)>z(a \vee b)+\frac{1}{2}.$$ However, $z'(a) \leq z'(a)+z'(b)=z'(a \vee b)<z(a \vee b)+\frac{1}{2}$, which is a contradiction. Thus, $z(a), z(b) < \infty.$
    
        Now, let $\varepsilon>0$ and consider the open neighborhood of $z$ defined by $U\coloneqq \prod_{x \in \cB}U_x$, where:

        \begin{itemize}
            \item $U_x\coloneqq (z(x)-\varepsilon,z(x)+\varepsilon)\cap[0,\delta]$ when $x\in\{a,b,a\vee b\}$,

            \item  $U_x\coloneqq [0,\delta]$ for any $x \in \cB \setminus \{ a, b, a \vee b \}$.
        \end{itemize}
        
        So there is some $z''\in U\cap F$. 
        Therefore $|z(x) -z''(x)|<\varepsilon$ for $x\in\{a,b,a\vee b\}$ and $z''(a \vee b)=z''(a)+z''(b)$. Hence
        \[|z(a \vee b)-z(a)-z(b)|=|z(a \vee b)-z''(a \vee b)+z''(a)-z(a)+z''(b)-z(b)|<3\varepsilon.\]
        Since $\varepsilon$ is arbitrary, $z(a \vee b)-z(a)-z(b)=0$, so $z\in F$.

\noindent\textbf{Case} $z(a \vee b)=\infty$ In this case we must have that either $z(a)=\infty$ or $z(b)=\infty$. Assume not, so $z(a),z(b)<\infty$ and consider $U'\coloneqq \prod_{x\in \cB}U'_x$ where,
        
        \begin{itemize}
            \item $U'_a\coloneqq \left( z(a)-\frac{1}{2},z(a)+\frac{1}{2} \right) \cap[0,\infty),$

            \item $U'_b\coloneqq \left( z(b)-\frac{1}{2},z(b)+\frac{1}{2} \right) \cap[0,\infty),$

            \item $U'_{a\vee b}\coloneqq (z(a)+z(b)+1,\infty],$

            \item $U'_x\coloneqq [0,\infty],$ when $x \in \cB \setminus \{a, b, a \vee b\}.$
        \end{itemize}
        
        It is clear that $U'$ is an open neighborhood of $z$ and, therefore, we can find some $z'\in F\cap U'$, so $|z'(x)-z(x)|<\frac{1}{2}$ for $x\in\{a,b\}$ and $z'(a)+z'(b)=z'(a \vee b)>z(a)+z(b)+1$. However, $z'(a)+z'(b)<z(a)+z(b)+1$, which is a contradiction. 
    
        Therefore $z(a)+z(b)=\infty=z(a \vee b)$, so $z\in F$.  Thus, $F$ is closed. 
\end{PROOF}

In the following result, we translate the compactness of finitely additive measures proved in \autoref{m40}, to a local property of finite subsets. This is a general criterion to ensure the existence of finitely additive measures.

\begin{theorem}\label{m44}
    Let $\cB$ be a Boolean algebra, $\langle b_i \colon  i \in I\rangle\subseteq\cB$, $0\leq\delta\leq\infty$, $K$ a closed subset of $[0,\delta]$,  and let $\langle K_i \colon  i\in I\rangle$ be a collection of closed subsets of $K$. 
    Assume that $\calX$ is a closed subset of ${}^{\cB'}K$  where $\cB'$ is the Boolean subalgebra generated by $\{b_i\colon\, i\in I\}$. 
    Then the following statements are equivalent.
    
    \begin{enumerate}[label=\normalfont(\roman*)]
        \item\label{m44i} There is a finitely additive measure 
        $\Xi \in \calX$ on $\cB'$ %the subalgebra generated by $\{b_i \colon i \in I\}$ 
        such that $\ran(\Xi) \subseteq K$ and $\Xi(b_i) \in K_i$ for all $i \in I$.
        
        \item\label{m44ii} For any finite $J \subseteq I$ and any collection $\langle G_i\colon i\in J\rangle$ of open subsets of $K$ such that $K_i\subseteq G_i$ for all $i\in J$, there is some
        $\Xi \in \calX$ such that $\Xi\frestr \la \{b_i \colon  i \in J\}\ra$ is a fam, 
        $\ran(\Xi) \subseteq K$ and $\Xi(b_i)\in G_i$ for all $i\in J$.
    \end{enumerate}
\end{theorem}

\begin{PROOF}{\ref{m44}}
    \index{finite intersection property}
    
    The implication~\ref{m44i} $\Rightarrow$~\ref{m44ii} is immediate. To prove the converse, %let $\cB'$ be the Boolean subalgebra generated by $\{b_i \colon i \in I \}$. 
    first fix $J\subseteq I$ finite. Since each $K_{i}$ is closed, for each $i\in J$ we can find a decreasing sequence $\langle F^i_k \colon  k<\omega\rangle$ of closed subsets of $K$ such that there is an open subset $G^i_k$ of $K$ such that $K_i\subseteq G^i_k\subseteq F^i_k$, and $\bigcap_{k<\omega} F^i_k = K_i$. Let $\cB_J$ be the subalgebra generated by $\{b_i\colon i\in J\}$. For any $k<\omega$ consider 
    \[C_{J,k} \coloneqq \{z\in \calX\cap{}^\cB K \colon  z  \, {\rest} \, \cB_J \text{ is a fam and }\forall i\in J\ (z(b_i)\in F^i_k) \}.\]

    By \autoref{m40}, $C_{J,k}$ is a compact set and, by~\ref{m44ii}, it is a non-empty set. Moreover, since $\langle C_{J,k} \colon  k<\omega\rangle$ is decreasing, it has the \emph{finite intersection property}, so 
    \[C_J\coloneqq \bigcap_{k<\omega} C_{J,k}=\{z\in \calX\cap{}^\cB K \colon z \, {\rest} \, \cB_J \text{ is a fam and }\forall i\in J\ (z(b_i)\in K_i) \}\] is a non-empty set. Also, $J\subseteq J'$ implies $C_{J'}\subseteq C_J$, so $\langle C_J \colon J\in[I]^{<\aleph_0}\rangle$ has the finite intersection property, thus it has non-empty intersection. If $\Xi$ is in this non-empty intersection, then it is as required.
\end{PROOF}

More concrete results can be obtained for Boolean subalgebras of $\pts(X)$.

\begin{lemma}\label{uapcl}
    Let $\cB$ be a Boolean subalgebra of $\pts(X)$ and $\delta\in [0,\infty)$. Then the set of fams $\Xi$ on $\cB$ with the uap and such that $\Xi(X)=\delta$ is closed in ${}^\cB [0,\delta]$.
\end{lemma}
\begin{PROOF}{\ref{uapcl}}
   Denote by $\calX$ the set of fams $\Xi$ on $\cB$ with the uap and such that $\Xi(X)=\delta$. If $\Xi$ is in the closure of $\calX$ then it is a fam on $\cB$ and $\Xi(X)=\delta$ by \autoref{m40}. We now show that $\Xi$ has the uap. For any $P\in\Pbf^\cB$ and $\varp>0$, we can find some $\Xi'\in \calX$ such that $|\Xi(b) - \Xi'(b)|<\frac{\varp}{2}$ for all $b\in P$. Then, there is some non-empty finite $u\subseteq X$ such that $\left| \delta\frac{|u\cap b|}{|u|} - \Xi'(b) \right| < \frac{\varp}{2}$ for all $b\in P$, so $\left| \delta\frac{|u\cap b|}{|u|} - \Xi(b) \right| < \varp$.
\end{PROOF}

It is clear that the lemma above still holds when considering fams such that all finite setsin $\cB$ have measure zero. However, this is not true for uniformly supported fams because, as the proof itself shows, any $\Xi$ with the uap gets approximated by a uniformly supported fam (namely, $\Xi^u$). 

The following existence criterion relates to \autoref{m390} and \autoref{m421}. 

\begin{theorem}\label{famgen-}
Let $X$ and $I$ be sets, $\delta\in(0,\infty)$, $\la a_i\colon i\in I\ra\subseteq\pts(X)$ and let $\la K_i\colon i\in I\ra$ be a collection of closed subsets of $[0,\delta]$. Then, the following statements are equivalent.
   \begin{enumerate}[label = \rm (\roman*)]
       \item\label{eqfamgen-} For every finite $J\subseteq I$ and any open set $G_i\subseteq\R$ containing $K_i$ for $i\in J$, there is some non-empty finite $u\subseteq X$ and some probability measure $\mu$ on $\pts(u)$ such that, for every $i\in J$, $\delta\mu(a_i\cap u)\in G_i$.
       \item\label{famgen2-} There is a fam $\Xi$ on the Boolean algebra generated by $\set{a_i}{i\in I}$ such that $\Xi(X)=\delta$ and $\Xi(a_i)\in K_i$ for all $i\in I$.
   \end{enumerate} 
   Even more, the set $u$ in~\ref{eqfamgen-} can be found of size ${\leq} \min\left\{2^{|J|},\left\lceil \frac{\delta}{\varp}\right\rceil\right\}$.
\end{theorem}
\begin{PROOF}{\ref{famgen-}}
%Wlog, we assume that $X=a_{i_0}$ for some $i_0\in I$ and $K_{i_0}=\{\delta\}$. 
\ref{eqfamgen-}${}\Rightarrow{}$\ref{famgen2-}: 
By \autoref{m44} and \autoref{uapcl} we need to show that, for any finite $J\subseteq I$ and $\la G_i\colon i\in J\ra$ as in~\ref{eqfamgen-}, there is some fam $\mu^*$ on $\pts(X)$ such that $\mu^*(a_i)\in G_i$ for all $i\in J$ and $\mu^*(X)=\delta$. Choose $u$ as in~\ref{eqfamgen-} and define $\mu^*(b)\coloneqq \delta\, \mu(u\cap b)$ for all $b\in\la\set{a_i}{i\in J}\ra$. 
%as in \autoref{m37}. 
Then $\mu^*$ is as required. %(since it satisfies \autoref{m400}~\ref{m400ii}).

\ref{famgen2-}${}\Rightarrow{}$\ref{eqfamgen-}: 
Assume that $\Xi$ is a fam as in~\ref{famgen2-}. Let $J\subseteq I$ be finite and let $G_i\subseteq \R$ be an open set containing $K_i$, for any $i\in J$. 
Let $a^0_i\coloneqq a_i$ and $a^1_i\coloneqq X\menos a_i$, and for any $q\in 2^J$ let $a^q\coloneqq \bigcap_{i\in J}a_i^{q(i)}$. Note that $\{a^q\colon q\in 2^J\}$ is a finite partition of $X$. %Let $E\coloneqq \{q\in 2^J\colon \Xi(a^q)>0\}$, so $X'\coloneqq \bigcup_{q\in E}a^q$ has measure $\delta$. 
Let $\varp>0$ such that $(\Xi(a_i)-\varp,\Xi(a_i)+\varp)\subseteq G_i$ for all $i\in J$.
Then, by \autoref{m421}, there is some non-empty finite $u\subseteq X$ of size ${\leq}\min\left\{2^{|J|},\left\lceil \frac{\delta}{\varp}\right\rceil\right\}$ and a probability measure $\mu$ on $\pts(u)$ such that \[\left|\delta\, \mu(u\cap a^q) - \Xi(a^q)\right| < \frac{\varp}{2^{|J|}},\]
for all $q\in 2^J$, which implies, for any $i\in J$,
\[\left|\delta\, \mu(u\cap a_i) - \Xi(a_i)\right| < \frac{\varp}{2},\]
so $\delta\,\mu(u\cap a_i)\in G_i$.
\end{PROOF}

Thanks to \autoref{m401}, we can reformulate the previous results for fams with the uap.

\begin{theorem}\label{famgen}
%Let $X$ and $I$ be sets, $\delta\in[0,\infty)$, $\la a_i\colon i\in I\ra\subseteq\pts(X)$ and let $\la K_i\colon i\in I\ra$ be a collection of closed subsets of $[0,\delta]$. Then, 
In the context of \autoref{famgen-}, 
the following statements are equivalent.
   \begin{enumerate}[label = \rm (\roman*)]
       \item\label{eqfamgen} For every finite $J\subseteq I$ and any open set $G_i\subseteq\R$ containing $K_i$ for $i\in J$, there is some non-empty finite $u\subseteq X$ such that, for every $i\in J$, $\delta\frac{|a_i\cap u|}{|u|}\in G_i$.
       \item\label{famgen2} There is a fam $\Xi$ on the Boolean algebra generated by $\set{a_i}{i\in I}$ such that $\Xi(X)=\delta$, $\Xi$ has the uap and $\Xi(a_i)\in K_i$ for all $i\in I$.
   \end{enumerate} 
   Even more, the set $u$ in~\ref{eqfamgen-} can be found of size ${\leq} \left\lceil \frac{\delta}{\varp}\right\rceil$ for small enough $\varp$.
\end{theorem}
\begin{PROOF}{\ref{famgen}}
    Proceed exactly as in the proof of \autoref{famgen-} with the following changes: in the proof of \ref{eqfamgen}${}\Rightarrow{}$\ref{famgen2}, set $\mu^*\coloneqq \delta\Xi^u$
    (see \autoref{m37}), which has the uap as a consequence of \autoref{m400}~\ref{m400ii}; in the proof of the converse, use \autoref{m401} instead of \autoref{m421}.
\end{PROOF}

In \autoref{famgen}, we cannot show \ref{famgen2}${}\Rightarrow{}$\ref{eqfamgen} without the demand that $\Xi$ has the uap. For example, let $\Xi$ be the probability ($\sigma$-additive) measure on $\pts(\omega)$ such that
\[\Xi(\{k\}) = 
   \left\{
      \begin{array}{ll}
         \frac{1}{3}  &  k=0,\\[1ex]
         \frac{2}{3}  &  k=1,\\[1ex]
         0  & k>1. 
      \end{array}
   \right.\]
Let $I\coloneqq \{0,1\}$, $a_i\coloneqq  \{i\}$, and $K_i\coloneqq \{\Xi(a_i)\}$. Then~\ref{eqfamgen} fails for $J=I$ and for $G_i\coloneqq (a_i-\varp,a_i+\varp)$ when $0<\varp<\frac{1}{6}$ (here, $\delta=1$).

On the other hand, concerning \ref{eqfamgen}${}\Rightarrow{}$\ref{famgen2}, we may not be able to find a $\Xi$ such that $\Xi(b) = 0$ for all finite $b$, even when $K_i = \{0\}$ for any finite $a_i$. For example, in the case when $X = \omega$, $\delta=1$, $I=\{0,1\}$, $a_0$ and $a_1$ are infinite disjoint sets such that $|\omega\menos(a_0\cup a_1)|=2$, and $K_0=K_1=\{\frac{1}{4}\}$.

To obtain fams where all finite sets have measure zero, we modify~\ref{eqfamgen} as follows.

\begin{theorem}\label{FAMchar}
   In the context of \autoref{famgen-}, the following are equivalent:
   \begin{enumerate}[label=\normalfont(\roman*)]
       \item\label{Fc1} For any finite $v\subseteq X$, any finite $J\subseteq I$ and any open set $G_i\subseteq\R$ containing $K_i$ for $i\in J$, there is some non-empty finite $u\subseteq X\smallsetminus v$ such that, for every $i\in J$, $\delta\frac{|a_i\cap u|}{|u|}\in G_i$.

       \item\label{Fc2} There exists a fam $\Xi$ on the Boolean algebra generated by $\set{a_i}{i\in I}$ such that $\Xi(X)=\delta$, all finite sets in the Boolean algebra have measure zero, and 
       $\Xi(a_i)\in K_i$ for every $i\in I$.
   \end{enumerate}
    Even more, the set $u$ in~\ref{eqfamgen-} can be found of size ${\leq} \left\lceil \frac{\delta}{\varp}\right\rceil$.
\end{theorem}
\begin{PROOF}{\ref{FAMchar}}
    To prove \ref{Fc1}${}\Rightarrow{}$\ref{Fc2},
    we first show that $0\in K_i$ for any $i\in I$ such that $a_i$ is finite. Towards a contradiction assume that $0\notin K_i$. We apply~\ref{Fc1} to $J=\{i\}$, $v=a_i$ and some open set $G_i$ in $\R$ such that $K_i\subseteq G_i$ and $0\notin G_i$, to get some $u\subseteq X\smallsetminus a_i$ such that $\delta\frac{|a_i\cap u|}{|u|}\in G_i$. %Assume that $a_i\cap u\neq\emptyset$, pick some $n\in a_i\cap u$. Then $n<n$ because $m$ is a upper bounds of $a_i$ and $u\subseteq\omega\smallsetminus m$. 
    Since $a_i\cap u=\emptyset$, $0\in G_i$ follows, which is a contraction.

    Let $\cB$ be the Boolean algebra generated by $\set{a_i}{i\in I}$ and 
    let $\{a_i\colon i\in I'\}$ be an indexation of all the finite sets in $\cB$ that are not in $\{a_i\colon i\in I\}$, where $I'\cap I=\emptyset$. Set $I_1\coloneqq I\cup I'$, $K'_i\coloneqq \{0\}$ for all $i\in I_1$ such that $a_i$ is finite, else $K'_i \coloneqq  K_i$. To prove~\ref{Fc2}, it is enough to show that $\la a_i\colon i\in I_1\ra$ and $\la K'_i\colon i\in I_1\ra$ satisfy~\ref{eqfamgen} of \autoref{famgen}.
    Let $J\subseteq I_1$ be finite and, for $i\in J$ let $G'_i\subseteq\R$ open containing $K'_i$. For $i\in I$ with  $a_i$ finite, since $K'_i \subseteq K_i$, extend $G'_i$ to some open set $G_i$ containing $K_i$, and set $G_i\coloneqq G'_i$ otherwise. 
    Let $v\coloneqq \bigcup \set{a_i\in[X]^{<\aleph_0}}{i\in J}$, which is finite. Then, by~\ref{Fc1}, there is some non-empty finite $u \subseteq X\menos v$ such that, for any $i\in J$, $\delta\frac{|a_i\cap u|}{|u|}\in G_i$. When $a_i$ is finite, $a_i \subseteq v$ and $a_i\cap u =\emptyset$, thus $\delta\frac{|u\cap a_i|}{|u|} = 0 \in G'_i$.
    
    For the converse, proceed as in the proof of \ref{famgen2-}${}\Rightarrow{}$\ref{eqfamgen-} of \autoref{famgen-}, but find $u \subseteq X\menos v$ by application of \autoref{m420}.
    %let $\Xi$ be a free fam as in~\ref{Fc2}. Under the assumptions of~\ref{Fc2}, for any  $i\in J$ let $a^0_i\coloneqq a_i\menos v$ and $a^1_i\coloneqq X\menos (a_i\cup v)$, and for any $q\in 2^J$ let $a^q\coloneqq \bigcap_{i\in J}a_i^{q(i)}$. Note that $\{a^q\colon q\in 2^J\}$ is a finite partition of $X\menos v$. Let $E\coloneqq \{q\in 2^J\colon \Xi(a^q)>0\}$, so $X'\coloneqq \bigcup_{q\in E}a^q$ has measure $\delta$. 
    %
    %Note that $\set{\bigcup_{q\in F}a^q}{F\subseteq E}$ is the Boolean subalgebra of $\pts(X')$ generated by $\{a_i\cap X'\colon i\in J\}$, which is finite, and that $a^q\cap X'$ is infinite for any $q\in E$ (because $\Xi$ is a free fam). For each $q\in E$, choose a rational number $r^q\in(0,1)$ such that $\delta\cdot r^q$ is very close to $\Xi(a^q)$, and $\sum_{q\in E}r^q=1$. Concretely, $\delta\cdot r^q$ is close enough to $\Xi(a^q)$ to guarantee that, for any $i\in J$,
    %\[\delta \sum\set{r_q}{q\in E,\ q(i)=0}\in G_i \text{ (i.e.\ this sum is very close to $\Xi(a_i)$).}\]
    %Find $0<d<\omega$ such that, for all $q\in E$, $r^q\coloneqq \frac{c^q}{d}$ for some natural number $0<c^q\leq d$, and choose some $u^q\subseteq a^q$ of size $c^q$. Then $u \coloneqq\bigcup_{q\in E}u^q$ is as required.
\end{PROOF}

\subsection{Compatibility}\label{3.2}

In this subsection, we are going to use the compactness results proved in \autoref{3.3} to establish extension criteria. First, we show when and how a function can be extended to a fam, which implies a general criterion to amalgamate many fams simultaneously. Afterward, we focus on the amalgamation of two fams, and also three under some restrictions.

We use the following connection between finite fams on a finitely generated Boolean algebra and positive linear maps on finite dimensional $\bbQ$-vector spaces.

\begin{lemma}\label{f700-}
    Let $\cB_0$ be a finite Boolean algebra and $A\coloneqq \At_{\cB_0}$ (its set of atoms). Then there is a canonical bijection between the finite fams on $\cB_0$ and the real-valued positive linear maps on the $\bbQ$-vector space $\bbQ^N$ where $N\coloneqq |A|$.
\end{lemma}
\begin{PROOF}{\ref{f700-}}
    Enumerate $A=\set{a_i}{i<N}$. 
    A finite fam $\Xi$ on $\cB_0$ determines the linear map  $L_\Xi\colon\bbQ^N\to\R$ such that $L_\Xi(e^i)=\Xi(a_i)$ for all $i<N$, where $\set{e^i}{i<N}$ denotes the standard basis of $\bbQ^N$. If $x\in\bbQ^N$ and $x\geq 0$ (i.e.\ $x_i\geq0$ for all $i<N$) then $L_\Xi(x)\geq 0$, that is, the linear map is \emph{positive}.

    We show that $\Xi\mapsto L_\Xi$ is the desired bijection. To define the inverse, consider a positive linear map $L\colon\bbQ^N\to\R$. For any $a\in\cB_0$ there is some unique $F_a\subseteq N$ such that $a=\bigvee_{i\in F_a}a_i$, so define $\Xi_L\colon \cB_0\to \R$ by $\Xi_L(a)\coloneqq  L\left(\sum_{i\in F_a} e^i\right)$. Using that $L$ is positive, it is easy to show that $\Xi_L$ is a finite fam. 

    Moreover, $\Xi_{L_\Xi} = \Xi$ and $L_{\Xi_L} = L$, which completes the proof.
\end{PROOF}

% Notice that, for a natural number $0<N<\omega$, there exists a bijection between finite fams on $\pts(N)$ and real-valued positive linear maps on $\Q^N$ (as a $\Q$-vector space). Indeed, if $\Xi$ is a finite fam on $\pts(N)$, 
% %where $N < \omega$ and we consider $\Q^N$ as a $\Q$-vector space, 
% then $\Xi$ defines the linear map $L_\Xi\colon\Q^N\to\R$ such that $L_\Xi(e^i)=\Xi(\{i\})$ for all $i<N$, where $e^i=\chi_{\{i\}}$ and $\chi_a$ denotes the characteristic function of $a\subseteq N$ (over $N$). Since $\la e^i\colon i<N\ra$ is a basis of the vector space $\Q^N$, we get $L_\Xi(\chi_a)=\Xi(a)$ for all $a\subseteq N$. Moreover, if $x\in\Q^N$ and $x\geq 0$ (i.e.\ $x_i\geq0$ for all $i<N$) then $L_\Xi(x)\geq 0$, that is, the linear map is \emph{positive}.
% Conversely, if $L\colon\Q^N\to\R$ is a positive linear map, then $\Xi_L\colon \pts(N)\to\R$, defined by $\Xi_L(a)\coloneqq L(\chi_a)$, is a finite fam on $\pts(N)$. 
% %Therefore, there is a bijection between finite fams on $\pts(N)$ and real-valued positive linear maps on $\Q^N$.

In the previous proof, notice that $a\mapsto F_a$ is an isomorphism from $\cB_0$ onto $\pts(N)$.

In the context of the previous result, for $F\subseteq\cB_0$ denote by $C_F$ the subspace of $\bbQ^N$ generated by $\set{v_a}{ a\in F}$ where $v_a\coloneqq \sum_{i\in F_a} e^i$. Notice that $v_a$ can be seen as the characteristic function of $F_a$ over $N$. 

We now look at the following problem. 
Given a function $f\colon F\to[0,\infty)$, we ask whether there is some positive linear map $L_f\colon C_F\to\R$ such that $L_f(v_a)=f(a)$. It is clear that there is at most one such map. 
%
%We now look at the following problem. For $A\subseteq\cB_0$ denote by $C_A$ the subspace of $\Q^N$ generated by $\{\chi_a\colon a\in A\}$. Given a function $f\colon A\to[0,\infty)$, we ask whether there is some positive linear map $L_f\colon C_A\to\R$ such that $L_f(\chi_a)=f(a)$. It is clear that there is at most one such map.
%
To approach the problem, we use the following fact. 
Recall from \autoref{b41} that $a_\sigma$ for $\sigma\in S\coloneqq \set{\sigma\colon F\to\{0,1\}}{a_\sigma\neq 0_{\cB_0}}$ forms the set of atoms of the subalgebra generated by $F$.

\begin{lemma}\label{preLf}
    In the context of \autoref{f700-}, let $F\subseteq\cB_0$, $S$ as above, and $g\colon F\to \R$. Then:
    \begin{enumerate}[label = \normalfont (\alph*)]
        \item\label{preLf-a} For $a\in F$, $\displaystyle v_a = \sum_{\substack{\sigma\in S\\ \sigma(a) = 0}} v_{a_\sigma}$.
        \item\label{preLf-b} $\displaystyle \sum_{a\in F}g(a)v_a = \sum_{\sigma \in S}\Bigg(\sum_{\substack{a\in F\\ \sigma(a) = 0}}g(a)\Bigg)v_{a_\sigma}$.
        \item\label{preLf-c} $\displaystyle \sum_{a\in F}g(a)v_a \geq 0$ iff $\displaystyle \sum_{\substack{a\in F\\ \sigma(a) = 0}}g(a) \geq 0$ for all $\sigma \in S$.
    \end{enumerate}
\end{lemma}
\begin{PROOF}{\ref{preLf}}
    \ref{preLf-a} follows because $\set{a_\sigma}{\sigma\in S,\ \sigma(a) = 0}$ partitions $a$ for all $a\in F$; for~\ref{preLf-b}:
    \[\sum_{a\in F}g(a) v_a = 
    \sum_{a\in F}g(a)\sum_{\substack{\sigma\in S\\ \sigma(a)=0}} v_{a_\sigma} =
    \sum_{\sigma \in S}\Bigg(\sum_{\substack{a\in F\\ \sigma(a) = 0}}g(a)\Bigg)v_{a_\sigma};\]
    and~\ref{preLf-c} follows by~\ref{preLf-b} and because $\set{v_{a_\sigma}}{\sigma\in S}$ form vectors with entries $0$ and $1$, where $1$ does not appear at the same position in two different vectors.
\end{PROOF}

\begin{lemma}\label{clm:exLf}
In the context of \autoref{f700-}, let $F\subseteq\cB_0$, $f\colon F\to[0,\infty)$, and $S$ as above. Then the following statements are equivalent.
\begin{enumerate}[label = \normalfont (\roman*)]
    \item\label{Lf1} $L_f$ exists.
    \item\label{Lf2} For any $h\colon F\to\Z$: if 
    %all $\sigma\colon F\to \{0,1\}$ such that $a_\sigma=\bigwedge_{a\in F}a^{\sigma(a)} \neq 0_\cB$ satisfies that 
    $\displaystyle \sum_{\substack{a\in F\\ \sigma(a)=0}}h(a)\geq 0$ for all $\sigma\in S$, then $\displaystyle \sum_{a\in F}h(a) f(a) \geq 0$.
\end{enumerate}
When $\cB_0$ is a field of sets over $X$, the following statement is equivalent with those above.
\begin{enumerate}[resume*]
    \item\label{Lf3} For any $h\colon F\to\Z$, if $\sum_{a\in F}h(a)\chi_{a}\geq0$, then $\sum_{a\in F}h(a)f(a)\geq0$.
\end{enumerate}
\end{lemma}
\begin{PROOF}{\ref{clm:exLf}}
    In the case that $\cB_0$ is a field of sets,
    \[\sum_{a\in F}h(a)\chi_a = \sum_{a\in F} h(a) \sum_{\substack{\sigma\in S\\ \sigma(a)=0}}\chi_{a_\sigma} = \sum_{\sigma\in S} \sum_{\substack{a\in F\\ \sigma(a) = 0}} h(a) \chi_{a_\sigma} = \sum_{\sigma\in S} \Bigg(\sum_{\substack{a\in F\\ \sigma(a) = 0}} h(a)\Bigg) \chi_{a_\sigma}.\]
    For any $a\in F$, 
    Since $\set{a_\sigma}{\sigma\in S}$ forms a partition of $X$, we get
    \[\sum_{a\in F}h(a)\chi_a\geq 0 \text{ iff } \sum_{\sigma\in S} \Bigg(\sum_{\substack{a\in F\\ \sigma(a) = 0}} h(a)\Bigg) \chi_{a_\sigma} \geq 0 \text{ iff } \sum_{\substack{a\in F\\ \sigma(a) = 0}} h(a) \geq 0 \text{\ for all $\sigma\in S$.}\]
    This shows~\ref{Lf2}${}\Leftrightarrow{}$\ref{Lf3}. 

    \ref{Lf1}${}\Rightarrow{}$\ref{Lf2}: Let $h\colon F\to\Z$ and assume that $\displaystyle \sum_{\substack{a\in F\\ \sigma(a)=0}}h(a)\geq 0$ for all $\sigma\in S$.
    Set $v\coloneqq \sum_{a\in F}h(a)v_a$. By \autoref{preLf}~\ref{preLf-c}, we get that $v\geq 0$, so $0\leq L_f(v) = \sum_{a\in F}h(a)f(a)$ by~\ref{Lf1}.

    \ref{Lf2}${}\Rightarrow{}$\ref{Lf1}: We show that, for any $g\colon F\to\bbQ$, if $\sum_{a\in F}g(a)v_a\geq0$ then $\sum_{a\in F}g(a)f(a)\geq0$. This implies that the definition $L_f(x)\coloneqq \sum_{a\in F}q_af(a)$ whenever $x = \sum_{a\in F}q_av_a$ for rationals $q_a$ ($a\in F$) is well defined (i.e.\ it does not depend on the linear combination of $x$ in terms of $v_a$ for $a\in F$) and that $L_f$ is a positive linear map.

    So assume $\sum_{a\in F}g(a) v_{a}\geq0$. We can find $0<d<\omega$ and $h\colon F\to \Z$ such that $g(a)=\frac{h(a)}{d}$ for any $a \in F$, so $\sum_{a\in F}h(a) v_{a}\geq0$. By \autoref{preLf}~\ref{preLf-c}, $\displaystyle \sum_{\substack{a\in F\\ \sigma(a) = 0}}h(a) \geq 0$ for all $\sigma\in S$.     
    Therefore, by~\ref{Lf2}, $\sum_{a\in F}h(a)f(a)\geq0$, which implies $\sum_{a\in F}g(a)f(a)\geq0$.
\end{PROOF}

Below, we establish a general criteria to extend arbitrary functions to fams, which results from \autoref{clm:exLf} and compactness (see \autoref{3.3}). 

\begin{theorem}[cf.~{\cite[Thm.~3.2.10]{BhaskaraRa}}]\label{f700}
    Let $\cB$ be a Boolean algebra and $f$ a partial function from $\cB$ into $[0,\infty)$ with $1_\cB\in \dom f$. Then, the following statements are equivalent.
    \begin{enumerate}[label = \normalfont (\roman*)]
        \item\label{f700.1} There is a finite fam $\Xi$ on the subalgebra generated by $\dom f$ extending $f$.
        \item\label{f700.2} For any finite $F\subseteq\dom f$ and $h\colon F\to \Z$:  if all $\sigma\colon F\to \{0,1\}$ such that $a_\sigma=\bigwedge_{a\in F}a^{\sigma(a)} \neq 0_\cB$ satisfies that $\displaystyle \sum_{\substack{a\in F\\ \sigma(a)=0}}h(a)\geq 0$, then $\displaystyle \sum_{a\in F}h(a) f(a) \geq 0$.
        %\item For any $\Seq{a_i}{i<m}, \Seq{b_j}{j<n}\in {}^{<\omega}\dom f$ (allowing repetitions): if all $\sigma\colon m\to \{0,1\}$ and $\tau\colon n\to\{0,1\}$ such that $\bigwedge_{i<m}a_i^{\sigma(i)}\wedge\bigwedge_{j<n}b_j^{\tau(j)}\neq 0_\cB$ satisfies that %$\sum_{i<m}(1-\sigma(i)) \leq \sum_{j<n}(1-\tau(j))$
        %$\big(\sum_{j<n}\tau(j)\big)-\big(\sum_{i<m}\sigma(i)\big) \leq n-m$, then $\sum_{i<m}f(a_i) \leq \sum_{j<n} f(b_j)$. 
    \end{enumerate}
    When $\cB$ is a field of sets over $X$, the following statement is equivalent with those above.
    \begin{enumerate}[resume*]
        \item\label{f700.3} For any finite $F\subseteq \dom f$ and $h\colon F\to\Z$, \[\sum_{a\in F}h(a)\chi_a\geq 0 \text{ implies } \sum_{a\in F}h(a) f(a)\geq 0.\]
        %\item For any $\Seq{a_i}{i<m}, \Seq{b_j}{j<n}\in {}^{<\omega}\dom f$, 
        %\[\sum_{i<m}\chi_{a_i} \leq \sum_{j<n}\chi_{b_j} \text{ implies } \sum_{i<m}f(a_i) \leq \sum_{j<n}f(b_j).\]
    \end{enumerate}
    %Moreover, in the case that $f$ allows the value $\infty$, then (iii) (and (v) in the case of field of sets) is equivalent to (i) but without the requirement that the fam is finite. 
\end{theorem}

\begin{PROOF}{\ref{f700}}
    The proof of~\ref{f700.2}${}\Leftrightarrow{}$\ref{f700.3} is exactly the same as~\ref{Lf2}${}\Leftrightarrow{}$\ref{Lf3} of \autoref{clm:exLf}. 

    Assume~\ref{f700.1} and let $F\subseteq \dom f$ finite and $h\colon F\to \Z$. Let $\cB_0$ be the subalgebra of $\cB$ generated by $F$ (which is finite) and let $\Xi_0 \coloneqq  \Xi\frestr \cB_0$, where $\Xi$ is as in~\ref{f700.1}. By \autoref{f700-}, we have a positive linear map $L_{\Xi_0}$, which shows the existence of $L_{f\frestr F}$, so~\ref{f700.2} follows by \autoref{clm:exLf}. This shows~\ref{f700.1}${}\Rightarrow{}$\ref{f700.2}. 

    Now assume~\ref{f700.2}. Thanks to \autoref{m44} with $\delta = f(1_\cB)$, $K=[0,\delta]$ and each $K_i$ a singleton, it is enough to prove~\ref{f700.1} when $\cB$ is finite. Let $A=\set{a_i}{i<N}$ be the set of atoms of $\cB$. By \autoref{clm:exLf} applied to $F=\dom f$ (which is finite now), we have that $L_f$ exists. The domain of $L_f$ is $C_F$. the subspace generated by $\set{v_a}{a\in F}$, and $L_f(a)=f(a)$ for $a\in F$. To find $\Xi$ as in~\ref{f700.1}, by \autoref{f700-} it is enough to extend $L_f$ to a positive linear map with domain $\bbQ^N$. This is a consequence of Hahn-Banach's Theorem, but we do not need to use AC in our application because we only need to extend the map on a base of $\bbQ^N$, which is finite. For illustration, we show that, whenever $W$ is a subspace of $\bbQ^N$, $L\colon W\to \R$ is a positive linear map and $x\notin W$, we can extend $L$ to a positive linear map on the subspace $W'$ generated by $W\cup\{x\}$. Every member of $W'$ is expressed, uniquely, as $w+qx$ for $w\in W$ and $q\in\bbQ$, so it is enough to define $x\mapsto z$ for some $z\in\R$ to extend $L$ to a linear map on $W'$. However, $z$ must be properly chosen so that the extension is also positive. In case $L'$ is such an extension, then $L(w)\leq L'(x)\leq L(v)$ for any $w,v\in W$ such that $w\leq x\leq v$. Note that
    \[\sup\{L(w)\colon w\in W,\ w\leq x\}\leq\inf\{L(v)\colon v\in W,\ x\leq v\},\]
    so we can choose any $z$ between these numbers and define $L'(x)\coloneqq z$. This gives us a (unique) linear map extending $L$ on $W'$. To show that it is positive, assume that $w'\in W'$ and $w'\geq 0$. So $w'=w+qx$ for some $w\in W$ and $q\in\bbQ$. If $q=0$ then clearly $w'=w\geq 0$, so $L'(w')=L(w)\geq0$; if $q<0$ then $-\frac{1}{q}w\geq x$, so $-\frac{1}{q}L(w)\geq z=L'(x)$ by the definition of $z$, which implies $L'(w')\geq 0$; and the case $q>0$ is proved similarly.
\end{PROOF}

As a consequence, we can deduce a general criterion to extend many fams simultaneously.

\begin{corollary}\label{f701}
    Let $\cB$ be a Boolean algebra and let $\set{\Xi_i}{i\in I}$  be a collection of finite fams on subalgebras of $\cB$. Then, the following statements are equivalent.
    \begin{enumerate}[label = \normalfont (\roman*)]
        \item\label{f701-1} There is some fam $\Xi$, on a subalgebra of $\cB$, extending $\Xi_i$ for $i\in I$.
        \item\label{f701-2} Any two fams in $\set{\Xi_i}{i\in I}$ are compatible as functions and $f\coloneqq \bigcup_{i\in I}\Xi_i$ satisfies~\ref{f700.2} of \autoref{f700}.
    \end{enumerate}
\end{corollary}

%We are ready to prove the \emph{Compatibility Theorem of finitely additive measures}. 
In the case of two fams, we have a better \emph{Compatibility Theorem}. Although this is presented in \cite[Ch.~3]{BhaskaraRa}, in \autoref{m48} below we offer not only a more direct proof but also a generalization to a broader context, where it is formulated for Boolean algebras rather than for fields of sets. 

\begin{theorem}\label{m48}
    Let $\cB$ be a Boolean algebra and, for $d\in\{0,1\}$, let $\cB_d$ be a Boolean Boolean subalgebra of $\cB$ with a finitely additive measure $\Xi_d\colon\cB_d\to[0,\infty)$. Then the following statements are equivalent.
    
    \begin{enumerate}[label=\normalfont(\roman*)]
        \item\label{m48-1} There is a finitely additive measure $\Xi$ on the Boolean subalgebra generated by $\cB_0\cup\cB_1$ extending $\Xi_d$ for $d\in\{0,1\}$. 

        \item\label{m48-2} $\Xi_0(1_\cB)=\Xi_1(1_\cB)$ and, for any $a\in\cB_0$ and $a'\in\cB_{1}$, if $a\leq a'$ then $\Xi_0(a)\leq\Xi_1(a')$.      

        \item\label{m48-3} For any $d,d'\in\{0,1\}$, $a\in\cB_d$ and $a'\in\cB_{d'}$, if $a\leq a'$ then $\Xi_d(a)\leq\Xi_{d'}(a')$.
    \end{enumerate}
\end{theorem}

%% There are two implications that do not require too much work: the implication (1) $\Rightarrow$ (2) is immediate and, to prove (2) $\Rightarrow$ (3), note that (3) is clear whenever $d=d'$ or $d=0$ and $d'=1$, so we need to prove it when $d=1$ and $d'=0$. Denote $\delta\coloneqq \Xi_0(1_\cB)=\Xi_1(1_\cB)$. If $a\in\cB_1$, $a'\in\cB_0$ and $a\leq a'$, then ${\sim}a'\leq{\sim}a$, so by (2) we obtain $\Xi_0({\sim}a')\leq\Xi_1({\sim}a)$, that is, \ $\delta-\Xi_0(a')\leq\delta-\Xi_1(a)$, so $\Xi_1(a)\leq\Xi_0(a')$, which proofs (3). However, (3) $\Rightarrow$ (1)  requires notions of linear algebra and  tools of functional analysis\footnote{The connection with linear algebra is given because there is a bijection between finite finitely additive measures on $\pts(N)$ and real-valued positive linear maps on $\bbQ^N$, for any $0<  N < \omega.$}, and is a fairly extensive proof. 

\begin{PROOF}{\ref{m48}}
The implication \ref{m48-1}${}\Rightarrow{}$\ref{m48-2} is obvious. To show \ref{m48-2}${}\Rightarrow{}$\ref{m48-3}, note that~\ref{m48-3} is obvious whenever $d=d'$ or $d=0$ and $d'=1$, so we need to prove it when $d=1$ and $d'=0$. Denote $\delta\coloneqq \Xi_0(1_\cB)=\Xi_1(1_\cB)$. If $a\in\cB_1$, $a'\in\cB_0$ and $a\leq a'$, then ${\sim}a'\leq{\sim}a$, so by~\ref{m48-2} we obtain $\Xi_0({\sim}a')\leq\Xi_1({\sim}a)$, i.e.\ $\delta-\Xi_0(a')\leq\delta-\Xi_1(a)$, so $\Xi_1(a)\leq\Xi_0(a')$.

We now prove \ref{m48-3}${}\Rightarrow{}$\ref{m48-1}. Note that~\ref{m48-3} implies that $\Xi_0$ and $\Xi_1$ coincide in $\cB_0\cap\cB_1$, which is a subalgebra of $\cB$. Set $f\coloneqq \Xi_0\cup \Xi_1$, which is a real-valued function with domain $\cB_0\cup\cB_1$. By \autoref{f701}, it is enough to show that $f$ satisfies~\ref{f700.2} of \autoref{f700}. So let $F\subseteq \cB_0\cup\cB_1$ finite, $h\colon F\to\Z$, and assume that $\displaystyle \sum_{\substack{a\in F\\ \sigma(a)=0}}h(a)\geq 0$ for all $\sigma\colon F\to \{0,1\}$ such that $a_\sigma=\bigwedge_{a\in F}a^{\sigma(a)} \neq 0_\cB$. By \autoref{preLf}, this assumption is equivalent to $\sum_{a\in F}h(a)v_a\geq 0$, where: 
\begin{itemize}
    \item $\cB^-$ is the subalgebra of $\cB$ generated by $F$ and $\At_{\cB^-}=\set{a_i}{i<N}$ ($N<\omega$);
    \item for $a\in \cB^-$, $F_a\coloneqq \set{i<N}{a_i\leq a}$; 
    \item $\set{e^i}{i<N}$ is the standard basis of the $\bbQ$-vector space $\bbQ^N$; and
    \item $v_a\coloneqq \sum_{i\in F_a}e^i$ for all $a\in \cB^-$.
\end{itemize}
Recall that $a\mapsto F_a$ is a isomorphism from $\cB^-$ onto $\pts(N)$, and that $v_a$ can be seen as the characteristic function of $F_a$ over $N$.

For $e\in\{0,1\}$, let $\cB_e^-$ be the subalgebra of $\cB_e$ generated by $A\cap\cB_e$, $A_0\coloneqq \At_{\cB^-_0}$ and $A_1\coloneqq \At_{\cB^-_1}\menos A_0$. Define $h_e\colon A_e\to\Z$ by
$h_e(b)\coloneqq \sum_{\substack{a\in F\\ b\leq a}} h(a)$. Hence,
\begin{align}
    \sum_{a\in F}h(a)v_a & = \sum_{a\in F}\sum_{e=0}^1\sum_{\substack{b\in A_e\\ b\leq a}} h(a)v_b = \sum_{e=0}^1\sum_{b\in A_e} \Bigg(\sum_{\substack{a\in F\\ b\leq a}} h(a)\Bigg) v_b = \sum_{e=0}^1\sum_{b\in A_e} h_e(b)v_b \nonumber\\ 
     & = \sum_{a\in A_0} h_0(a)v_a + \sum_{b\in A_1} h_1(b)v_b. \label{eqm48-1}
\end{align}
Thus, in $\bbQ^N$,
\begin{equation}\label{eqm48-1.1}
    \sum_{a\in A_0} h_0(a)v_a + \sum_{b\in A_1} h_1(b)v_b \geq 0
\end{equation}
To finish the proof, it is enough to show that
\begin{equation}\label{eqm48-2}
    \sum_{a\in A_0} h_0(a)\Xi_0(a) + \sum_{b\in A_1} h_1(b)\Xi_1(b) \geq 0,
\end{equation}
as this implies $\sum_{a\in F}h(a)f(a)$ by using similar calculations as in~\eqref{eqm48-1}.

Define $h'_e(a)\coloneqq |h_e(a)|$, so, by moving terms with negative coefficients to the right side in~\eqref{eqm48-1.1}, we can partition $A_e=A_e^+\cup A_e^-$ such that 
\[\sum_{a\in A^+_0}h'_0(a)v_a+\sum_{b\in A^+_1}h'_1(b)v_b\geq\sum_{a\in A^-_0}h'_0(a)v_a+\sum_{b\in A^-_1}h'_1(b)v_b.\]
This implies that
\[\sum_{a\in A^+_0}h'_0(a)v_a \geq\sum_{b\in A^-_1}h'_1(b)v_b\quad \text{ and }\quad
\sum_{b\in A^+_1}h'_1(b)v_b\geq\sum_{a\in A^-_0}h'_0(a)v_a.\]
We show the first one (the second is analogous). If $i\in F_b$ for some $b\in A_1^-$ then\footnote{$(u,v)$ denotes the canonical inner product in $\bbQ^N$.}
\[\sum_{b\in A^-_1}h'_1(b)(v_b,e^i)\leq \sum_{a\in A^+_0}h'_0(a) (v_a,e^i) +\sum_{b\in A^+_1}h'_1(b) (v_b,e^i)=\sum_{a\in A^+_0}h'_0(a)(v_a,e^i)\] 
because $b$ is disjoint with any member of $A^+_1$; else, if $i\notin \bigcup_{b\in A^-_1}F_b$ then the inequality is clear.

Order $A^+_0=\{a_j\colon j<m\}$ such that the map $j\mapsto h'_0(a_j)$ is monotone increasing, and order $A^-_1=\{b_j\colon j<n\}$ in a similar way with respect to $h'_1$. Set $a'_j\coloneqq \bigvee_{\ell=j}^{m-1}a_\ell$ and $b'_j\coloneqq \bigvee_{\ell=j}^{n-1}b_\ell$ for any $j<\omega$. %(note that $a'_j=b'_{j'}=\emptyset$ whenever $j\geq m$ and $j'\geq n$). 
For $k<\omega$ define $a^*_0\coloneqq 1_\cB$, $a^*_k\coloneqq a'_j$ when $h'_0(a_{j-1})<k\leq h'_0(a_j)$ (set $h'_0(a_{-1})=h'_1(b_{-1})\coloneqq 0$), and $a^*_k\coloneqq 0_\cB$ when $k\geq\sup_{j<m}\{h'_0(a_j)+1\}$. The motivation of this definition is as follows. Consider the finite sequence $\la c_\ell\colon \ell<m'\ra$ obtained by repeating $a_j$ $h'_0(a_j)$-many times, following the order of the $j<m$. Then, for $i<N$, $i\in F_{a^*_k}$ iff $i$ is in at least $k$-many of the $F_{c_\ell}$, i.e.\ $\sum_{\ell<m'} (v_{c_\ell},e^i) \geq k$. This allows to calculate \[\sum_{a\in A^+_0}h'_0(a)\Xi_0(a)=\sum_{\ell<m'}\Xi_0(c_\ell)=\sum_{k=1}^\infty\Xi_0(a^*_k).\] 
The first equality is clear. For the second, let $c^0_\ell\coloneqq c_\ell$ and $c^1_\ell\coloneqq {\sim}c_\ell$ for each $\ell<m'$, and for $p\in {}^{m'} 2$ let $c^p\coloneqq \bigwedge_{j<\ell}c^{p(\ell)}_\ell$. Note that $a^*_k=\bigvee\{c^p\colon |p^{-1}[\{1\}]|\geq k\}$. Thus
\[\begin{split}
  \sum_{k\geq 1}\Xi_0(a^*_k) &=\sum_{k\geq 1} \sum_{\substack{p\in {}^{m'}2\\ |p^{-1}[\{1\}]|\geq k}}\Xi_0(c^p)
    = \sum_{k\geq 1}\sum_{k'\geq k} \sum_{\substack{p\in {}^{m'}2\\ |p^{-1}[\{1\}]|= k'}}\Xi_0(c^p)
    = \sum_{k'\geq 1}\sum_{k=1}^{k'} \sum_{\substack{p\in {}^{m'}2\\ |p^{-1}[\{1\}]|= k'}}\Xi_0(c^p)\\
 &  = \sum_{k'\geq 1} k'\cdot\sum_{\substack{p\in {}^{m'}2\\ |p^{-1}[\{1\}]|= k'}}\Xi_0(c^p)
    = \sum_{k'\geq 1} \sum_{\substack{p\in {}^{m'}2\\ |p^{-1}[\{1\}]|= k'}}k'\Xi_0(c^p)
    = \sum_{p\in {}^{m'}2}|p^{-1}[\{1\}]|\Xi_0(c^p)\\
 &  = \sum_{p\in {}^{m'}2}\sum_{\substack{\ell<m'\\ p(\ell)=1}}\Xi_0(c^p)
    = \sum_{\ell<m'}\sum_{\substack{p\in {}^{m'}2\\ p(\ell)=1}}\Xi_0(c^p)
    = \sum_{\ell<m'}\Xi_0(c_\ell).
\end{split}\]

Define $b^*_k$ similarly with respect to $h'_1$, so we also have $\sum_{b\in A^-_1}h'_1(b)\Xi_1(b)=\sum_{k=1}^\infty\Xi_1(b^*_k)$.
Note that $a^*_k\in\cB_0$ and $b^*_k\in\cB_1$. We show that $b^*_k\leq a^*_k$. The cases $k=0$ and $k\geq \sup_{j<n}\{h'_1(b_j)+1\}$ are clear, so assume $h'_1(b_{j-1})<k\leq h'_1(b_j)$ for some (unique) $j<n$. Hence $b^*_k=b'_j$. If $i\in F_{b^*_k}$ then $i\in F_{b_\ell}$ for some $j\leq\ell<n$, so
\[k\leq h'_1(b_\ell)=\sum_{b\in A^-_1}h'_1(b)(v_b,e^i)\leq \sum_{a\in A^+_0}h'_0(a)(v_a,e^i)=h'_0(a_{j'})\]
where $j'<m$ is unique such that $i\in F_{a_{j'}}$ (it exists because $k>0$). This implies that $i\in F_{a'_{j''}}$ where $j''<m$ is the smallest such that $k\leq h'_0(a_{j''})$ (so $j''\leq j'$). Thus $i\in F_{a^*_k}$. This shows $F_{b^*_k}\subseteq F_{a^*_k}$, which is equivalent to $b^*_k\leq a^*_k$.

Therefore, by the hypothesis~\ref{m48-3}, $\Xi_1(b^*_k)\leq\Xi_0(a^*_k)$ for all $k<\omega$, so
\[\sum_{a\in A^+_0}h'_0(a)\Xi_0(a)=\sum_{k=1}^\infty\Xi_0(a^*_k) \geq 
\sum_{k=1}^\infty\Xi_1(b^*_k)=\sum_{b\in A^-_1}h'_1(b)\Xi_1(b).\]

In the same way, we can show $\sum_{b\in A^+_1}h'_1(b)\Xi_1(b) \geq \sum_{a\in A^-_0}h'_0(a)\Xi_0(a)$. Therefore, we can conclude that~\eqref{eqm48-2} holds.\qedhere

\end{PROOF}

Immediately, we obtain a well-know result (see e.g \cite{Bachman80} and \cite{Extension}). 

\begin{corollary}\label{m52}
    Let $\cB$ be a Boolean algebra, $\cC \subseteq\cB$ a sub-Boolean with a finitely additive measure $\Xi\colon \cC \to [0,\infty)$, and let $b\in\cB$. If $z\in[0,\infty)$ is between $\sup\{\Xi(a) \colon  a\leq b,\ a\in \cC \}$ and $\inf\{\Xi(a) \colon  b\leq a,\ a \in \cC\},$ then there is a finitely additive measure $\Xi'$ on the Boolean algebra generated by $\cC \cup\{b\}$, extending $\Xi$, such that $\Xi(b)=z$.
\end{corollary}

\begin{PROOF}{\ref{m52}}
    Note that $\{b\}$ generates the Boolean subalgebra $\cC' \coloneqq \{0_\cB,b,{\sim}b,1_\cB\}$. By the hypothesis, $z\leq\Xi(1_\cB)$, so we can define the finitely additive measure $\Xi' \colon \cC' \to [0,\infty)$ such that $\Xi'(b) \coloneqq z$ and $\Xi'(1_\cB) \coloneqq \Xi(1_\cB)$. The result follows by \autoref{m48}~\ref{m48-2} and the hypothesis on $z$.
\end{PROOF}

As a direct consequence of Zorn's Lemma, we obtain:

\begin{corollary}\label{famext}
Let $\cB$ be a Boolean algebra and let $\cC\subseteq\cB$ be a subalgebra with a fam $\Xi\colon\cC\to[0,\infty)$. Then there is a fam on $\cB$ that extends $\Xi$.
\end{corollary}

As a consequence, we can strengthen \autoref{m48} by allowing in~\ref{m48-1} that $\Xi$ is a fam on the whole $\cB$. We can strengthen \autoref{m44} in the same way when $\delta<\infty$ and $K=[0,\delta]$, but if we want to allow $\delta\leq\infty$ and $K$ an arbitrary closed subset of $[0,\delta]$, we need to do a bit more of work.

\begin{lemma}\label{ex:fam+1}
   Let $0\leq\delta\leq\infty$ and let $K$ be a closed subset of $[0,\delta]$, $\cB$ a Boolean algebra and let $\cC\subseteq\cB$ be a subalgebra with a fam $\Xi\colon\cC\to K$. If $b\in\cB$ then there is a fam $\Xi'$ on the Boolean algebra generated by $\cC\cup\{b\}$, extending $\Xi$, such that $\ran\Xi'\subseteq K$.
\end{lemma}
\begin{PROOF}{\ref{ex:fam+1}}
 By \autoref{m44} it is enough to prove this when $\cB$ is finite and $\cC\cup\{b\}$ generates $\cB$. Let $A$ be the set of atoms of $\cC$. Note that they may not be atoms in $\cB$, actually, all the atoms of $\cB$ are precisely the non-zero objects of the form $a\wedge b$ or $a\wedge{\sim}b$ for $a\in A$, and all these are pairwise incompatible. 
 In the case that $a\wedge b$ and $a\wedge{\sim}b$ are non-zero, define $\Xi'(a\wedge b)\coloneqq \Xi(a)$ and $\Xi'(a\wedge{\sim}b)\coloneqq 0$, but if one of them is zero, then we set the $\Xi'$-measure of the other object as $\Xi(a)$.
 Since any function on the atoms extends to a (unique) fam, we obtain a fam $\Xi'$ into $\ran\Xi\subseteq K$ extending $\Xi$.
\end{PROOF}

As a consequence of \autoref{ex:fam+1}:

\begin{corollary}\label{famgen+}
\autoref{m44} is valid when, in~\ref{m44i}, $\Xi$ is a fam into $K$ on the whole $\cB$.
\end{corollary}

\begin{corollary}\label{extfamK}
If $\cB$ is a Boolean algebra, then any fam on a subalgebra into $K$ can be extended to a fam on $\cB$ into $K$.
\end{corollary}

Since any filter is represented by a fam onto $\{0,1\}$, \autoref{extfamK} applied to $K=\{0,1\}$ implies the well-known result that any filter on a Boolean algebra $\cB$ can be extended to an ultrafilter.

\autoref{m48} is valid when we demand that the fams are into $K=\{0,\delta\}$ for $0\leq\delta<\infty$. The case $\Xi_0(1_{\cB})=0$ is trivial; when $\Xi_0(1_{\cB})>0$, \autoref{m48}~\ref{m48-2} is equivalent to the finite intersection property of the corresponding filters. 

\begin{remark}
    In \autoref{m48}, when $K\subseteq[0,\infty)$ is compact and $\ran\Xi_e\subseteq K$ for $e\in\{0,1\}$, the conclusion is not valid when we add in~\ref{m48-1} that $\ran\Xi\subseteq K$. 
    
    Indeed, consider for instance, $X = \{ 0, 1\} \times \{ 0, 1 \}$ and the algebras generated by the vertical and horizontal sections of $X$, respectively, that is, $\cB_0\coloneqq \{\emptyset,X_0,X_1,X\}$ and $\cB_1\coloneqq \{\emptyset,X^0,X^1,X\}$, where $X_e\coloneqq \{e\}\times\{0,1\}$ and $X^e\coloneqq \{0,1\}\times\{e\}$. Consider the probability fams $\Xi_0$ and $\Xi_1$ on $\cB_0$ and $\cB_1$, respectively, determined by $\Xi_0(X_0)=\frac{1}{3}$ and $\Xi_1(X^0)\coloneqq \frac{1}{2}$. Clearly, $\ran\Xi_0 \cup \ran \Xi_1\subseteq K\coloneqq  \left\{0,\frac{1}{3},\frac{1}{2},\frac{2}{3},1\right\}$.
    
    Any fam $\Xi$ on $\pts(X)$ extending both $\Xi_0$ and $\Xi_1$ must satisfy:
    \begin{align*}
        \Xi(\{0,0\}) & =\Xi(\{1,1\})-\frac{1}{6}, & \Xi(\{(1,0)\}) & = \frac{2}{3}-\Xi(\{1,1\}), & \Xi(\{(1,0)\}) & = \frac{1}{2}-\Xi(\{1,1\}),\\[1ex]
        \frac{1}{6} & \leq \Xi(\{1,1\}) \leq \frac{1}{2}. & & & & 
    \end{align*}
    It is easy to see that $\ran \Xi \nsubseteq K$ for any such $\Xi$.
\end{remark}

We apply \autoref{m48} to obtain another extension criteria.

%\red{The following is one of the main extension results of this work:}

\begin{theorem}\label{m69}
   Let $\Xi_0$ be a finitely additive measure on a Boolean subalgebra $\cB_0$ of $\cB$ and $\{ b_i \colon i\in I\} \subseteq\cB$. Assume that $0 < \delta\coloneqq \Xi_0(1_\cB)< \infty$ and, for every finite $J\subseteq I$ and $b \in \cB_0$, if $\Xi_0(b)>0$ then $b \wedge\bigwedge_{i\in J} b_i\neq0_\cB$. Then, there exists a unique finitely additive measure $\Xi$ on $\la\cB_0\cup\set{b_i}{i\in I}\ra$ extending $\Xi_0$ such that $\Xi(b_i)=\delta$ for every $i\in I$.
\end{theorem}

\begin{PROOF}{\ref{m69}}
    By the hypothesis, and using that  $\Xi_0(1_\cB) > 0,$ we have that $\{b_i \colon  i\in I\}$ generates a filter $F$ on $\cB.$ Let %$\cB_0\coloneqq \dom(\Xi_0)$, 
    $\cB_1$ be the Boolean subalgebra generated by $F$ and define  $\Xi_{1} \coloneqq \delta \, \Xi_{F}$. %where $\Xi_{F}$ is as in \autoref{m70}. 
    So $\Xi_{1} \colon \cB_{1} \to \{0, \delta\}$ and, for any $b \in \cB_{1}$, $\Xi_{1}(b) = \delta$ iff $b \in F$. %since $b \in F \Leftrightarrow \Xi_{F}(b) = 1 \Leftrightarrow \Xi_{1}(b) = \delta.$  
    
    Now, it is enough to show that $\Xi_0$ and $\Xi_1$ satisfy \autoref{m48}~\ref{m48-2}. So let $a\in\cB_0$, $b\in\cB_1$ and assume $a\leq b$. If $b\in F$ then %, by \autoref{m15}, 
    $\Xi_0(a)\leq\delta=\Xi_1(b)$; otherwise ${\sim}b \in F$ and $\Xi_1(b)=0$, so we must show that $\Xi_0(a)=0$. If this is not the case and $\Xi_{0}(a) > 0,$ then by hypothesis, we get that $a\wedge{\sim}b\neq0_\cB$, but since $a\leq b$ it follows that $b \wedge{\sim}b\neq 0_\cB$, a contradiction. 

    To prove uniqueness, we show how $\Xi$ is uniquely determined by $\Xi_0$ and $\{ b_i\colon i\in I\}$. If $a\in \la\cB_0\cup\set{b_i}{i\in I}\ra$ then there are some $J\subseteq I$ finite and $g\colon {}^J2\to \cB_0$ such that $a=\bigvee_{\sigma \in {}^J2} g(\sigma)\wedge \bigwedge_{i\in J}b^{\sigma(i)}_i$, where $b^e_i$ is defined as in \autoref{b41}. If $\sigma(i)=0$ for all $i\in J$, then $\Xi\left(g(\sigma)\wedge \bigwedge_{i\in J}b^{\sigma(i)}_i\right) =\Xi_0(g(\sigma))$ (because $\Xi_1\left(\bigwedge_{i\in \dom\sigma}b_i\right)=\delta$); otherwise, if $\sigma(i')=1$ for some $i'\in J$, then $\Xi_1(b_{i'}^{\sigma(i')})= \Xi_1({\sim}b_{i'})=0$, so $\Xi\left(g(\sigma)\wedge \bigwedge_{i\in J}b^{\sigma(i)}_i\right)=0$. In conclusion, $\Xi(a) = \Xi_0(g(\sigma^0))$, where $\sigma^0$ is the constant map onto $0$. 
    %This shows that $\Xi$ is uniquely determined. 
\end{PROOF}

We wonder whether the results of this section are valid when restricting to fams with the uap. These are not in general. For example, let $a\subseteq\omega$ be the set of even numbers, $\cB\coloneqq \pts(\omega)$, let $\cB_0$ be the field of sets over $\omega$ generated by $a$ and let $F\subseteq\pts(\omega)$ be the filter of subsets of $\omega$ containing $\{0,1\}$. Let $\Xi_0$ be the (unique) probability measure on $\cB$ such that $\Xi(a)=\frac{2}{3}$ (which is clearly uniformly supported). By \autoref{m69}, there is a unique probability fam $\Xi$ on $\la \cB_0 \cup F \ra$ extending $\Xi_0$ such that $\Xi(b)=1$ for all $b\in F$, in particular, $\Xi(\{0,1\})=1$. Then, $\Xi(\{0\})=\Xi(a\cap \{0,1\}) = \Xi_0(a)=\frac{2}{3}$ and, naturally, $\Xi(\{1\})=\frac{1}{3}$. However, $\Xi$ cannot have the uap because no singleton can have measure $\frac{2}{3}$ under a probability fam with the uap (by \autoref{m400}). This also shows that \autoref{m48} is not valid in general when restricting to the uap: both $\Xi_0$ and $\Xi_F$ are uniformly supported, but any common extension extends $\Xi$ and thus cannot have the uap.

%As a consequence of \autoref{m400}, when $\Xi_{0}$ in \autoref{m69} has the uap, we can extend it to a fam that still has the uap.  

We can validate \autoref{m69} to fams with the uap under further conditions.

\begin{theorem}\label{t96}
    In \autoref{m69}, assume that $\cB=\pts(X)$ and that $\cB_0$ contains all the finite subsets of $X$. If $\Xi_0$ has the uap, then so does $\Xi$. Likewise for ``uniformly supported".
\end{theorem}
\begin{PROOF}{\ref{t96}}
    If all finite sets have $\Xi_0$-measure zero, then the same is true for $\Xi$. So assume that $\Xi_0$ is uniformly supported with support $P$. Since $\cB_0$ contains all finite sets, all finite sets in $P$ with positive measure must be singletons (by minimality). 
    It is not hard to show that $P$ is a support of $\Xi$, so $\Xi$ is uniformly supported.
\end{PROOF}

Similarly, we can validate the extension results when adding one set to the subalgebra.

\begin{lemma}\label{t91}
    In the context of \autoref{ex:fam+1}, assume that $\cB=\pts(X)$ and that $\cC$ contains all the singletons. Then, if $\Xi$ has the uap, then $\Xi'$ can be found with the uap. As a consequence, $\Xi$ can be extended to a fam into $K$ with the uap on the whole $\pts(X)$.

    Here, ``uap" can be replaced by ``uniformly supported", even more, extensions of $\Xi$ can be found with the same support.
\end{lemma}
\begin{PROOF}{\ref{t91}}
    The result is trivial when all singletons have $\Xi$-measure zero. So assume that $\Xi$ is uniformly supported. The argument in the proof of \autoref{ex:fam+1}, when replacing $A$ by a support of $\Xi$, allows to find a uniformly supported $\Xi'$ with the same support as $\Xi$. The rest follows by Zorn's Lemma.
\end{PROOF}

We ask whether~\ref{m48-3}${}\Rightarrow{}$\ref{m48-1} of \autoref{m48} is valid when $\Xi_0$ and $\Xi_1$ have the uap for finding $\Xi$ with the uap, under the condition that both $\cB_0$ and $\cB_1$ contain all singletons. If, say, $\Xi_0$ is free then $\Xi_1$ and $\Xi$ must be free, but we do not know how to obtain an uniformly supported $\Xi$ when  $\Xi_0$ and $\Xi_1$ are uniformly supported. 

\begin{question}\label{Qusup}
    Assume that $\cB_0$ and $\cB_1$ are fields of sets over $X$, containing all the singletons, and $\Xi_0$ and $\Xi_1$ are uniformly supported fams with the uap on $\cB_0$ and $\cB_1$, respectively. Assume that $\Xi_0$ and $\Xi_1$ are compatible (in the sense of \autoref{m48}~\ref{m48-3}). Is there a uniformly supported fam that extends both $\Xi_0$ and $\Xi_1$?
\end{question}

A strategy to solve this problem may suggest a short proof of \autoref{m48}. 

We finish this section with one more extension theorem. \autoref{m48} indicates simple conditions to extend two fams to a single one. But, are there conditions simpler than \autoref{f701} to extend three or more fams? The case of filters is easy thanks to the finite intersection property, but it is unknown how does it work with fams in general. Below, we prove a concrete case of extending three fams when one of them is associated with a filter.

\begin{theorem}\label{m80}
    For $e\in\{0,1\}$, let $\Xi_e$ be a finitely additive measure on a Boolean subalgebra $\cB_e$ of $\cB$, and let $\{  b_i \colon i\in I\}\subseteq\cB$. Assume that $\Xi_0(1_\cB)= \Xi_1(1_\cB)=\delta\in(0,\infty)$. Then, the following statements are equivalent:
    \begin{enumerate}[label=\normalfont(\roman*)]
        \item\label{m80i} There is a finitely additive measure $\Xi$ extending $\Xi_0\cup\Xi_1$ such that $\Xi(b_i)=\delta$ for all $i\in I$.

        \item\label{m80ii} The following conditions are satisfied:
          \begin{enumerate}[label = $(\bullet_{\arabic*})$]
              \item\label{m80b1} For any $e\in\{0,1\}$, $b\in\cB_e$ and $J\subseteq I$ finite, if $\Xi_e(b)>0$ then $b\wedge\bigwedge_{i\in J}b_i\neq 0_\cB$.
              \item\label{m80b2} For any $a\in\cB_0$, $b\in \cB_1$ and $J\subseteq I$ finite, if $\Xi_0(a)>0$ and $a\wedge\bigwedge_{i\in J}b_i \leq b\wedge\bigwedge_{i\in J}b_i$ then $\Xi_0(a)\leq \Xi_1(b)$.
          \end{enumerate}
    \end{enumerate}
\end{theorem}

\begin{PROOF}{\ref{m80}}
    We prove \ref{m80ii}${}\Rightarrow{}$\ref{m80i} (the converse is clear). Fix $e\in\{0,1\}$ and let $\cB'_e\coloneqq  \la\cB_e\cup\set{b_i}{i\in I}\ra$. By~\ref{m80b1}, there is a finitely additive measure $\Xi'_e$ on $\cB'_e$ extending $\Xi_e$ such that $\Xi'_e(b_i)=\delta$, uniquely determined as in the end of the proof of \autoref{m69}. To prove the theorem, is is enough to show that $\Xi'_0$ and $\Xi'_1$ satisfy \autoref{m48}~\ref{m48-2}. It is clear that $\Xi'_0(1_\cB) =\Xi'_1(1_\cB) = \delta$, so it remains to show that $\Xi'_0(a_0) \leq \Xi'_1(a_1)$ whenever $a_0\in\cB'_0$, $a_1\in\cB'_1$ and $a_0\leq a_1$. There are a finite $J\subseteq I$ and functions $g_e\colon {}^J2\to\cB'_e$ for $e\in\{0,1\}$ such that $a_e=\bigvee_{\sigma\in {}^J 2}g_e(\sigma)\wedge\bigwedge_{i\in J}b_i^{\sigma(i)}$. Recall from the proof of \autoref{m69} that $\Xi'_e(a_e)=\Xi_e(g_e(\sigma^0))$ where $\sigma^0$ is the constant map onto $0$. Then, $a_0\leq a_1$ implies that $g_0(\sigma^0)\wedge\bigwedge_{i\in J}b_i\leq g_1(\sigma^0) \wedge\bigwedge_{i\in J}b_i$. If $\Xi_0(g_0(\sigma^0))=0$ then it is clear that $\Xi'_0(a_0) \leq \Xi'_1(a_1)$, otherwise it follows by~\ref{m80b2}.
    %    
    %There are $a\in\cB_0$, $b\in \cB_1$ and $\sigma,\tau\in \Fn(I,2)$ such that $a'=a\wedge\bigwedge_{i\in \dom\sigma} b_i^{\sigma(i)}$ and $b'=b\wedge\bigwedge_{i\in \dom\tau} b_i^{\tau(i)}$. In the case that either $\Xi_0(a)=0$ or $\sigma(j)=1$ for some $j\in \dom\sigma$, we get $\Xi'_0(a')=0\leq \Xi'_1(b')$; otherwise, if $\sigma(i)=0$ for all $i\in \dom\sigma$, then $\tau(i)=0$ for all $i\in\dom\tau$ (because $\set{b_i}{i\in I}$ generates a filter), so $a'=a\wedge\bigwedge_{i\in \dom\sigma} b_i$, $b'=b\wedge\bigwedge_{i\in \dom\tau} b_i$, $\Xi'_0(a')=\Xi_0(a)$ and $\Xi'_1(b') =\Xi_1(b)$. Now $a'\leq b'$ implies that $a\wedge\bigwedge_{i\in J}b_i \leq b\wedge\bigwedge_{i\in J}b_i$ where $J\coloneqq \dom\sigma\cup\dom\tau$, so by~\ref{m80b2} we obtain $\Xi'_0(a')=\Xi_0(a)\leq \Xi_1(b)= \Xi'_1(b')$.
\end{PROOF}

We wonder if \autoref{m80} is valid in the context of the uap when both $\cB_0$ and $\cB_1$ contain all the singletons. To show this, it is enough to solve the following problem.

\begin{question}\label{m800}
   In the context of \autoref{m80}, assume that $\cB=\pts(X)$ and that both $\cB_0$ and $\cB_1$ contain all the singletons. If both $\Xi_0$ and $\Xi_1$ are uniformly supported, can $\Xi$ in~\ref{m80i} be found uniformly supported?
\end{question}

Notice that a positive answer to \autoref{Qusup} solves the problem above (with the same proof of \autoref{m80}).

\section{Riemann integration on Boolean algebras}\label{3.5}

In this section, we develop a Riemann-type integral on Boolean algebras with respect to finitely additive measures and study some of its fundamental properties. Notably, these properties largely coincide with those of the Riemann integral on the real numbers. In particular, we show that our integral extends the Riemann integral on rectangles in $ \mathbb{R}^n $.

For this section, fix a Boolean subalgebra $\cB$ of $\calP(X)$ for some non-empty set $X$ and a finitely additive measure $\Xi \colon \cB \to [0, \delta],$ where $\delta$ is a non-negative real number.
We start defining partitions and its refinements. 

\subsection{Integration over fields of sets}\label{4.1}

\begin{definition}\label{t5}
    Recall that $\Pbf^\Xi$ is the set of finite partitions of $X$ into sets in $\dom(\Xi)=\cB$. 
    
    \begin{enumerate}[label=\normalfont(\arabic*)]
        \item If $P, Q \in \bfP^{\Xi}$, we say that \emph{$Q$ is a refinement of $P$}, denoted by $Q \ll P$, if every element of $P$ can be finitely partitioned into elements of $Q.$ \index{$\ll$} \index{refinement}

        \item If $P$ and $Q$ are in $\bfP^{\Xi},$ we define $P \sqcap Q \coloneqq \bigcup \{ a\cap b \colon a\in P,\ b\in Q \}.$
        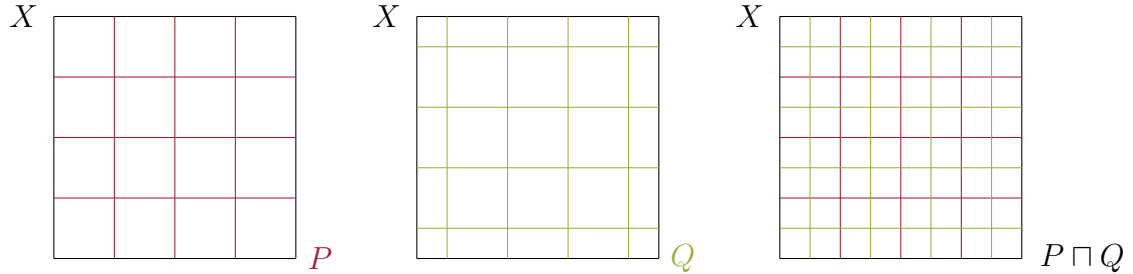
\begin{figure}[ht]
        \centering
        \begin{tikzpicture}[scale=0.8]
            % primero 
            \draw (0, 0) -- (0, 4);
            \draw (0, 0) -- (4, 0);
            \draw (4, 0) -- (4, 4);
            \draw (0, 4) -- (4, 4);

            % segundo
            \draw (6, 0) -- (6, 4);
            \draw (6, 0) -- (10, 0);
            \draw (10, 0) -- (10, 4);
            \draw (6, 4) -- (10, 4);

            %tercero 
            \draw (12, 0) -- (12, 4);
            \draw (12, 0) -- (16, 0);
            \draw (16, 0) -- (16, 4);
            \draw (12, 4) -- (16, 4);     

            % P horizontales
            \draw[redun] (0, 1) -- (4, 1);
            \draw[redun] (0, 2) -- (4, 2);
            \draw[redun] (0, 3) -- (4, 3);
            % P verticales
            \draw[redun] (1, 0) -- (1, 4);
            \draw[redun] (2, 0) -- (2, 4);
            \draw[redun] (3, 0) -- (3, 4);

            % Q horizontales
            \draw[greenun] (6, 0.5) -- (10, 0.5);
            \draw[greenun] (6, 1.5) -- (10, 1.5);
            \draw[greenun] (6, 2.5) -- (10, 2.5);
            \draw[greenun] (6, 3.5) -- (10, 3.5);
            
            % P horizontales
            \draw[greenun] (6.5, 0) -- (6.5, 4);
            \draw[greenun] (7.5, 0) -- (7.5, 4);
            \draw[greenun] (8.5, 0) -- (8.5, 4);
            \draw[greenun] (9.5, 0) -- (9.5, 4);

            % PQ horizontales
            \draw[redun] (12, 1) -- (16, 1);
            \draw[redun] (12, 2) -- (16, 2);
            \draw[redun] (12, 3) -- (16, 3);
            % PQ verticales
            \draw[redun] (13, 0) -- (13, 4);
            \draw[redun] (14, 0) -- (14, 4);
            \draw[redun] (15, 0) -- (15, 4);

            % PQ horizontales
            \draw[greenun] (12, 0.5) -- (16, 0.5);
            \draw[greenun] (12, 1.5) -- (16, 1.5);
            \draw[greenun] (12, 2.5) -- (16, 2.5);
            \draw[greenun] (12, 3.5) -- (16, 3.5);
            
            % PQ horizontales
            \draw[greenun] (12.5, 0) -- (12.5, 4);
            \draw[greenun] (13.5, 0) -- (13.5, 4);
            \draw[greenun] (14.5, 0) -- (14.5, 4);
            \draw[greenun] (15.5, 0) -- (15.5, 4);

            \node at (-0.5, 4) {$X$};
            \node at (5.5, 4) {$X$};
            \node at (11.5, 4) {$X$};

            \node[redun] at (4.4,0) {$P$};
            \node[greenun] at (10.4, 0) {$Q$};
            \node at (17, 0) {$P \sqcap Q$};
        \end{tikzpicture}
        
        \caption{A graphic example of $P \sqcap Q.$}
        \label{f46}
    \end{figure}    
    \end{enumerate}
\end{definition}

For example, it is clear that $\{ X \} \in \bfP^{\Xi}$ and if $P \in \bfP^{\Xi},$ then $P \ll \{ X \}$ and $P \ll P.$ Moreover, $\ll$ is a partial order on $\bfP^{\Xi}.$ Also, for $P, Q \in \bfP^{\Xi}$, $P \sqcap Q$ is a common refinement of $P$ and $Q$.  

Now we can define the integral with respect to $\Xi$:

\begin{definition}\label{t14}
  \index{$\supsum(f, P)$} \index{$\infsum(f, P)$}
  
  Let $f\colon X\to\bbR$ be a bounded function. We define:
  
  \begin{enumerate}   [label=\normalfont(\arabic*)]   
      \item For any $P\in\Pbf^\Xi$, 
      \[\supsum^{\Xi}(f,P) \coloneqq \sum_{b\in P}\sup(f[b])\Xi(b)\text{\ and }\infsum^{\Xi}(f,P) \coloneqq \sum_{b\in P}\inf(f[b])\Xi(b).\] 
      
      \item $\overline{\int_{X}} fd\Xi \coloneqq \inf\{\supsum(f,P) \colon  P\in\Pbf^\Xi\}$ \ and \
      $\underline{\int_{X}}fd\Xi \coloneqq \sup\{\infsum(f,P) \colon  P\in\Pbf^\Xi\}$. \index{$\overline{\int}_{X} d \Xi$} \index{$\underline{\int}_{X} f d \Xi$}
      
      \item We say that $f$ is \emph{$\Xi$-integrable}, denoted by $f  \in \cI(\Xi)$ if, and only if, $\overline{\int_{X}}fd\Xi=\underline{\int_{X}}fd\Xi,$ in which case this value is denoted  by $\int_{X} fd\Xi$. \index{$\int f d \Xi$} \index{$\cI(\Xi)$}
  \end{enumerate}
\end{definition}

Naturally, when the context is clear, we omit the superscript ``$\Xi$'' in ``$\supsum^{\Xi}(f, P)$'' and ``$\infsum^{\Xi}(f, P)$''. 

For example, it is clear that any constant function is $\Xi$-integrable. Concretely, if for all $x \in X, $ $ f(x) = c \in \bbR,$ then $\int_{X} f(x) d \Xi = c\, \Xi(X)$.

% \blueq{
% \begin{question}
%     Más aún, si fijamos un conjunto $X$, entonces las funciones constantes son $\Xi$-integrables para toda fam $\Xi$ cuyo dominio es un campo sobre $X$. A una función que sea integrable para toda fam definida sobre un campo sobre $\calP(X)$ le podríamos llamar \emph{$X$-universal} o algo así. La pregunta es si podemos caracterizar todas las funciones $X$-universales para un $X$ dado.  Creería que precisamente solo las constantes son $X$-universales. \DM{Eso que dices equivale a integrar sobre el campo $\{\emptyset,X\}$, pero solo las funciones constantes son integrables ah\'i (a menos que $\Xi(X)=0$). Quiz\'a sea m\'as interesante pensar en $\cB$-universales para un \'algebra $\cB$ dada. Por ejemplo, las $\pts(X)$-universales son las funciones acotadas.}
% \end{question}
% }

The following result is~\cite[Prop.~4.5.2]{BhaskaraRa}, but the proof is exactly the same as in the case of the Riemann integral.

\begin{lemma}\label{t18}
    Let $f \colon X \to \bbR$ a bounded function. If $P, Q \in \bfP^{\Xi}$ and $Q \ll P,$ then: $$ \infsum(f, P) \leq \infsum(f, Q) \leq \supsum(f, Q) \leq \supsum(f, P).$$ As a consequence, $\supsum(f, Q) - \infsum(f, Q) \leq \supsum(f, P) - \infsum(f, P).$
\end{lemma}

Using $P \sqcap Q$, we obtain:

\begin{corollary}\label{t22}
    If $P, Q \in \bfP^{\Xi},$ then $\infsum(f, P) \leq \supsum(f, Q).$ As a consequence, 
    $\underline{\int_{X}} f d  \Xi \leq \overline{\int_{X}} f d  \Xi.$ 
\end{corollary}

We now show that our $\Xi$-integrability generalizes the Riemann integral over rectangles in $\bbR^{n}$. Recall from \autoref{n364} that, for $E\subseteq X$, we denote $\cB|_E\coloneqq \set{E\cap b}{b\in\cB}$, which is a Boolean subalgebra of $\calP(E)$. In the case $E\in\cB$, we denote $\Xi|_E\coloneqq \Xi\frestr \cB|_E$. 

% We often use the following notation in our integration theory.

% \begin{notation}\label{t24}
% For $E\subseteq X$, denote $\cB|_E\coloneqq \set{E\cap b}{b\in\cB}$, which is a Boolean subalgebra of $\calP(E)$. In the case $E\in\cB$, we denote $\Xi|_E\coloneqq \Xi\frestr \cB|_E$. 
% \end{notation}

\begin{example}\label{t16}
    Let $n\geq 1$ be a natural number. For $a,b\in \bbR^n$, we write $a\leq b$ for $a_i\leq b_i$ for all $i<n$. We denote
    \[\begin{split}
        [a,b] & \coloneqq \prod_{i<n}[a_i,b_i],\\
        [a,b) & \coloneqq \prod_{i<n}[a_i,b_i),
    \end{split}\]
    and the rectangles $(a,b]$ and $(a,b)$ are defined in the natural way.

    Let $\cC^n$ be the collection of subsets of $\bbR^n$ that can be written as a finite union of rectangles of the form $[x,y)$ with $x\leq y$ in $\bbR^n$. Note that any set in $\cC^n$ can be written as a finite disjoint union of such rectangles. Although $\cC^n$ has $\emptyset$ as an element and it is closed under $\cup$, $\cap$ and $\smallsetminus$ (set difference), it is not a Boolean algebra because $\bbR^n \notin \cC^n$.

    There is a unique function $\Xi^n\colon\cC^n\to [0,\infty)$ that satisfies the properties of a fam and such that
    \[\Xi^n([x,y))=\prod_{i<n}(y(i)-x(i))\]
    for any $x\leq y$ in $\bbR^n$. In fact, $\Xi^n$ is $\sigma$-additive. If $a\leq b$ in $\bbR^n$ then, allowing the use of \autoref{n364}, we get that $\cC^n|_{[a,b]}$ is a Boolean subalgebra of $\calP([a,b])$ and $\Xi^n|_{[a,b]}$ is a fam on $\cC^n|_{[a,b]}$, which is naturally defined as the volume of the rectangles.

    It is not hard to show that a bounded function $f\colon [a,b]\to \bbR$ is Riemann integrable in $\bbR^{n}$ iff it is $\Xi^n|_{[a,b]}$-integrable. 
\end{example}

Another example comes from ultrafilters. 

\begin{example}\label{t17}
    Let $U$ be an ultrafilter on $\cB$ and consider $\Xi_U$ as its associated measure (see \autoref{2.1}). If $f\colon X\to\bbR$ is bounded, then
    \[\underline{\int_X} f d \Xi_U = \sup_{a\in U} \{\inf 
    (f[a])\} \text{ and }\overline{\int_X} f d \Xi_U = \inf_{a\in U} \{\sup(f[a])\}.\]
    When $f$ is $\Xi_U$-integrable we say that $f$ is \emph{$U$-integrable} and its $\Xi_U$-integral is any of the values above. The above are basically $\liminf_U f$ and $\limsup_U f$, and $U$-integrable means that $\lim_U f$ (the ultrafilter limit) exists, and the integral is this limit.

    For $x\in X$ let $u_x\coloneqq \set{b\in\cB}{x\in b}$, which is clearly an ultrafilter on $\cB$. Clearly $\underline{\int_X} f d \Xi_{u_x} \leq f(x) \leq \overline{\int_X} f d \Xi_{u_x}$, so $\int_X f d \Xi_{u_x} = f(x)$ whenever $f$ is $u_x$-integrable. 
\end{example}

The previous example motivate the following definition, which we will use to characterize integrability with respect to ultrafilters in \autoref{t25}.

\begin{definition}\label{t20}
    Let $U$ be an ultrafilter on $\cB$ and $f\colon X\to\bbR$ bounded. We define the \emph{oscillation of $f$ on $U$} by
    $$\osc_f(U)\coloneqq \inf_{a\in U}\{ \sup(f[a]) - \inf(f[a]) \}.$$
\end{definition}

Our notion of $\Xi$-integrability can also be expressed in terms of sums similar to Riemann sums.

\begin{definition}
    Let $f\colon X\to \bbR$ (not necessarily bounded). When $\tau\in\prod P\coloneqq \prod_{b\in P}b$ for some $P\in \bfP^{\Xi}$, we denote
    \[\Sm^\Xi(f,\tau)\coloneqq \sum_{b\in P}f(\tau(b))\Xi(b).\]
    Such a sum is called a \emph{$\Xi$-sum of $f$}.
\end{definition}

Now, we prove what we call the \emph{Criterion of $\Xi$-Integrability}: \index{criterion of $\Xi$-integrability}

\begin{theorem}\label{t29}
    Let $f\colon X\to \bbR$ be a bounded function. Then the following statements are equivalent.
    \begin{enumerate}[label = \normalfont(\roman*)]
        \item\label{t29i} $f$ is $\Xi$-integrable.
        \item\label{t29ii} For all $\varp > 0,$ there exists a partition $P \in \bfP^{\Xi}$ such that $ \supsum(f, P) - \infsum(f, P) < \varp$.
        \item\label{t29iii} There is some $L \in\bbR$ such that, for any $\varp>0$, there is some $P\in\bfP^\Xi$ such that, for all $Q\ll P$ in $\bfP^\Xi$,
        \[|\Sm^\Xi(f,\tau)-L|<\varp \text{ for all }\tau\in\prod Q.\]
        \item\label{t29iv} There is some $L \in\bbR$ such that, for any $\varp>0$, there is some $P\in\bfP^\Xi$ such that
        \[|\Sm^\Xi(f,\tau)-L|<\varp \text{ for all }\tau\in\prod P.\]
    \end{enumerate}
    When~\ref{t29iii} and~\ref{t29iv} hold, we actually have $L=\int_X f d \Xi$ (so such an $L$ is unique in both statements).
\end{theorem}

\begin{PROOF}{\ref{t29}}
    We first show \ref{t29i}${}\Rightarrow{}$\ref{t29ii}. Assume that $f \, {\in} \, \cI(\Xi)$ and let $\varp > 0.$ By basic properties of $\sup$ and $\inf,$ there are $P, Q \in \bfP^{\Xi}$ such that: $$ \int_{X} f d \Xi - \frac{\varp}{2} < \infsum(f, P) \text{ and } \supsum(f, Q) < \int_{X} f d \Xi + \frac{\varp}{2}. $$ 
    Consider $R \coloneqq P \sqcap Q$, which is a common refinement of $P$ and $Q$. So, by virtue of  \autoref{t18}, 
    $$\infsum(f, P) \leq \infsum(f, R) \text{ and } \supsum(f, R) \leq \supsum(f, Q).$$
    Therefore, 
    $$ \int_{X} f - \frac{\varp}{2} < \infsum(f, R) \text{ and } \supsum(f, R) < \int_{X} f d \Xi + \frac{\varp}{2}. $$ 
    Thus, $\supsum(f, R) - \infsum(f, R) < \varp.$ 

    For \ref{t29ii}${}\Rightarrow{}$\ref{t29i}, let $P \in \bfP^{\Xi}$ such that $\supsum(f, P) - \infsum(f, P) < \varp.$ Hence, by the definition of $\overline{\int}$ and $\underline{\int}$, we have that 
    $$ \overline{\int_{X}} f d \Xi \leq \supsum(f, P) < \infsum(f, P) + \varp \leq \underline{\int_{X}} f d \Xi + \varp.$$ 
    Since $\varp$ is arbitrary, by \autoref{t22} it follows that $f \in \cI(\Xi).$

    We now show \ref{t29i}${}\Rightarrow{}$\ref{t29iii}. Assume~\ref{t29i} and let $L\coloneqq \int_X f d \Xi$. For $\varp>0$ we have some $P\in\bfP^\Xi$ such that, for all $Q\ll P$ in $\bfP^\Xi$, $L-\varp<\infsum(f,Q)\leq \supsum(f,Q)<L+\varp$. For any such $Q$, $\infsum(f,Q)\leq\Sm^\Xi(f,\tau)\leq\supsum(f,Q)$ for all $\tau \in\prod Q$, so $|\Sm^\Xi(f,\tau)-L|<\varp$.

    The implication \ref{t29iii}${}\Rightarrow{}$\ref{t29iv} is obvious.
    To finish the proof, we assume~\ref{t29iv} and prove that $f$ is integrable and $L=\int_X f d \Xi$. Fix $\varp>0$ and choose a partition $P\in\bfP^\Xi$ as in~\ref{t29iii}. For any $b\in P$, pick some $\tau(b)\in b$ such that $\sup(f[b])-f(\tau(b))<\varp$. Then
    \[
        L+\varp > \Sm(f,\tau) > \sum_{b\in P}(\sup(f[b])-\varp)\Xi(b) = \supsum(f,P) - \varp \Xi(X),
    \]
    so $\overline{\int_X}f d \Xi\leq \supsum(f,P)< L+(1+\Xi(X))\varp$.
    Since $\varp$ is arbitrary, we conclude $\overline{\int_X}f d \Xi \leq L$. A similar argument yields $L\leq \underline{\int_X} f d \Xi$, therefore $\overline{\int_X}f d \Xi= \underline{\int_X} f d \Xi = L$.
\end{PROOF}

In the case of ultrafilters as discussed in \autoref{t17}, we obtain:

\begin{corollary}\label{t25}
    Let $U$ be an ultrafilter on $\cB$ and $f\colon X\to\bbR$ bounded. Then $f$ is $U$-integrable iff $\osc_f(U)=0$.
\end{corollary}

% unterl here

\begin{example}\label{t33}
    Assume that $A$ is a set of atoms of $\cB$ with $\Xi$-positive measure such that $\bigcup A \in\cB$ and $X\menos\bigcup A$ has measure zero. Notice that, by \autoref{m50}, $A$ is countable.

    If $f\colon X\to \R$ is bounded, then $f\in\cI(\Xi)$ implies that $f\frestr a$ is constant for all $a\in A$. Indeed, for $a\in A$ and $\varp>0$, we can find some $P\in\Pbf^\Xi$ including $a$ (after refining, also because $a$ is an atom) such that $\supsum^\Xi(f,P)-\infsum^{\Xi}(f,P)<\varp$. Then, for $x,y\in a$, 
    \[|f(x)-f(y)|\, \Xi(a) \leq (\sup(f[a])-\inf(f[a]))\, \Xi(a) \leq \supsum^\Xi(f,P)-\infsum^{\Xi}(f,P)<\varp,\]
    so $|f(x)-f(y)|<\frac{\varp}{\Xi(a)}$. Since $\varp$ is arbitrary, we get that $f(x)=f(y)$. Therefore, $f$ is constant on $a$. 

    The converse requires that either $A$ is finite or $\Xi$ is $\sigma$-additive. In such cases, if $f$ has constant value $c_a$ for $a\in A$, then $\int_X f d\Xi = \sum_{a\in A}c_a \Xi(a)$. The latter sum is absolutely convergent because $\sum_{a\in A}\Xi(a) \leq \Xi(X)$ and the sequence $\Seq{c_a}{a\in A}$ is bounded.
    
    To show this, since $f$ is bounded, we can find some $M>0$ such that $|f(x)|\leq M$ for all $x\in X$. Let $\varp>0$ and pick a finite $F\subseteq A$ such that $\sum_{a\in A\menos F}\Xi(a) < \frac{\varp}{2M}$, but in the case that $A$ is finite set $F\coloneqq A$. 
    Let $a'\coloneqq X\menos \bigcup F$ and $P\coloneqq F\cup \{ a' \}$, which is in $\Pbf^\Xi$. Notice that $\Xi(a')<\frac{\varp}{2M}$ (by cases on whether $F=A$ or $\Xi$ is $\sigma$-additive). Then, 
    for any $\tau\in \prod P$, $S^\Xi(f,\tau) = f(\tau(a')) + \sum_{a\in F}c_a \Xi(a)$, so
    \begin{align*}
        \left|S^\Xi(f,\tau) - \sum_{a\in A}c_a\Xi(a)\right| & = \left| \sum_{a\in A\menos F}c_a\Xi(a) + f(\tau(a'))\Xi(a') \right| \leq \left| \sum_{a\in A\menos F}c_a\Xi(a)\right| + M \Xi(a')\\ 
         & \leq M\sum_{a\in A\menos F}\Xi(a) + M \Xi(a') < \frac{\varp}{2} + \frac{\varp}{2} = \varp.
    \end{align*}
    Therefore, by \autoref{t29}, $f\in \cI(\Xi)$ and $\int_X f d\Xi = \sum_{a\in A}c_a \Xi(a)$.
\end{example}

As a particular case of \autoref{t33}, we can explicitly calculate $\int_{X} f d \Xi^{u}$, where $\Xi^{u}$ is as in \autoref{m37}. 

\begin{example}\label{t42}
    Let $f \colon X \to \bbR$ be a bounded function. If $u\subseteq X$ is finite and non-empty, then $f$ is $\Xi^u$-integrable and
    $$\int_{X} fd\Xi^u=\frac{1}{|u|}\sum_{k\in u} f(k).$$ 

    Indeed, consider $A \coloneqq \{ \{ x \} \colon x \in u \}$, which is a finite set of atoms of $\calP(X)$ with positive measure, $\bigcup A = u \in \calP(X)$ and $\Xi^{u}(X \setminus \bigcup A) = 0$. Since $A$ is finite, and for any $a = \{ x \} \in A$ the value $c_{a}$ from \autoref{t33} is $f(x)$, we can use it to conclude that $f$ is $\Xi^u$-integrable and calculate: 
    $$
    \int_{X} f d \Xi^{u} = \sum_{a \in A} c_{a}  \Xi(a) = \frac{1}{\vert u \vert} \sum_{x \in u} f(x). 
    $$
\end{example}

We could have had defined integrability as in \autoref{t29}~\ref{t29iii} without demanding that $f$ is bounded. However, it is very curious that this condition actually implies that $f$ must be bounded except on a measure zero set, as we prove below.

\begin{lemma}\label{t30}
    Assume that $f\colon X\to \bbR$ satisfies~\ref{t29iv} of \autoref{t29}. Then there is some $b_0\in\cB$ with $\Xi(b_0)=0$ such that $f{\upharpoonright}(X\smallsetminus b_0)$ is bounded.
\end{lemma}
\begin{PROOF}{\ref{t30}}
    Let $L$ be as in \autoref{t29}~\ref{t29iv}, and choose $P\in\bfP^\Xi$ as in this condition for $\varp=1$, i.e.\
    \begin{equation}\label{eqt30}
        |\Sm^\Xi(f,\tau)-L|<1 \text{ for all }\tau\in\prod P.
        \tag{$\blacktriangle$}
    \end{equation}
    It is enough to show that any $b\in P$ such that $f{\upharpoonright}b$ is unbounded has $\Xi$-measure $0$. For such a $b$, wlog assume that $\sup(f[b]) = \infty$ (the case $\inf(f[b])= -\infty$ is dealt with similarly), and towards a contradiction, assume that $\Xi(b)>0$. For any $b'\neq b$ in $P$, choose any $\tau(b')\in b'$, and pick $\tau(b)\in b$ such that
    \[f(\tau(b))>\frac{L+1-\sum_{b'\in P\menos\{b\}} f(\tau(b'))\Xi(b')}{\Xi(b)}.\]
    This implies $\Sm^\Xi(f,\tau)-L>1$, which contradicts~\eqref{eqt30}.
\end{PROOF}

Our notion of integration can be extended to functions that are bounded $\Xi$-almost everywhere: if $f\colon X\to\bbR$ and it is bounded outside some $\Xi$-measure zero set $b_0\in\cB$, we can say that $f$ is $\Xi$-integrable iff $f{\upharpoonright}(X\menos b_0)$ is $\Xi|_{X\menos b_0}$-integrable, and in such a case $\int_X f d\Xi \coloneqq \int_{X\menos b_0} f d \Xi|_{X\menos b_0}$. This discussion makes sense by \autoref{t31}.

\subsection{Basic properties of the integral}\label{4.2}

We now deal with some basic properties of $\Xi$-integrability. We start with the linearity and order properties. 

\begin{lemma}\label{t47}
    Let $f, g \in \cI(\Xi)$ and $c \in \bbR.$ Then $cf, f+g \in \cI(\Xi)$ and:  
    \begin{enumerate}[label=\normalfont(\alph*)]
        \item\label{t47a} $\int_{X} (cf) d \Xi = c \int_{X} f d \Xi.$ 

        \item\label{t47b} $\int_{X}(f + g) d \, \Xi = \int_{X} f d \, \Xi + \int_{X} g d \Xi.$ 
    \end{enumerate}
\end{lemma}

\begin{PROOF}{\ref{t47}}
    The reader can verify that the proof follows exactly the same procedure as the analogous result for the Riemann integral on the real line, as it is carried out, for example, in~\cite[Ch.VI,~\S2]{Rosenlicht}. 
\end{PROOF}

% So, inductively, we get: 

% \begin{corollary}\label{t50}
%     Let $\{ f_{i} \colon i < n \}$ a finite sequence of $\Xi$-integrable functions. Then $\sum_{i < n} f_{i} \in \cI(\Xi)$ and $$\int_{X} \left( \sum_{i < n} f_{i}\right) d \Xi = \sum_{i < n} \left( \int_{X} f_{i} d \Xi \right).$$
% \end{corollary}

\begin{corollary}\label{t53}
    If $f, g \colon X \to \bbR$ are $\Xi$-integrable functions and $f \leq g,$ then $ \int_{X} f d \Xi \leq \int_{X} g d \Xi.$
\end{corollary}

\begin{PROOF}{\ref{t53}}
    This follows an standard argument, see e.g. \cite[\S 7.6]{apostol1974}.
    % Prueba detallada
    % For all $x \in X,$ define $h(x) \coloneqq g(x) - f(x).$ So $h \geq 0.$ It is clear that, for every $P \in \bfP^{\Xi}, \, \supsum(h, P) \geq 0.$ Finally, since by \autoref{t47} $h \in \cI(\Xi),$ we have that $$\int_{X} g d \Xi - \int_{X} f d \Xi = \int_{X} h d \Xi \geq 0.$$ Thus $\int_{X} f d \Xi \leq \int_{X} g d \Xi.$
\end{PROOF}

\begin{lemma}\label{t55}
    If $f \in \cI(\Xi)$, then $\vert f \vert \in \cI(\Xi)$ and $ \left \vert \int_{X} f d \Xi \right \vert \leq \int_{X} \vert f \vert d \Xi. $ 
\end{lemma}

\begin{PROOF}{\ref{t55}}
    This follows an standard argument using the positive and the negative parts of $f$. For more details, see e.g. ~\cite[\S 7.6]{apostol1974}. 
    %
    % For any integrable function $f \in \cI(\Xi)$, the positive and negative parts 
    % $f_{+}$ and $f_{-}$ are also $\Xi$-integrable, and 
    % $\lvert f \rvert = f_{+} + f_{-}$ satisfies 
    % \[
    % \left \vert \int_{X} f d\Xi \right \vert \le \int_{X} \lvert f \rvert d\Xi.
    % \]
    % For more details see e.g. ~\cite[Sec.~7.6]{apostol1974}. \qedhere
    % Prueba detallada
    % Define $f_{+} \colon X \to \bbR$ such that, for any $x \in X$, $f_{+}(x) \coloneqq \max \{ f(x) , 0 \}$. Similarly, consider $f_{-} \colon X \to \bbR$ such that $f_{-}(x) \coloneqq f_{+}(x) - f(x)$. We can write $\vert f(x) \vert = f_{+}(x) + f_{-}(x)$ for any $x \in X$.   It is clear that, if $P \in \bfP^{\Xi}$, then 
    % $$ \supsum(f_{+}, P) - \infsum(f_{+}, P) =\sum_{b\in P}(\sup f_+[b]-\inf f_+[b])\Xi(b) \leq \supsum(f, P) - \infsum(f, P),$$ and therefore, by \autoref{t29} and \autoref{t47}, $f \in \cI(\Xi)$ implies $f^{+}, f_{-} \in \cI(\Xi)$. On the other hand, since $\vert f \vert = f_{+} + f_{-}$, we can conclude that $\vert f \vert \in \cI.$  Finally, 
    % $$ \left \vert \int_{X} f d \Xi \right \vert = \left \vert \int_{X} f_{+} d \Xi - \int_{X} f_{-} d \Xi \right \vert \leq  \int_{X} f_{+} d \Xi  +  \int_{X} f_{-} d \Xi = \int_{X} \vert f \vert d \Xi.\qedhere$$
\end{PROOF}

Towards showing that the product of two $\Xi$-integrable functions is $\Xi$-integrable, we prove the following result. 

\begin{lemma}\label{t58}
    If $f \in \cI(\Xi)$, then $f^{2} \in \cI(\Xi).$  As a consequence, the product of $\Xi$-integrable functions is $\Xi$-integrable. 
\end{lemma}

\begin{PROOF}{\ref{t58}}
    %Prueba corta:  

    Since $f$ is bounded, there exists some $0 < N < \omega$ such that $\vert f(x) \vert \leq N$ for all $x \in X$. 
    On the other hand, since $f$ is $\Xi$-integrable, the usual argument based on the estimate 
    $|f^{2}(c) - f^{2}(d)| \le 2N\,|f(c) - f(d)|$ shows that $f^{2}$ is $\Xi$-integrable (see e.g. \cite[\S 7.6]{apostol1974}). 
\end{PROOF}

We now study the relation between the lower and upper integrals of $f$ with respect to two fams $\Xi_{0} \subseteq \Xi_{1}$. %such that $\Xi_{0}$ is defined on $\cB_{0}$, $\Xi_{1}$ is defined on $\cB_{1}$ and $\cB_{0} \subseteq \cB_{1}$. 

\begin{lemma}\label{t93}
    Let $\cB_{0}, \cB_{1} \subseteq \calP(X)$ be Boolean algebras such that $\cB_{0} \subseteq \cB_{1},$ and let $\Xi_{0}, \Xi_{1}$ be finite finitely additive measures on $\cB_{0}, \cB_{1}$ respectively, such that $\Xi_{0} \subseteq \Xi_{1}.$ Let $f \colon X \to \bbR$ be a bounded function. Then, 
    \begin{equation*}
        \underline{\int_{X}} f d \Xi_{0}\leq\underline{\int_{X}} f d \, \Xi_{1}\leq\overline{\int_{X}} f d \, \Xi_{1}\leq \overline{\int_{X}} f d \Xi_{0}.
    \end{equation*}
    As a consequence, if $f \in \cI(\Xi_{0})$ then $f \in \cI(\Xi_{1})$, and %\DM{Incluir desigualdades de las integrales superiores e inferiores.}
    $$
    \int_{X} f d \Xi_{0} = \int_{X} fd \Xi_{1}.
    $$ 
\end{lemma}

\begin{PROOF}{\ref{t93}}
    Since $\Xi_{0} \subseteq \Xi_{1}$, for any $P \in \bfP^{\Xi_{0}}$, $\{\supsum^{\Xi_{0}}(f, P) \colon P \in \bfP^{\Xi_{0}}\} \subseteq \{ \supsum^{\Xi_{1}}(f, P) \colon P \in \bfP^{\Xi_{1}} \}$ and $ \{ \infsum^{\Xi_{0}}(f, P) \colon P \in \bfP^{\Xi_{0}} \} \subseteq \{ \infsum^{\Xi_{1}}(f, P) \colon P \in \bfP^{\Xi_{1}} \}.$ Therefore,
    \begin{equation*}
           \underline{\int_{X}} f d \Xi_{0}\leq\underline{\int_{X}} f d \, \Xi_{1}\leq\overline{\int_{X}} f d \, \Xi_{1}\leq \overline{\int_{X}} f d \Xi_{0}.
        \end{equation*}

    %Therefore, $$ \overline{\int_{X}} f d \, \Xi_{1} = \inf \{ \supsum^{\Xi_{1}}(f, P) \colon P \in \bfP^{\Xi_{1}} \} \leq \inf \{ \supsum^{\Xi_{0}}(f, P) \colon P \in \bfP^{\Xi_{0}} \} = \overline{\int_{X}} f d \Xi_{0}, $$ and $$ \underline{\int_{X}} f d \Xi_{0}  = \sup \{ \infsum^{\Xi_{0}}(f, P) \colon P \in \bfP^{\Xi_{0}} \} \leq \sup \{ \infsum^{\Xi_{1}}(f, P) \colon P \in \bfP^{\Xi_{1}} \} = \underline{\int_{X}} f d \, \Xi_{1}.$$ Whence, it follows that: 
    
        % \begin{equation}\label{e45}
        %     \text{$ \underline{\int_{X}} f d \, \Xi_{0}  \leq \underline{\int_{X}} f d \Xi_{1} \leq \overline{\int_{X}} f d \Xi_{1} \leq \overline{\int_{X}} f d \, \Xi_{0}.$}
        % \end{equation}
    
    Clearly, if $f$ is $\Xi_{0}$-integrable, then $f$ is $\Xi_{1}$-integrable and the value of the integrals is the same.
    %
    %For the converse, notice that if $\cB_{0} = \cB_{1},$ then $\bfP^{\Xi_{0}} = \bfP^{\Xi_{1}}$. Hence, inclusions in the sets above are really equalities, which implies that inequalities above are equalities as well. Thus, $f \in \cI(\Xi_{0}) \Leftrightarrow f \in \cI(\Xi_{1}),$ and  it is clear that both integrals coincide. 
    %
    % \begin{equation*}
    %         \text{$ \underline{\int_{X}} f d \, \Xi_{0}  = \underline{\int_{X}} f d \Xi_{1} = \overline{\int_{X}} f d \Xi_{1} = \overline{\int_{X}} f d \, \Xi_{0}.$}
    %     \end{equation*}
    %
\end{PROOF}

We study some results about $\Xi$-integration when we transfer information between Boolean algebras. Recall from \autoref{s16} that, if $h \colon X \to Y$ is a function and $\cB$ is a Boolean subalgebra of $\calP(X)$, then the collection $h^{\to}(\cB) \coloneqq \{ A \subseteq Y \colon h^{-1}[A] \in \cB \}$ is a Boolean subalgebra of $\calP(Y)$ and, if $\Xi$ is a finitely additive measure on $\cB$, then $\Xi_{h}$ is a finitely additive measure on $h^{\to}(\cB)$. 

\begin{lemma}[{\cite[Lem.~3.16]{PU2025}}]\label{s5}
    Let $X$, $Y$ be non-empty sets, $h \colon X \to Y$ a function, $\cB_{X}$ a Boolean subalgebra of $\calP(X)$, $\cB_{Y} \coloneqq h^{\to}(\cB_{X})$, and $\Xi$ a finitely additive measure on $\cB_{X}$. Let $f \colon Y \to \bbR$ be a bounded function. Then:

    \begin{enumerate}[label = \normalfont (\alph*)]
        \item\label{s5.a} For any $P \in \bfP^{\Xi_{h}}$, there are $P^{\bullet} \in \bfP^{\Xi_{h}}$ and $Q \in \bfP^{\Xi}$ such that $P^{\bullet} \ll P$ and $$ \supsum^{\Xi_{h}}(f, P^{\bullet}) = \supsum^{\Xi}(f \circ h, Q) \text{ and }  \infsum^{\Xi_{h}}(f, P^{\bullet}) = \infsum^{\Xi}(f \circ h, Q).$$

        \item\label{s5.b} If $h$ is one-to-one, then for any $Q \in \bfP^{\Xi}$ there exists some $P \in \bfP^{\Xi_{h}}$ such that $$\supsum^{\Xi}(f \circ h, Q) = \supsum^{\Xi_{h}}(f, P) \text{ and } \infsum^{\Xi}(f \circ h, Q) = \infsum^{\Xi_{h}}(f, P).$$ 
    \end{enumerate}
\end{lemma}

\begin{PROOF}{\ref{s5}}
    We include the proof for completeness. First notice that $\ran h\in h^\to(\cB_X)$.
    
    \ref{s5.a}: Let $P \in \bfP^{\Xi_{h}}$. Define $P^{\bullet} \coloneqq P \sqcap \{ \ran h, Y\menos \ran h\}$. Hence, $P^{\bullet} \in \bfP^{\Xi_{h}}$, it is a refinement of $P$ and, for any $A \in P^{\bullet}$, either $A \subseteq \ran h$ or $A \cap \ran h = \emptyset$. Consider the set $R \coloneqq \{ A \in P^{\bullet} \colon A \subseteq \ran h \}$. Notice that, if $A \in R$ then $A = h[h^{-1}[A]]$ and, otherwise, $\Xi_{h}(A) = \emptyset$. Define $Q \coloneqq \{ h^{-1}[A] \colon A \in R \}$. Clearly $Q \in \bfP^{\Xi}$. As a consequence,  
    \begin{equation*}
        \begin{split}
            \supsum^{\Xi_{h}}(f, P^{\bullet}) & = \sum_{A \in P^{\bullet}} \sup(f[A]) \Xi_{h}(A) = \sum_{A \in R} \sup(f[A]) \Xi_{h}(A) + \sum_{A \in P^{\bullet} \setminus R} \sup(f[A]) \Xi_{h}(A)\\
            & = \sum_{A \in R} \sup(f \left[h[h^{-1}[A]] \right]) \Xi(h^{-1}[A]) = \sum_{B \in Q} \sup(f \circ h[B]) \Xi(B)\\
            & = \supsum^{\Xi}(f \circ h, Q). 
        \end{split}
    \end{equation*}

    Similarly, $ \infsum^{\Xi_{h}}(f, P^{\bullet}) = \infsum^{\Xi}(f \circ h, Q)$. 

    \ref{s5.b}: Assume that $h$ is one-to-tone. Let $Q \in \bfP^{\Xi}$ and define $P \coloneqq \{ h[B] \colon B \in Q \} \cup \{ Y\menos \ran h \}$. Clearly, $P \in \bfP^{\Xi_{h}}$ and
    \begin{equation*}
        \begin{split}
            \supsum^{\Xi_{h}}(f, P) & = \sum_{A \in P} \sup(f[A]) \Xi_{h}(A)  = \sum_{B \in Q} \sup(f[h[B]]) \Xi_{h}(h[B])\\
            & = \sum_{B \in Q} \sup(f \circ h [B]) \Xi(B) = \supsum^{\Xi}(f \circ h, Q). 
        \end{split}
    \end{equation*}
    The proof for the lower sum follows similar lines. 
\end{PROOF}

As a consequence, under the conditions  in \autoref{s5}, integrability is preserved under composition as well as the value of the integrals. 

\begin{theorem}[{\cite[Thm.~3.17]{PU2025}}]\label{s2}
   % Let $X$, $Y$ be non-empty sets, $h \colon X \to Y$ a function, $\cB_{X}$ a Boolean subalgebra of $\calP(X)$, $\cB_{Y} \coloneqq h^{\to}(\cB_{X})$, and $\Xi$ a finitely additive measure on $\cB_{X}$. Let $f \colon Y \to \bbR$ be a bounded function. Then: 
   In the context of \autoref{s5}: 
    \begin{enumerate}[label = \normalfont (\alph*)]
        \item\label{s2.1} If $f \in \cI(\Xi_{h})$ then $f \circ h \in \cI(\Xi)$ and 
         $$
        \int_{Y} f d \Xi_{h} = \int_{X} f \circ h\, d \Xi. 
        $$

        \item\label{s2.2} If $h$ is one-to-one, then $f \circ h \in \cI(\Xi)$ implies $f \in \cI(\Xi_{h})$ and 
         $$  \int_{X} f \circ h\, d \Xi = 
        \int_{Y} f d \Xi_{h}. 
        $$
    \end{enumerate}
\end{theorem}

\begin{PROOF}{\ref{s2}}
    \ref{s2.1}: Assume that $f \in \cI(\Xi_{h})$ and $\varp > 0$. So we can find a partition $P \in \bfP^{\Xi_{h}}$ such that $\supsum^{\Xi_{h}}(f, P) - \infsum^{\Xi_{h}}(f, P) < \varp$. Let $P^{\bullet} \in \bfP^{\Xi_{h}}$ and $Q \in \bfP^{\Xi}$ as in \autoref{s5}~\ref{s5.a}. By \autoref{t18} %, 
    %$$
    % \supsum^{\Xi}(f \circ h, Q) -  \infsum^{\Xi}(f \circ h, Q) = \supsum^{\Xi_{h}}(f, P^{\bullet}) - \infsum^{\Xi_{h}}(f, P^{\bullet}) \leq \supsum^{\Xi_{h}}(f, P) - \infsum^{\Xi_{h}}(f, P) < \varp.
    % $$ 
    % This shows that $f \circ h$ is $\Xi_{h}$-integrable by applying \autoref{t29}. 
    %
     % We now deal with the value of the integral. On the one hand, let $P \in \bfP^{\Xi_{h}}$. By applying \autoref{s5}~\ref{s5.a}, we can find $P^{\bullet} \ll P$  in $\bfP^{\Xi_{h}}$ and $Q \in \bfP^{\Xi}$ such that 
     and \autoref{s5}~\ref{s5.a},
    $$
       \infsum^{\Xi_{h}}(f, P^\bullet) = \infsum^{\Xi}(f \circ h, Q) \leq \underline{\int_{X}} f \circ h\, d \Xi \leq \overline{\int_{X}} f \circ h\, d \Xi \leq  \supsum^{\Xi}(f \circ h, Q) = \supsum^{\Xi_{h}}(f, P^{\bullet}).
    $$
    By \autoref{t18}, since $P^{\bullet}$ is a refinement of $P$, both $\underline{\int_{X}} f \circ h\, d \Xi$ and $\overline{\int_{X}} f \circ h\, d \Xi$ are in the interval $[\infsum^{\Xi_h}(f, P),\supsum^{\Xi_h}(f, P)]$, which has lenght ${<}\varp$. This shows that $f\circ h$ is $\Xi$-integrable and and $\int_{X} f \circ h\, d \Xi = 
        \int_{Y} f d \Xi_{h}$.

    %As a consequence, $\int_{X} f \circ h d \Xi = \int_{Y} f d \Xi_{h}$.
    
    \ref{s2.2}: Assume that $h$ is one-to-one, $f \circ h$ is $\Xi$-integrable, and let $\varp > 0$. By \autoref{t29}, there exists some $Q \in \bfP^{\Xi}$ such that $\supsum^{\Xi}(f \circ h, Q) - \infsum^{\Xi}(f \circ h, Q) < \varp$. Consider $P \in \bfP^{\Xi_{h}}$ as in \autoref{s5}~\ref{s5.b}. As a consequence, 
    $$
        \supsum^{\Xi_{h}}(f, P) - \infsum^{\Xi_{h}}(f, P) = \supsum^{\Xi}(f \circ h, Q) - \infsum^{\Xi}(f \circ h, Q) < \varp.
    $$
    Thus, by \autoref{t29}, $f \in \cI(\Xi_{h})$. Finally, the value of the integral follows by applying~\ref{s2.1}. 
\end{PROOF}

\begin{corollary}[{\cite[Cor.~3.18]{PU2025}}]\label{da1}
    Let $X, Y$ be non-empty sets such that $X \subseteq Y$, %$\iota \colon X \to Y$ the inclusion function, $\cB_{X}$, 
    $\cB_{Y}$ Boolean algebras on $\calP(X)$ and $\calP(Y)$, respectively, $\Xi^{X}$, $\Xi^{Y}$ finitely additive measures on $\cB_{X}$  and $\cB_{Y}$, respectively, and $g \colon Y \to \bbR$ a bounded function. Assume $\cB_{Y} \subseteq \set{b\subseteq Y}{b\cap X\in\cB_X}$ and $\Xi^{Y}(b) = \Xi^X(b\cap X)$ for all $b\in\cB_Y$. If $g \in \cI(\Xi^{Y})$ then $g \rest X \in \cI(\Xi^{X})$ and 
    \begin{equation*}\label{e66}
        \text{$\int_{Y} g d \Xi^{Y} = \int_{X} g \rest X d \Xi^{X}.$}
    \end{equation*}
\end{corollary}
\begin{PROOF}{\ref{da1}}
    Apply \autoref{s2} to the inclusion map $\iota\colon X\to Y$ and \autoref{t93}. Notice that $\iota^\to(\cB_X) = \set{b\subseteq Y}{b\cap X\in\cB_X}$ and $\Xi_\iota(b) = \Xi^X(b\cap X)$ for all $b\in \iota^\to(\cB_X)$.
\end{PROOF}

We complete this section by introducing a notion of \emph{measurability} of real-valued functions in the context of fields of sets.

\begin{definition}\label{t32}
    A function $h \colon X \to \bbR$ is \emph{$\cB$-semi-measurable} if it is $\cB$-$\cC^1$-measurable (see \autoref{t16}) i.e.\ for any $y, z \in \bbR$,  $\{ x \in X \colon y \leq h(x) < z  \} \in \cB$. Notice that, if $\cB$ is $\sigma$-complete this is equivalent to $h$ being $\cB$-measurable, that is, for any Borel set $E \subseteq \bbR$, $h^{-1}[E] \in \cB$. 
    %Notice that this is equivalent to say that $\set{x\in X}{h(x)<y}\in\cB$ for all $y\in\bbR$.
\end{definition}

For example, it is easy to see that, if $E \in \cB,$ then $\chi_{E}$ is $\cB$-semi-measurable. %We have chosen the name \emph{$\cB$-semi-measurable}, since, as shown in \autoref{t38.1}, obtaining an analogous measurability in the usual sense of analysis requires the algebra to be $\sigma$-complete.

The notion of $\cB$-semi-measurability give us a useful condition of $\Xi$-integrability. 
 
\begin{theorem}\label{t38}
    If $h \colon X \to \bbR$ is bounded and $\cB$-semi-measurable, then $h$ is $\Xi$-integrable. As a consequence, if $\cB = \calP(X)$, then any bounded function $h \colon X \to \bbR$ is $\Xi$-integrable. 
\end{theorem}

\begin{PROOF}{\ref{t38}}
    %By \autoref{t560}, it is enough to show that $h$ is $\hat\Xi$-integrable.
    Let $\varepsilon > 0.$ Since $h$ is bounded, there are $a, b \in \bbQ$ such that $a < b$ and $f[X] \subseteq [a,b)$. Partition $[a,b)$ into finitely many intervals $\{[a_k,b_k) \colon k < m \}$ with rational endpoints and with length ${<}\frac{\varepsilon}{\Xi(X)+1}$. Consider $P \coloneqq \{h^{-1}\left[[a_k,b_k)\right] \colon  k<m\}$, which is in $\Pbf^{\Xi}$ because $h$ is $\cB$-semi-measurable. Hence,
    \begin{equation*}
        \begin{split}
            \supsum(h,P)-\infsum(h,P) & = \sum_{b\in P}(\sup h[b]-\inf h[b]) \Xi(b)\\
            & \leq \sum_{b\in P}\frac{\varepsilon}{\Xi(X)+1}\Xi(b)=\frac{\varepsilon}{\Xi(X)+1}\Xi(X)<\varepsilon.
        \end{split}
    \end{equation*}    
    Thus, by the criterion of $\Xi$-integrability, $h \in \cI(\Xi).$ 
\end{PROOF}

The converse of the previous theorem is not always true. See \autoref{u20}.

If $\cB$ is a $\sigma$-algebra and $\Xi$ is $\sigma$-additive, then any bounded and $\cB$-measurable $f\colon X\to \R$ is both $\Xi$-integrable and Lebesgue integrable (the latter is well-known), and both integrals of $f$ coincide. One may ask whether $\Xi$-integrable implies $\cB$-measurable (and hence Lebesgue integrable). This is true when the measure space $(X,\cB,\Xi)$ is complete, which we prove in \autoref{t46}, but there are silly counter-examples without completeness: when $\cB$ is the Borel $\sigma$-algebra of $\R$ and $\Xi$ is the probability measure on $\cB$ such that $\Xi(\{0\})=1$, then any bounded $f\colon \R\to\R$ is $\Xi$-integrable and $\int_\R f d\Xi = f(0)$, but there are bounded functions that are not $\cB$ measurable. 

% \blueq{el problema es que se necesitan hipotewsis parar cosas sencillas: por ejemplo, la suma de medibles necesita sigma completes

% \begin{lemma}\label{t41}
%     If $f$ is $\cB$-semi-measurable, if $\cB$ is $\sigma$-contable  then $\forall y \in \bbR (\{ x \in X \colon f(x) < a \}, \{ x \in X \colon f(x) > a\} \in \cB).$
% \end{lemma}

% \begin{PROOF}{\ref{t41}}
%     dfdfdfdf
% \end{PROOF}

% \begin{theorem}\label{t44}
%     If $f$ is $\cB$-semi-measurable and $\cB$ is $\sigma$-complete, then $\{ x \in X \colon f(x) \neq g(x) \} \in \cB.$ 
% \end{theorem}

% \begin{PROOF}{\ref{t44}}
%     fgfg
% \end{PROOF}}

% \blueq{Así que, las funciones de tal parte son integrables en algebras booleanas sigma completas. Sin embargo, este criterio tiene ina debilidad mas alla de la sigma completes. Si iteramos por ejemplo con fams sobre $\calP(\omega),$ entonces pueden aparecer nuevos reales y funciones que antes eram medibles, dejan de serlo en las extensiones genericas. ASi que necesitamos resultados de integrabilidad menos vulnerables antes el terror de las iteraciones.  Para los propositos de este trabajo, solo necesitamos integrar funciones sencillas. Debemos trabajar un poco mas: }

\subsection{Extension and approximation criteria with integrals}\label{7}

We prove some extension and approximation criteria for finitely additive measures involving integrals, generalizing results of \autoref{3} and results presented (without proof) in~\cite[Sec.~1]{Sh00}.

First, we use \autoref{m40} to prove the following lemma:

\begin{lemma}\label{t80}
    Let $\cB$ a Boolean subalgebra of $\pts(X)$,
    $f\colon X\to\R$ a bounded function and let $\delta\in[0,\infty)$. Define $Z \coloneqq {}^{\cB} [0, \delta]$ and $F^\cB_\delta \coloneqq  \set{ \Xi\in Z}{\Xi\text{ is a fam and } \Xi(X)=\delta }$, (which is closed in $Z$ by \autoref{m40}). Then: 
    \begin{enumerate}[label = \rm (\alph*)]
        \item\label{t80a} The map $\Xi\mapsto \overline{\int_{X}} fd\Xi$ is upper semicontinuous on $F^\cB_\delta$. 
        \item\label{t80b} The map $\Xi\mapsto \underline{\int_{X}} fd\Xi$ is lower semicontinuous on $F^\cB_\delta$.
        \item\label{t80c} For any $r\in\R$, the set
          \[\overline{A}(f,r)\coloneqq \set{ \Xi\in F^\cB_\delta}{ \overline{\int_{X}} fd\Xi \geq r }\] is closed in $Z$.
        \item\label{t80d} For any $r\in\R$ the set 
         \[\underline{A}(f,r)\coloneqq \set{ \Xi\in F^\cB_\delta}{ \underline{\int_{X}} fd\Xi \leq r }\] is closed in $Z$.
    \end{enumerate}
\end{lemma}

\begin{PROOF}{\ref{t80}}
We only show~\ref{t80c} (\ref{t80d} can be proved similarly,~\ref{t80a} and~\ref{t80c} are equivalent, and likewise~\ref{t80b} and~\ref{t80d}). 
Suppose that  $\Xi\in \mathrm{cl}_Z(\overline{A}(f,K))$. By \autoref{m40}, $\Xi\in F^\cB_\delta$. It remains to prove that $\overline{\int_{X}} fd\Xi\geq r$. 

It is enough to show that $\sum\limits_{b\in P}\sup(f[b])\,\Xi(b)\geq r$ for any $P\in\Pbf^\Xi$.         
Find $M>0$ such that $|f(x)|<M$ for all $x\in X$. Fix $\varepsilon>0$ and 
let $\varepsilon'\coloneqq \frac{\varepsilon}{M|P|}$. Now, for $b\in P$, set $U_{b}\coloneqq(\Xi(b)-\varepsilon',\,\Xi(b)+\varepsilon')\cap[0,\delta]$. Then we consider the open neighborhood $U\coloneqq \prod\limits_{a\in \cB}U_a,$ of $\Xi$ where $U_a=[0,\delta]$ for each $a\notin P$. Choose some $\Xi' \in\overline{A}(f,r)$ such that $|\Xi(b)-\Xi'(b)|<\varepsilon'$ for every $b\in P$, so $\sum\limits_{b\in P}\sup(f[b])\,\Xi'(b)\geq \overline{\int_X}f d\Xi' \geq r$. Now, observe that 
\begin{equation*}
        \begin{split}
            \left|\sum\limits_{b\in P}\sup(f[b])\Xi(b)-\sum\limits_{b\in P}\sup(f[b])\Xi'(b)\right|&\leq \sum\limits_{b\in P}|\sup(f[b])\,\Xi(b)-\sup(f[b])\,\Xi'(b)|\\
            & = \sum\limits_{b\in P}|\sup(f[b])||\,\Xi(b)-\Xi'(b)|\\
            & \leq\sum\limits_{b\in P}|\sup(f[b])|\varepsilon' \leq \varepsilon.\\
           % & \leq \sum \limits_{b\in P}M\varepsilon'=|P|M\varepsilon'\\
           % & = \varepsilon.
        \end{split}
    \end{equation*}
On the other hand, since $\sum\limits_{b\in P}\sup(f[b])\,\Xi'(b)\geq r$, we obtain \[r-\varepsilon\leq\sum\limits_{b\in P}\sup(f[b])\,\Xi'(b)-\varepsilon\leq \sum\limits_{b\in P}\sup(f[b])\,\Xi(b).\] 
Since $\varepsilon>0$ was taken arbitrary, $\sum\limits_{b\in P}\sup(f[b])\,\Xi(b)\geq r$. %This implies that $\Xi\in A(f,r)$. 
%Since $U$ is an open neighborhood of $\Xi$,    
\end{PROOF}

\begin{corollary}\label{t805}
    Let $\cB$, $f$, $\delta$, $Z$ and $F^\cB_\delta$ be as in \autoref{t80}. If $f$ is $\Xi$-integrable for all $\Xi \in F^\cB_\delta$ then, for any closed $K\subseteq \R$, the set
    \[A(f,K)\coloneqq  \set{\Xi\in F^\cB_\delta}{ \int_X fd\Xi \in K}\]
    is closed in $Z$.
\end{corollary}
\begin{PROOF}{\ref{t805}}
    Immediate by \autoref{t80} because, by assumption, $\Xi \mapsto \int_X fd\Xi$ is a continuous map on $F^\cB_\delta$.
\end{PROOF}

\begin{example}
    In the context of the previous results, the set of fams $\Xi$ such that $f$ is $\Xi$-integrable may not be closed in ${}^\cB [0,\delta]$. Consider $X=\omega$, $\delta = 1$ and let $\cB$ be the Boolean subalgebra of $\pts(\omega)$ composed of all the finite and co-finite sets. Note that the map $\Xi\mapsto \la \Xi(\{n\})\colon n<\omega\ra$ is a bijection from $F^\cB_1$ onto $\set{x\in {}^\omega [0,1]}{\sum_{k<\omega} x_k \leq 1}$, and it is continuous and hence a homeomorphism (because $F^\cB_1$ is compact).

    Any finite partition of $\omega$ into sets in $\cB$ can be refined to a partition of the form 
    \[P_n\coloneqq \set{\{k\}}{k<n}\cup \{\omega\menos n\},\] 
    and $P_{n+1}$ refines $P_n$. Therefore, for any bounded $f\colon \omega\to \R$ and any $\Xi\in F^\cB_1$,
    \begin{align*}
        \overline{\int_X}f d\Xi & = \lim_{n\to \infty}\supsum(f,P_n) \text{ and }\\
        \underline{\int_X}f d\Xi & = \lim_{n\to \infty}\infsum(f,P_n).
    \end{align*}
    Also, by setting $x_k\coloneqq \Xi(\{k\})$,
    \begin{align*}
        \supsum(f,P_n) & = \sum_{k<n}f(k) x_k + \sup_{k\geq n}f(k) \left(1-\sum_{k< n}x_k\right), \text{ and }\\
        \infsum(f,P_n) & = \sum_{k<n}f(k) x_k + \inf_{k\geq n}f(k) \left(1-\sum_{k< n}x_k\right),
    \end{align*}
    so
    \begin{align*}
        \overline{\int_X}f d\Xi & = \sum_{k<\omega}f(k) x_k + \limsup_{k\to \infty}f(k) \left(1-\sum_{k< \omega}x_k\right), \text{ and }\\
        \underline{\int_X}f d\Xi & = \sum_{k<\omega}f(k) x_k + \liminf_{k\to \infty}f(k) \left(1-\sum_{k< \omega}x_k\right).
    \end{align*}
    Therefore, $f$ is $\Xi$-integrable iff either $\lim_{k\to\infty} f(k)$ exists or $\sum_{k<\omega}x_k = 1$. In the first case, $f$ is $\Xi$-integrable for any $\Xi\in F^\cB_1$, but in the case that $\lim_{k\to\infty} f(k)$ does not exist, the set of $\Xi\in F^\cB_1$ such that $f$ is $\Xi$-integrable is not closed. The reason is that its corresponding set via the previously mentioned homeomorphism is $\set{x\in {}^\omega [0,1]}{\sum_{k<\omega} x_k = 1}$, and this is not closed in ${}^\omega [0,1]$. For example, let $x^n\in {}^\omega [0,1]$ defined by
    \[
      x^n_k \coloneqq \left\{
         \begin{array}{ll}
            \frac{1}{n+1}  &  \text{if $k\leq n$,}\\
            0  & \text{if $k>n$.}
         \end{array}
      \right.
    \]
    Although $\sum_{k<\omega}x^n_k =1$, $\lim_{n\to\infty} x^n$ (in ${}^\omega [0,1]$) is the constant sequence $0$.
\end{example}

The following result gives us a criterion to approximate integrals.

%{\color{blue} Incluso con $F=\emptyset$, el siguiente teorema no es cierto en general cuando $\Xi$ no es libre. Considere el mismo contraejemplo presentado luego de \autoref{famgen} con $E=\omega$, $i^*=1$ y $f_0(x)=x+1$: la integral es $\frac{5}{3}$, y cuando $|u|=k>0$, la suma es mayor que $2$ cuando $u\neq k$, y cuando $u=k$ la suma es igual a $\frac{k+1}{2}$.}

\begin{lemma}\label{t90}
    Let $\Xi$ be a fam on a Boolean subalgebra $\cB$ of $\calP(X)$ such that $\Xi(X) = \delta \in (0,\infty)$ and let 
    $i^{\ast} < \omega$. For any $i < i^{\ast},$ let $f_{i} \colon X \to \bbR$ be a bounded function. Then, for all $\varp > 0$ and any partition $P^*\in \Pbf^\Xi $, there exists a finite set $u \subseteq X$ of size ${\leq}\left\lceil\frac{\delta}{\varp}\right\rceil$ and a probability measure on $\pts(u)$ such that, 
    \begin{enumerate}[label = \rm (\alph*)]
        \item\label{t90a} for any $b\in P^*$, $
        |\delta \mu(u \cap b) - \Xi(b)| < \varp$, and
        \item\label{t90b} for any $i < i^{\ast}$, $$\underline{\int_{X}} f_{i} d \Xi -\varp  < \delta\int_u f_i d\mu = \delta \sum_{x\in u}f_i(x) \mu(\{x\}) < \overline{\int_{X}} f_{i} d \Xi + \varp $$    
    \end{enumerate}
    Moreover, when $\Xi$ has the uap and $\varp$ is small enough, $\mu$ in~\ref{t90a} can be found uniform (so $\int_u f_i d\mu = \frac{1}{|u|}\sum_{k\in u}f_i(x)$) and, whenever $\Xi(b) = 0$ for any finite $b\in\cB$, we can even find $u\subseteq X\menos F$ when fixing an arbitrary finite $F\subseteq X$.
\end{lemma}

%\frac{\delta}{\vert u \vert}  \sum_{x \in u} f_{i}(x)

\begin{PROOF}{\ref{t90}}
Find an $0<M < \omega$ such that, for any $i < i^{\ast}, \, \vert f_{i}(x) \vert < M$ for all $x \in X.$  Let $\varp > 0$ and $P^* \in \Pbf^\Xi$. For each $i < i^{\ast}$, there is some $P^{i} \in \Pbf^\Xi$ refining $P^*$ such that 
\[ \underline{\int_{X}} f_{i} d \Xi - \frac{\varp}{2} < \infsum(f_i,P^i) \text{ and } \supsum(f_i,P^i) < \overline{\int_{X}} f_{i} d \Xi + \frac{\varp}{2}.\]

Let $P \coloneqq \sqcap_{i < i^{\ast}} P^{i}.$ By \autoref{t18}, for each $i < i^{\ast},$ it follows that 
$$\underline{\int_{X}} f_{i} d \Xi - \frac{\varp}{2} < \infsum(f_i,P) \text{ and } \supsum(f_i,P) < \overline{\int_{X}} f_{i} d \Xi + \frac{\varp}{2}.$$ 
By \autoref{m421}, we can find some non-empty finite $u\subseteq X$ of size ${\leq}\min\left\{|P|,\left\lceil\frac{\delta}{\varp}\right\rceil\right\}$ and a probability measure on $\pts(u)$ such that, for any $b\in P$,
\[\left|\delta\mu(u\cap b) - \Xi(b)\right|<\frac{\varp}{4 M |P|}.\]
Since $P$ is a refinement of $P^*$, we can conclude~\ref{t90a}. 

For any $i < i^{\ast}$, let $g_{i} \coloneqq f_{i} + M$. Then, 
 \begin{equation*}
        \begin{split}
            % line 1
            \delta\int_u f_i d\mu +\delta M  & =  \int_u g_{i} d\mu  = \sum_{b\in P} \delta\int_{u\cap b} g_{i} d\mu   
             \leq  \sum_{b\in P}  \sup(g_{i}[b])\, \delta \mu(u\cap b)\\
            %line 2
            & \leq  \sum_{b\in P} \left[ \sup(g_{i}[b])  \left ( \Xi(b) + \frac{\varp}{4 M |P|} \right) \right]  = \supsum(g_i,P) + \sum_{b \in P} \sup(g_{i}[b]) \frac{\varp}{4 M |P|}\\
            & \leq \supsum(g_i,P) + \frac{\varp}{2} = \supsum(f_{i}, P) + \delta M + \frac{\varp}{2} < \overline{\int_{X}} f_{i} d \Xi + \varp + \delta M. 
        \end{split}
    \end{equation*}
    As a consequence, 
    $$
    \delta\int_u f_i d\mu < \overline{\int_{X}} f_{i} d \Xi + \varp.
    $$
    %
    % \begin{equation*}
    %     \begin{split}
    %         % line 1
    %         \sum_{x \in u} \frac{ \delta \, f_{i}(x)}{\vert u \vert} & = \sum_{b\in P} \left( \sum_{x \in u\cap b} \delta \frac{f_{i}(x)}{|u|} \right)  
    %          \leq  \sum_{b\in P}  \sup(f_{i}[b])\, \delta\,  \frac{|u\cap b|}{|u|}\\
    %         %line 2
    %         & \leq  \sum_{b\in P} \left[ \sup(f_{i}[b])  \left ( \Xi(b) + \frac{\varp}{2 M |P|} \right) \right] \\
    %         % line 3
    %         & = \supsum(f_i,P) + \sum_{b \in P} \sup(f_{i}[b]) \frac{\varp}{2 M |P|} \leq \supsum(f_i,P) + \frac{\varp}{2} 
    %         %line 4
    %         %& \leq \frac{\delta}{\Xi(E)} \supsum(f_{i}, P) + \sum_{j < j^{\ast}} \frac{\varp}{2j^{\ast}} \\
    %         %line 5
    %         %& < \frac{1}{\Xi(E)} \left( \frac{\varp\Xi(E)}{2} + \delta \infsum(f_{i}, P) \right) + \frac{\varp j^{\ast}}{2j^{\ast}}\\
    %         %line 6
    %         %& \leq \frac{1}{\Xi(E)} \left( \frac{\varp \Xi(E)}{2} +  \delta \int_{E} f_{i} d \Xi \right) + \frac{\varp}{2} \\
    %         % line 7
    %         %& = \frac{\varp}{2} + \frac{\delta}{\Xi(E)} \int_{E} f_{i} d \Xi + \frac{\varp}{2} 
    %         % line 8
    %     < \overline{\int} f_{i} d \Xi + \varp. 
    %     \end{split}
    % \end{equation*}
    Similarly, we obtain
    \[\delta\int_u f_i d\mu > \underline{\int_{X}} f_{i} d \Xi - \varp.\]
    For the latter part of the theorem, use \autoref{m401} and \autoref{m420}.
\end{PROOF}

%Notice that this is an approximation result in the following sense: by \autoref{t42}, if $u$ is the finite set from the conclusion of \autoref{t90}, then $\int_{X} f_{i} d \Xi^{u} = \frac{1}{\vert u \vert} \sum_{x \in u} f_{i}(x),$ so this allows us to approximate $\int_{X} f_{i} d \Xi$ using the integral with respect to $\delta \Xi^{u}.$  

Assume that $I$ is an index set and $f_{i} \colon X \to \bbR$ is a function for any $i \in I$. In the following result we use \autoref{t90} to prove an extension criterion similar to \autoref{famgen-} and~\ref{famgen}, but including some control on the value of the integrals of the functions $f_{i}$. 

\begin{theorem}\label{t85}
   Let $\Xi_{0}$ be a finitely additive measure on some Boolean subalgebra $\cB \subseteq \mathcal{P}(X)$ and let $\delta \coloneqq \Xi_0(X)\in (0,\infty)$.
   Let $I$ be an index set and for each $i\in I$, let $K_i$ be a closed subset of $\R$ and $f_i\colon X\to\R$ bounded.  
   Then the following statements are equivalent.
   \begin{enumerate}[label = \rm (\roman*)]
   \item\label{t85I} For any $P\in\Pbf^{\Xi_0}$, $\varepsilon>0$, any finite set $J\subseteq I$, and any open $G_i\subseteq \R$ containing $K_i$ for $i\in J$, there is some non-empty finite $u\subseteq X$  and a probability measure on $\pts(u)$ such that:
   \begin{enumerate}[label=\normalfont(\alph*)]
       \item\label{t85a} $ \left |\Xi_0(b)-\delta\mu(b\cap u) \right | <\varepsilon$ for any $b\in P$, and
       
       \item\label{t85b} $\delta\int_u f_i d\mu \in  G_i$ for any $i\in J$.
   \end{enumerate}
   \item\label{t85II} There is some fam $\Xi$ on $\pts(X)$ extending $\Xi_0$ such that, for any $i \in I$, $\int_{X} f_id \Xi \in K_i$.
   \end{enumerate}
   Even more, $u$ can be found of size ${\leq}\left\lceil\frac{\delta}{\varp}\right\rceil$. 

   Also, the equivalence can be modified so that $\mu$ is uniform in~\ref{t85I} and $\Xi$ has the uap in $\ref{t85II}$, but the upper bound of the size of $u$ holds for small enough $\varp$. In addition, $\Xi$ can be restricted to be free in~\ref{t85II} as long as $u$ in~\ref{t85I} can be found disjoint with any fixed (but arbitrary) finite $F\subseteq X$. 
\end{theorem}

\begin{PROOF}{\ref{t85}}
    \ref{t85I}${}\Rightarrow{}$\ref{t85II}: 
    Let $Z\coloneqq  {}^{\mathcal{P}(X)}[0,\delta]$.
    Given $J\subseteq I$, $P\in\Pbf^{\Xi_0}$ and $\varepsilon\geq0$, for $i\in J$ let $K_i^\varp$ be the set of points in $\R$ with (euclidean) distance ${\leq}\varp$ to $K_i$ (so $K^0_i = K_i$), and define 
    \begin{multline*}
        F_{J,P,\varepsilon} \coloneqq  \bigg\{\Xi\in Z \colon  \Xi \text{\ is a fam, } \Xi(X) = \delta,\ \forall b\in P \ ( |\Xi(b)-\Xi_0(b)|\leq \varepsilon ) \text{ and }\\
       \forall i\in J\ \bigg(  \int_X f_id\Xi \in K_i^\varp \bigg)  \bigg\}.
    \end{multline*}

    If $u$ and $\mu$ satisfy the conditions~\ref{t85a} and~\ref{t85b} (using the interior of $K^\varp_i$), then $\delta \mu\in F_{J,P,\,\varepsilon}$ for $\varp>0$ and finite $J\subseteq I$. By \autoref{t805} and \autoref{m40}
    $F_{J,P,\varepsilon}$ is a compact set. Even more, 
    %Notice that $\varepsilon'\leq\varepsilon$ implies that $F_{J,P,\varepsilon'}\subseteq F_{J,P,\varepsilon}$, so we get that 
    $\{F_{J,P,\varepsilon}\colon \varepsilon>0\}$ has the finite intersection property. Then $F_{J,P,0}$ is non-empty and closed by compactness. 

    It is clear that $J\subseteq J'$ and $P'$ is a refinement of $P$ in $\Pbf^{\Xi_0}$ imply $F_{J',P',0}\subseteq F_{J,P,0}$. Then the family $\{F_{P,J,0} \colon P\in\Pbf^{\Xi_0},\ J\in[I]^{<\aleph_0}\}$ has the finite intersection property, so there is some $\Xi$ in its intersection.

     \ref{t85II}${}\Rightarrow{}$\ref{t85I}: Let $P$, $\varp$, $J$ and $G_i$ as in~\ref{t85I}. Wlog we may assume that $\varp$ is small enough such that the interval with center in $\int_X f_i d\Xi$ and length $2\varp$ is contained in $G_i$ for any $i\in J$.  By \autoref{t90} we can find the required $u$ and $\mu$ satisfying~\ref{t85a} and $\int_X f_{i} d \Xi -\varp  < \delta \int_u f_{i} d\mu < \int_X f_{i} d \Xi + \varp$ for all $i\in J$. Hence $\delta \int_u f_{i} d\mu \in G_i$.

     The rest follows the same proof, using \autoref{uapcl}.
\end{PROOF}

In the previous result, even if $\Xi_0(b) = 0$ for all finite $b\in \cB$, it cannot be guaranteed that $\Xi$ can be found free without the requirement of the finite $F$, for example, when $\cB = \{\emptyset,X\}$ with $0\in X$, $f_0 = \chi_{\{0\}}$ and $K_0 = \{\delta\}$.

When $\Xi_0$ is a fam associated with an ultrafilter on $\cB$, we may not obtain a fam $\Xi$ associated with an ultrafilter on $\pts(X)$. To guarantee extension to an ultrafilter, the theorem is modified as follows.

\begin{theorem}\label{t851}
   Let $U_0$ be an ultrafilter on some Boolean subalgebra $\cB \subseteq \mathcal{P}(X)$.
   Let $I$ be an index set and for each $i\in I$, let $K_i$ be a closed subset of $\R$ and $f_i\colon X\to\R$ bounded.  
   Then the following statements are equivalent.
   \begin{enumerate}[label = \rm (\roman*)]
   \item\label{t851i} For any $b\in U_0$, any finite set $J\subseteq I$, and any open $G_i\subseteq \R$ containing $K_i$ for $i\in J$, there is some $z\in b$ such that $f_i(z)\in G_i$ for all $i\in I$.
   
   \item\label{t851ii} There is some ultrafilter $U\subseteq\pts(X)$ containing $U_0$ such that, for any $i \in I$, $\lim_U f \in K_i$.
   \end{enumerate}
\end{theorem}
\begin{PROOF}{\ref{t851}}
    Note that~\ref{t851i} implies that the set
    \[U_0 \cup \bigcup_{i\in I}\set{f_i^{-1}[G]}{G\supseteq K_i \text{ open}}\]
    has the finite intersection property, so there is an ultrafilter $U$ on $\pts(X)$ containing it. Fix $i\in I$. If $\lim_U f_i \notin K_i$, then we can find disjoint open sets $G$ and $G'$ such that $K_i\subseteq G$ and $\lim_U f_i\in G'$, so $f_i^{-1}[G] \in U$ and $f_i^{-1}[G'] \in U$ (the latter by the definition of ultrafilter limit), hence $f_i^{-1}[G\cap G']\in U$, which implies that $G\cap G'\neq \emptyset$, a contradiction. Therefore, $\lim_U f_i \in K_i$.

    For the converse, assume~\ref{t851ii}, $b\in U_0$, $J\subseteq I$ finite and $G_i\subseteq\R$ open containing $K_i$ for all $i\in J$. Since $\lim_U f_i\in G_i$, $f_i^{-1}[G_i] \in U$, so $b\cap \bigcap_{i\in J}f_i^{-1}[G_i] \in U$, so any $z$ in this set is as required for~\ref{t851i}.
\end{PROOF}

%%%%%%%%%%%%%%%%%%%%%%%%%%%%%%%%%%%%%%%%%%%%%%%%%%%%%%%%%%%%%%%%%%%%%%%%%%%%%%%%%%%%%%%%%%%%%%%%%%%%%%%%%%%%%%%%%%%%%%%%%%%%%%%%%%%%%%%%%%%%%%%%%%%

\section{The Jordan measure}\label{5}

Fix a Boolean subalgebra $\cB$ of $\calP(X)$ for some non-empty set $X$ and a finitely additive measure $\Xi \colon \cB \to [0, \delta],$ where $\delta$ is a non-negative real number. In this section, we are going to define expansions $\cJ^{\Xi}$ and $\hat{\Xi}$ of $\cB$ and $\Xi$, known as the \emph{Jordan algebra} and the \emph{Jordan measure}, respectively. Additionally, we prove a limit theorem that will allow us to establish the main result in this section: the equivalence between $\Xi$-integrability and $\hat{\Xi}$-integrability.

We start looking at integration over subsets of $X$. For $E\subseteq X$, recall that $\chi_E\colon X\to\{0,1\}$ denotes the \emph{characteristic function of $E$ (over $X$)}, i.e.\ $\chi_E(x)=1$ iff $x\in E$.

\begin{definition}\label{t67.e}
    Let $E \subseteq X$ and $f\colon X\to\bbR$ bounded on $E$. We say that \emph{$f$ is $\Xi$-integrable on $E$} when  $f \chi_{E} \in \cI(\Xi)$, in which case we denote 
    $$\int_{E} f d \Xi \coloneqq \int_{X} f \chi_{E}  d \Xi.$$
    This definition makes sense when the domain of $f$ contains $E$ but is not necessarily all $X$. In this case, it is natural to consider $f\chi_E$ with domain the whole $X$ (points in $X\menos \dom f$ are sent to $0$).
\end{definition}

It is not hard to prove that we have a monotonicity property over subsets. 

\begin{lemma}\label{t70}
    Let $E, F\subseteq X$ and assume that $f\colon X\to \bbR$ is bounded, non-negative, and $\Xi$-integrable on $E$ and $F$. If $E \subseteq F,$ then $\int_{E} f d \Xi \leq \int_{F} f d \Xi.$
\end{lemma} 

\begin{PROOF}{\ref{t70}}\
    If $E \subseteq F$ then $\chi_{E} f \leq \chi_{F}f$, so the result follows by \autoref{t53}. %we have that: \[ \int_{E} f d \Xi = \int_{X} \chi_{E} f d \Xi \leq \int_{X} \chi_{F} f d \Xi = \int_{F} f d \Xi.\qedhere\] 
\end{PROOF}

\begin{lemma}\label{t31}
    Let $E\in\cB$ and $f\colon X\to \bbR$ bounded on $E$. Then $f{\upharpoonright}E \in \cI(\Xi|_{E})$ iff $f\chi_E\in\cI(\Xi)$. In addition, when $f$ is $\Xi$-integrable on $E$,
    \[\int_E f d \Xi = \int_E f{\upharpoonright}Ed \Xi|_{E}.\]
\end{lemma}

To prove this and more results, we use the following resource.

\begin{definition}\label{t62}
    Let $E \, {\in} \, \cB.$ For $P \in \bfP^{\Xi}$ we define $P_{E} \coloneqq \{ E \cap A \colon A \in P \}$ and $\hat{P}_{E} \coloneqq P \sqcap \{ E, X\menos E \}.$ We call $P_{E}$  the \emph{partition of $E$ induced by $P$}.  \index{induced partition} \index{$P_{E}$} \index{$\hat{P}_{E}$}
\end{definition}

It is clear that, for every $E\in\cB$ and $P \in \bfP^{\Xi}, \, P_{E}\in \bfP^{\Xi|_E}$ and $\hat{P}_{E} \in \bfP^{\Xi}.$   Also, $P_{E} \subseteq \hat{P}_{E}.$ 

Using \autoref{t62} we can prove \autoref{t31}. 

\begin{PROOF}[Proof of \autoref{t31}]{\ref{t31}}
   Assume that $f$ is $\Xi$-integrable on $E$. Let $\varp>0$. By \autoref{t29}~\ref{t29iii} with $L\coloneqq \int_E f d\Xi$, there is some $P\in \bfP^\Xi$ such that, for any $Q\ll P$ in $\bfP^\Xi$,
   \[|\Sm^\Xi(f\chi_E,\tau)-L|<\varp \text{ for any }\tau\in \prod Q.\]
   Let $Q\coloneqq \hat{P}_E$.
   Then, for any $\sigma\in\prod P_E$, if we pick some $\tau\in\hat{P}_e$ extending $\sigma$, then 
   $$S^{\Xi|_E}(f\frestr E,\sigma) = S^\Xi(f\chi_E,\tau),$$
   so $|S^{\Xi|_E}(f\frestr E,\sigma)-L|<\varp$. Therefore, by \autoref{t29}~\ref{t29iv}, $f\frestr E$ is $\Xi|_E$-integrable and $\int_E f{\upharpoonright}Ed \Xi|_{E}= L =\int_E f d\Xi$.

   The proof of the converse follows a similar idea: if $f$ is $\Xi|_E$-integrable then, for $\varp>0$, there is some $P\in \bfP^{\Xi|_E}$ as in \autoref{t29}~\ref{t29iv}. So the partition $P'\coloneqq P\cup\{X\menos E\}$ can be used to show that $f$ is $\Xi$-integrable on $E$.
\end{PROOF}

We now deal with characteristics functions. 

\begin{lemma}\label{t63}
Let $E\subseteq X$. Then:
\begin{enumerate}[label = \normalfont(\alph*)]
    \item\label{t63.a} $\overline{\int_X} \chi_Ed \Xi = \inf\{\Xi(b)\colon E\subseteq b,\ b\in\cB\}$.

    \item\label{t63.b} $\underline{\int_X} \chi_Ed \Xi = \sup\{\Xi(a)\colon a\subseteq E,\ a\in\cB\}$.

    \item\label{t63.c} If $E \in \cB,$ then $\chi_{E} \in \cI(\Xi)$ and $\int_{X} \chi_{E} d \Xi = \Xi(E)$.
\end{enumerate}    
\end{lemma}

\begin{PROOF}{\ref{t63}}
    \ref{t63.a}: Let $b \in \cB$ such that $E \subseteq b$ and consider $P \coloneqq \{ b, X \setminus b \}$, which is a partition of $X$. Then, 
    $$ \supsum(\chi_{E}, P) = \sum_{b' \in P} \sup(\chi_{E}[b']) \Xi(b') = \sup(\chi_{E}[X \setminus b]) \Xi(X \setminus b)  + \sup(\chi_{E}[b]) \Xi(b) \leq \Xi(b),$$
    
    and therefore, $\overline{\int_{X}} \chi_{E} d \Xi \leq \Xi(b)$. As a consequence, $ \overline{\int_{X}} \chi_{E} d \Xi \leq \inf\{\Xi(b)\colon E\subseteq b,\ b\in\cB\}.$ 

    For the converse inequality, let $P \in \bfP^{\Xi}$. Then, $$ \supsum(\chi_{E}, P) = \sum_{b \in P} \sup(\chi_{E}[b]) \Xi(b) = \sum_{b \in P, \  b \cap E \neq \emptyset} \Xi(b) = \Xi(b_{E}), $$ where $b_{E}$ is the union of all elements $b \in P$ such that $b \cap E \neq \emptyset,$ which is in $\cB$ because $P$ is a finite partition contained in $\cB$. It is clear that $E\subseteq b_E$. 
    As a consequence, $\supsum(\chi_{E}, P) \geq \inf\{\Xi(b)\colon E\subseteq b,\ b\in\cB\}$ and, therefore, $\inf\{\Xi(b)\colon E\subseteq b,\ b\in\cB\} \leq \overline{\int_{X}} \chi_{E} d \Xi.$  

    \ref{t63.b}: Easy by applying~\ref{t63.a} to $X\menos E$.

    \ref{t63.c}: Let $P \in \bfP^{\Xi}$ and notice that  
    $$\supsum(\chi_{E}, \hat{P}_{E}) = \sum_{b \in \hat{P}_{E}} \sup(\chi_{E}[b]) \Xi(b) = \sum_{b \in P_{E}} \Xi(b) = \Xi(E),$$ 
    where $P_{E}$ is the partition of $E$ induced by $P.$ Similarly, \, $\infsum(\chi_{E}, \hat{P}_{E}) = \Xi(E).$ Thus by the criterion of $\Xi$-integrability, $\chi_{E} \in \cI(\Xi)$ and $\int_{X} \chi_{E} d \Xi = \Xi(E).$ 
\end{PROOF}

Notice that, if $E \in \cB$ and $f \in \cI(\Xi)$ then, by \autoref{t63}, $\chi_{E}$ is $\Xi$-integrable and therefore, by \autoref{t58}, $\chi_{E} f$ is $\Xi$-integrable. %In general, this is the context in which we will use \autoref{t67.e}. 

The first two items in \autoref{t63} have a special role as inner and outer measures. So we adopt the following notation.

\begin{definition}\label{t311}
    The map $\Xi^*\colon \calP(X)\to [0,\delta]$ defined by 
    \[\Xi^*(E)\coloneqq \overline{\int_X} \chi_Ed \Xi = \inf\{\Xi(b)\colon E\subseteq b,\ b\in\cB\}\]
    is called the \emph{outer measure of $\Xi$}. Dually, 
    the map $\Xi_*\colon \calP(X)\to [0,\delta]$ defined by
    \[\Xi_*(E)\coloneqq \underline{\int_X} \chi_Ed \Xi = \sup\{\Xi(a)\colon a\subseteq E,\ a\in\cB\}\]
    is called the \emph{inner measure of $\Xi$}.

    We say that a set $E\subseteq X$ is \emph{Jordan $\Xi$-measurable} if $\Xi^*(E)=\Xi_*(E)$. We denote by $\cJ^\Xi$ the collection of Jordan $\Xi$-measurable subsets of $X$, and $\hat{\Xi}\colon \cJ^\Xi\to[0,\delta]$, defined by $\hat{\Xi}(E)=\Xi^*(E)=\Xi_*(E)$, is called the \emph{Jordan $\Xi$-measure}. 
\end{definition}

% The notion of Jordan $\Xi$-measurable generalizes the notion of Jordan measurability for bounded subsets of $\bbR^n$, which corresponds to Jordan $\Xi^n|_{[a,b]}$-measurability (see~\autoref{t16}).

The notion of Jordan $\Xi$-measurability generalizes Jordan measurability for bounded subsets of $\mathbb{R}^n$, which corresponds to Jordan $\Xi^n|_{[a,b]}$-measurability (see~\autoref{t16}).

% \begin{lemma}\label{t317}
%     Let $\langle E_{n} \colon n < \omega \rangle$ be a sequence of subsets of $X$ and $E \coloneqq \bigcup_{n < \omega} E_{n}.$ If for any $n < \omega$, $\Xi^{\ast}(E_{n}) = 0,$ then $E \in \cJ^{\Xi}$ and $\hat{\Xi}(E) = 0.$  
% \end{lemma}

As a consequence of \autoref{t63}:

\begin{corollary}\label{t510}
    $\cB\subseteq\cJ^\Xi$ and $\hat\Xi$ extends $\Xi$. Even more, $E\subseteq X$ is Jordan $\Xi$-measureable iff $\chi_E$ is $\Xi$-integrable, in which case $\hat{\Xi}(E)=\int_X \chi_Ed \Xi$.
\end{corollary}

It is well-known that in the case of $\bbR$, $A$ is Jordan-measurable iff $\lambda(\partial A) = 0,$ where $\partial A \coloneqq \bar{A} \cap A^{\circ}$ is the boundary of $A$ and $\lambda$ denotes the Lebesgue measure on $\R$.  Let $\calB(\bbR)$ the $\sigma$-algebra of Borel sets in $\bbR$. We can find a measurable set $A \subseteq \bfC$ such that $A \notin \calB(\bbR)$, where $\bfC$ is the ternary Cantor set (see e.g.~\cite{Jones}). Since $\bfC$ is closed, $\lambda(\bar{A}) = 0$ and clearly $A^{\circ} = \emptyset$, therefore, $\lambda(\partial A) = 0$, that is, $A \in \cJ^{\lambda}$ and $A \notin \calB(\bbR)$, where $\lambda$ is the Lebesgue measure on $\bbR$.  

We can use elements in $\cB$ to approximate Jordan-measurable sets in the following way.

\begin{lemma}\label{t509}
    Let $E \subseteq X$. Then, $E \in \cJ^{\Xi}$ iff for any $\varp > 0,$ there are $A, B \in \cB$ such that $A \subseteq E \subseteq B$ and $\Xi(B \setminus A) < \varp.$ 
\end{lemma}

\begin{PROOF}{\ref{t509}}
    Assume that $E \in \cJ^{\Xi}$ and let $\varp > 0$. Using \autoref{t63}, we can find $A, B \in \cB$ such that $\hat{\Xi}(E) + \frac{\varp}{2} > \Xi(B)$, $\hat{\Xi}(E) - \frac{\varp}{2} < \Xi(A)$ and $A \subseteq E \subseteq B$. Since $A \subseteq B$, we have that: $$ \Xi(B \setminus A) = \Xi(B) - \Xi(A) < \hat{\Xi}(E) + \frac{\varp}{2} - \hat{\Xi}(E) + \frac{\varp}{2} = \varp.$$ 

    To see the converse, let $\varp > 0$ and $A, B \in \cB$ such that $A \subseteq E \subseteq B$ and $\Xi(B \setminus A) < \varp.$ Let us show that $E \in \cJ^{\Xi}$. By applying \autoref{t63}, it follows that $\overline{\int_{X}} \chi_{E} d \Xi \geq \Xi(B)$ and $\underline{\int}_{X} \chi_{E} d \Xi \leq \Xi(A)$. Hence, $$ \varp > \Xi(B \setminus A) = \Xi(B) - \Xi(A) \geq \overline{\int_{X}} \chi_{E} d \Xi - \underline{\int_{X}} \chi_{E} d \Xi,$$ therefore, $$ \overline{\int_{X}} \chi_{E} d \Xi \leq \underline{\int_{X}} \chi_{E} d \Xi + \varp.$$ Since $\varp$ is arbitrary, we can conclude that $\chi_{E} \in \cI(\Xi)$, that is, $E \in \cJ^{\Xi}$. 
\end{PROOF}

% \Andres{Si definimos $\cA^{\Xi}$ como la colección de los subconjuntos de $X$ que satisfacen la propiedad de \auotref{509}, tenemos que $\cB \subseteq \cJ^{\Xi} \subseteq \cA^{\Xi}$.}

We use the following easy claim to prove that $\cJ^\Xi$ is a Boolean algebra and that $\hat{\Xi}$ is indeed a fam.

\begin{lemma}\label{t73}
    %If $f \in \cI(\Xi)$ and $\langle E_{i} \colon i < n \rangle \in \bfP^{\Xi},$ then:  $$\int_{X} f d \Xi = \sum_{i < n} \left( \int_{E_{i}} f d \Xi \right).$$ 
    If $E_1,E_2\in\cJ^\Xi$, $E_1\cap E_2=\emptyset$, and $f\colon X\to \bbR$ is bounded on $E_1\cup E_2$, then $f$ is $\Xi$-integrable on $E_1\cup E_2$ iff $f$ is $\Xi$-integrable on both $E_1$ and $E_2$, in which case
    $$\int_{E_1\cup E_2} f d \Xi = \int_{E_1} f d \Xi + \int_{E_2} f d \Xi.$$
\end{lemma}

\begin{PROOF}{\ref{t73}}
    On the one hand, assume that $f$ is $\Xi$-integrable on $E_{1} \cup E_{2}$, that is, $\chi_{E_{1}}f + \chi_{E_{2}}f$ is $\Xi$-integrable. Since $E_{1} \cap E_{2} = \emptyset$, we can write $\chi_{E_{1}} f = \chi_{E_{1}}(\chi_{E_{1}}f + \chi_{E_{2}}f)$, which is $\Xi$-integrable by virtue of \autoref{t58} and the hypothesis. Similarly, $\chi_{E_{2}}f \in \cI(\Xi)$. On the other hand, assume that $f$ is $\Xi$-integrable on $E_{1}$ and $E_{2}$. Since $E_{1}$ and $E_{2}$ are disjoint, we can write $\chi_{E_{1} \cup E_{2}}f =  \chi_{E_{1}}f + \chi_{E_{2}}f$, which is $\Xi$-integrable by \autoref{t47}. Finally,
    \begin{equation*}
        \begin{split}
        \int_{E_{1} \cup E_{2}} f d \Xi & = \int_{X} \chi_{E_{1} \cup E_{2}} f d \Xi = \int_{X} ( \chi_{E_{1}} + \chi_{ E_{2}}f) d \Xi \\
        & = \int_{X} \chi_{E_{1}}f d \Xi + \int_{X} \chi_{E_{2}} f d \Xi = \int_{E_{1}} f d \Xi + \int_{E_{2}} f d \Xi.\qedhere 
        \end{split}
    \end{equation*}             
\end{PROOF}

As a consequence, $\cJ^{\Xi}$ is a Boolean algebra extending $\cB$ and $\hat{\Xi}$ is a fam. 

\begin{theorem}\label{t540}
The family $\cJ^\Xi$ is a Boolean subalgebra of $\calP(X)$ and $\hat{\Xi}$ is a fam on $\cJ^\Xi$.
\end{theorem}
\begin{PROOF}{\ref{t540}}
    Since $\cJ^{\Xi}$ is endowed with the set-theoretic operations, to prove that it is a Boolean algebra, it is enough to show that it is closed under such operations. For this, let $E_{0}, E_{1} \in \cJ^{\Xi}$. By virtue of \autoref{t47}, it is clear that $X\menos E_{0} \in \cJ^{\Xi}$. On the other hand, to prove that $E_{0} \cap E_{1}, E_{0} \cup E_{1} \in \cJ^{\Xi}$, we will use \autoref{t509}. We are going to deal with the intersection. The union is completely analogous. 

    Let $\varp > 0.$ By \autoref{t509}, we can find $A_{0}, A_{1}, B_{0}, B_{1} \in \cB$ such that $A_{0} \subseteq E_{0} \subseteq B_{0}$, $A_{1} \subseteq E_{1} \subseteq B_{1}$, $\Xi(B_{0} \setminus A_{0}) < \frac{\varp}{2}$ and $\Xi(B_{1} \setminus B_{0}) < \frac{\varp}{2}$. Define $A \coloneqq A_{0} \cap A_{1}$, $B \coloneqq B_{0} \cap B_{1}$. It is clear that $A \subseteq E_{0} \cap E_{1} \subseteq B$. Also, $B \setminus A \subseteq (B_{1} \setminus A_{1}) \cup (B_{0} \setminus A_{0})$, so $\Xi(B\menos A)<\varp$.  
    %Finally, $A \cap B \in \cJ^{\Xi}.$ 

    As a direct consequence of \autoref{t73}, $\hat\Xi$ is a fam on $\cJ^\Xi$.
\end{PROOF}

We now aim to study the integration theory of $\hat\Xi$. The conclusion will be that it does not differ from the integration over $\Xi$. 

\begin{lemma}\label{t520.0}
    Let $E \subseteq X$. Then, $\overline{\int_{X}} \chi_{E} d \Xi = \overline{\int_{X}} \chi_{E} d \hat{\Xi}$ and $\underline{\int}_{X} \chi_{E} d \Xi = \underline{\int}_{X} \chi_{E} d \hat{\Xi}$. 
\end{lemma}

\begin{PROOF}{\ref{t520.0}}
    We just deal with the result for the upper integral because the proof for the lower integral is completely analogous. The inequality $\overline{\int_{X}} \chi_{E} d \hat{\Xi} \leq \overline{\int_{X}} \chi_{E} d \Xi$ is a direct consequence of \autoref{t93}.
    
    %For this, notice that, $\bfP^{\Xi} \subseteq \bfP^{\hat{\Xi}}$. Since $\cB \subseteq \cJ^{\Xi}$, it follows that $ \{ \supsum^{\Xi}(\chi_{E}, P) \colon P \in \bfP^{\Xi} \} \subseteq \{ \supsum^{\hat{\Xi}}(\chi_{E}, P) \colon P \in \bfP^{\hat{\Xi}} \},$ hence $\overline{\int_{X}} \chi_{E} d \hat{\Xi} \leq \overline{\int_{X}} \chi_{E} d \Xi$. \DM{Esto se puede generalizar en la secci\'on previa (con \cite[Thm.~3.13]{PU} i.e.\ \autoref{t93}).}
    
    For the converse inequality,  let $P \in \bfP^{\hat{\Xi}}$. Since $P$ is a finite partition of $X$, using \autoref{t70}, we have that: $$ \supsum^{\hat{\Xi}}(\chi_{E}, P) = \sum_{b \in P} \sup(\chi_{E}[b]) \hat{\Xi}(b) = \sum_{ \substack{  b \in P \\ b \cap E \neq \emptyset}} \hat{\Xi}(b) = \hat{\Xi}(b_{E}) = \overline{\int_{X}} \chi_{b_{E}} d \Xi \geq \overline{\int_{X}} \chi_{E} d \Xi,$$ where $b_{E} \in \cJ^{\Xi}$ is the union of all elements in $P$ whose intersection with $E$ is non-empty, which includes $E$. Therefore, the latter upper integral above is a lower bound of the set $\{ \supsum^{\hat{\Xi}}(\chi_{E}, P) \colon P \in \bfP^{\hat{\Xi}} \}.$ Thus, $\overline{\int_{X}} \chi_{E} d \Xi  \leq \overline{\int_{X}} \chi_{E} d \hat{\Xi}.$  
\end{PROOF}

As a consequence of \autoref{t520.0}, we get that Jordan $\Xi$-measurability is complete in the following sense.

\begin{corollary}\label{t520}
    $(\hat{\Xi})^*=\Xi^*$, $(\hat{\Xi})_*=\Xi_*$, $\cJ^{\hat{\Xi}}= \cJ^\Xi$ and $\hat{\hat{\Xi}}=\hat{\Xi}$.
\end{corollary}
\begin{PROOF}{\ref{t520}}
    Let $E \subseteq X$.  Using \autoref{t520.0}, $(\hat{\Xi})^{\ast}(E) = \overline{\int_{X}} \chi_{E} d \hat{\Xi} = \overline{\int_{X}} \chi_{E} d \Xi = \Xi^{\ast}(E).$ Thus,  $(\hat{\Xi})^{\ast} = \Xi^{\ast}$. Analogously, $(\hat{\Xi})^{\ast} = \Xi_{\ast}$. Using this, we have: $$ E \in \cJ^{\hat{\Xi}} \Leftrightarrow (\hat{\Xi})^{\ast}(E) = (\hat{\Xi})_{\ast}(E) \Leftrightarrow \Xi^{\ast}(E) = \Xi_{\ast}(E) \Leftrightarrow E \in \cJ^{\Xi}. $$ Thus, $\cJ^{\Xi} = \cJ^{\hat{\Xi}}.$  Finally, for any $B \in \cJ^{\Xi}$, $\hat{\hat{\Xi}}(B) \coloneqq (\hat{\Xi})^{\ast}(B) = \Xi^{\ast}(B) = \hat{\Xi}(B)$. Consequently, $\hat{\hat{\Xi}} = \hat\Xi$.  
\end{PROOF}

% Consequently, we get the equivalence between $\Xi$-integrability and $\hat{\Xi}$-integrability for characteristic functions, serving as the first step toward the main result of this section (see \autoref{t560}).

\begin{corollary}\label{t522}
    Let $E \subseteq X$. Then $\chi_{E} \in \cI(\Xi)$ iff $\chi_{E} \in \cI(\hat{\Xi})$, and 
    $$
    \int_{X} \chi_{E} d \Xi = \int_{X} \chi_{E} d \hat{\Xi}.$$  
\end{corollary}

The main result in this section is a generalization of \autoref{t522}: $\Xi$-integrability and $\hat{\Xi}$-integrability are equivalent. Concretely:

\begin{theorem}\label{t560}
    Let $f\colon X\to\bbR$ be bounded. Then $f$ is $\Xi$-integrable iff it is $\hat\Xi$-integrable, and $\int_X f d \Xi = \int_X f d \hat\Xi$.
\end{theorem}

In particular, Riemann integration over rectangles in $\bbR^{n}$ is equivalent to integration with respect to the Jordan measure on $\bbR^{n}$. 

The easy direction in \autoref{t560}, namely ``$\Rightarrow$", is a direct consequence of \autoref{t93}, \autoref{t510} and \autoref{t540}. However, to show the converse, we need to develop some important tools. In particular, we need to prove a limit theorem for $\Xi$-integration. For this, we use the notion of \emph{hazy convergence} from~\cite[Sec.~4.3]{BhaskaraRa} but with a more natural name.

% \begin{theorem}\label{t38}
%     If $h \colon X \to \bbR$ is bounded and $\cJ^\Xi$-measurable, then $h$ is $\Xi$-integrable.
% \end{theorem}

% \begin{PROOF}{\ref{t38}}
%     By \autoref{t560}, it is enough to show that $h$ is $\hat\Xi$-integrable.
%     Let $\varepsilon > 0.$ Since $h$ is bounded, there are $a, b \in \bbQ$ such that $a < b$ and $f[X] \subseteq [a,b)$. Partition $[a,b)$ into finitely many intervals $\{[a_k,b_k) \colon k < m \}$ with rational endpoints and with length ${<}\frac{\varepsilon}{\hat\Xi(X)+1}$. Consider $P \coloneqq \{h^{-1}\left[[a_k,b_k)\right] \colon  k<m\}$, which is in $\Pbf^{\hat\Xi}$ because $h$ is $\cJ^\Xi$-measurable. Hence,
%     \begin{equation*}
%         \begin{split}
%             \supsum(h,P)-\infsum(h,P) & = \sum_{b\in P}(\sup h[b]-\inf h[b])\hat\Xi(b)\\
%             & \leq \sum_{b\in P}\frac{\varepsilon}{\hat\Xi(X)+1}\hat\Xi(b)=\frac{\varepsilon}{\hat\Xi(X)+1}\hat\Xi(X)<\varepsilon.
%         \end{split}
%     \end{equation*}    
%     Thus, by the criterion of $\hat\Xi$-integrability, $h \in \cI(\hat\Xi).$ 
% \end{PROOF}

\begin{definition}[Convergence in outer measure]\label{t64}
    Let $f_n,f\colon X\to \bbR$ ($n<\omega$). We say that the sequence $\la f_n\colon n<\omega\ra$ is \emph{$\Xi^*$ convergent to $f$}, denoted $f_n \rightarrow^{\Xi^*} f$, if 
    % \[\forall\, \varp,\varp'>0\ \exists\, N<\omega\, \forall\, n\geq N\ \exists\, b_n\in\cB\colon \Xi(b_n)<\varp' \text{ and } \forall\, x\in X\menos b_n\colon |f_n(x)-f(x)|<\varp.\]
    % In terms of the outer measure of $\Xi$ introduced in \autoref{t311}, the previous means that
    \[\lim_{n\to\infty}\Xi^*\left(\set{x\in X}{|f_n(x)-f(x)|\geq\varp}\right)=0 \text{ for all $\varp>0$.}\]
\end{definition}

Under the conditions of \autoref{s2}, convergence in outer measure is preserved under composition as follows. 

\begin{lemma}\label{da600}
    Let $X$, $Y$ be non-empty sets, $h \colon X \to Y$ a function, $\cB \subseteq \calP(X)$ a Boolean algebra, $\Xi$ a finite fam on $\cB$, and $f_n,f \colon Y \to \bbR$ for $n<\omega$. If $f_{n} \to^{\Xi^{\ast}_{h}} f$ then $f_{n} \circ h \to^{\Xi^{\ast}} f \circ h$. 
\end{lemma}

\begin{PROOF}{\ref{da600}}
    For any $\varp > 0$ and $n < \omega$, consider de the sets $Y_{n\varp} \coloneqq \{ y \in Y \colon \vert f_{n}(y) - f(y) \vert \geq \varp \}$ and  $X_{n, \varp} \coloneqq \{ x \in X \colon \vert g_{n}(x) - g(x) \vert \geq \varp  \}$, where $g_{n} \coloneqq f_{n} \circ h$ and $g \coloneqq f \circ h$. 
    Clearly, $X_{n,\varp} = h^{-1}[Y_{n,\varp}]$. 
    
    Let $\varp' > 0$. Since $\langle f_{n} \colon n < \omega \rangle$ $\Xi^{\ast}_{h}$-converges to $f$, we can find an $N < \omega$ such that, for any $n \geq N$, $\Xi_{h}^{\ast}(Y_{n, \varp}) < \varp'$. By \autoref{s5}~\ref{s5.a}, considering that $\chi_{X_{n,\varp}} = \chi_{Y_{n,\varp}} \circ h$, $\Xi^{\ast}(X_{n, \varp}) \leq \Xi_{h}^{\ast}(Y_{n, \varp}) < \varp'$, which proves that $\langle g_{n} \colon n < \omega \rangle$ $\Xi^{\ast}$-converges to $g$. 
\end{PROOF}

From here, we use simple functions as building blocks for our integration theory. 

\begin{definition}\label{t65.0}
    A function $f \colon X \to \bbR$ is said to be a \emph{simple function} if there are a finite partition $P$ of $X$ and a sequence $\langle c_{b} \colon b \in P \rangle$ of real numbers such that $f = \sum_{b \in P} c_{b} \chi_{b}.$ When $P\subseteq \cB$, we say that $f$ is a \emph{simple function on $\cB$}.
\end{definition}

\begin{lemma}\label{t65.01}
    Let $f = \sum_{b \in P} c_{b} \chi_{b}$ be a simple function. If $P \subseteq \cJ^{\Xi}$ then $f$ is $\Xi$-integrable. The converse is true when $b\neq b'$ implies $c_b\neq c_{b'}$ for $b,b'\in P$.
\end{lemma}

\begin{PROOF}{}
    \renewcommand{\qed}{}
    Assume that $P \subseteq \cJ^{\Xi}$. Then, by the definition, each $\chi_{b}$ is integrable and therefore, by applying \autoref{t47}, it follows that $f$ is $\Xi$-integrable.

     To see the converse, it suffices to prove the following fact. 
     \end{PROOF}

     \begin{claim}\label{t65.01clm}
       Let $A$ be a finite disjoint family of subsets of $X$ and $f = \sum_{b \in A} c_{b} \chi_{b}$ where, for $b,b'\in A$, $c_b\neq 0$, and $b\neq b'$ implies $c_b\neq c_{b'}$. If $f$ is $\Xi$-integrable, then $A\subseteq \cJ^{\Xi}$
     \end{claim}
    \begin{PROOF}{\ref{t65.01clm},\ \ref{t65.01}}
     We proceed to prove the claim by induction on $n\coloneqq |A|\geq 1$.  
     When $n=1$, i.e.\ $A=\{b_0\}$, if $f=c_{b_0}\chi_{b_{0}}$ is $\Xi$-integrable and $c_{b_0}\neq 0$, then $\chi_{b_0}=\frac{1}{c_{b_0}}f$ is $\Xi$-integrable, so $b_0\in \cJ^\Xi$.
     
     For the successor step, assume that $f = \sum_{i=0}^n c_{b_i} \chi_{b_i}$ and define $g\coloneqq f - c_{b_0}$. Notice that, for $x\in b_0$, $f(x)=c_{b_0}$, so $g(x)=0$. Since 
     $g=\sum_{i=1}^n (c_{b_i}-c_{b_0}) \chi_{b_i}$ is integrable, by the inductive hypothesis we get $\{b_1, b_2,\ldots, b_n\}\subseteq\cJ^{\Xi}$. It is only left to see that $b_0\in \cJ^{\Xi}$. To this end, let $h \coloneqq f - c_{b_1}$. 
     %Once observe that if $x\in b_1$, then $h(x)=0$. 
     Since $h=\sum_{\substack{i=0\\ i\neq 1}}^n(c_{b_i}-c_{b_1}) \chi_{b_i}$, we obtain  $\{b_0, b_2, b_3,\ldots, b_n\}\subseteq\cJ^{\Xi}$ as before. Thereby,  $\{b_0, b_1, b_2,\ldots, b_n\}\subseteq\cJ^{\Xi}$, which ends the induction. 
    \end{PROOF}

\begin{corollary}\label{t526}
    Let $f$ be a simple function. Then $f \in \cI(\Xi)$ iff $f \in \cI(\hat{\Xi})$, and $\int_{X} f d \Xi = \int_{X} f d \hat{\Xi}$.
\end{corollary}

\begin{PROOF}{\ref{t526}}
    We deal with the non-trivial implication. First, notice that, without loss of generality, we can write  $f = \sum_{b \in P} c_{b} \chi_{b}$ for some $P\in\Pbf^{\pts(X)}$ such that, for any $b \neq b'$ in $P$, $c_{b} \neq c_{b'}$. If $f \in \cI(\hat{\Xi})$ then, by using \autoref{t65.01}, $P\subseteq \cJ^{\hat{\Xi}} = \cJ^{\Xi}$. Using \autoref{t65.01} again, it follows that $f \in \cI(\Xi).$ 
\end{PROOF}

% \begin{PROOF}{\ref{t526}}
%     Assume that $f$ is a simple function. Then, there are $n^{\ast} < \omega,$ $c_{n} \in \bbR$ and $S_{n} \subseteq X$ such that $f = \sum_{n < n^{\ast}} c_{n} \chi_{S_{n}}$ and the sequence of $S_{n}$-s is pairwise disjoint. If $f \in \Xi$ then...............
% \end{PROOF}

% \begin{lemma}
%     A simple function $f = \sum_{n < n^{\ast}} c_{n} \chi_{S_{n}}$ is $\Xi$-integrable iff for any $n < n^{\ast}$, $S_{n} \in \cJ^{\Xi}.$ \Andres{No estoy seguro si esta es la caracterización que queremos. }  
% \end{lemma}

We now prove that every $\Xi$-integrable function is the $\Xi^*$-limit of simple functions on $\cB$. Recall that a family (or sequence) $\set{f_i\colon X\to \R}{i\in I}$ of functions is \emph{uniformly bounded} if there is some $M\geq 0$ such that $|f_i|\leq M$ for all $i\in I$. In this case, we say that $M$ is an \emph{uniform bound} of the family.

\begin{lemma}\label{t65}
     Every bounded $\Xi$-integrable function $f\colon X\to \R$ is the $\Xi^*$-limit of some uniformly bounded sequence of simple $\Xi$-integrable functions.\footnote{This property is called \emph{$\mathrm{T}_1$-measurable} in~\cite[Def.~4.4.5]{BhaskaraRa}. In fact,~\cite[Thm.~4.5.7]{BhaskaraRa} proves that $\mathrm{T}_1$-measurable and $\Xi$-integrable on $X$ (which they call \emph{$\mathrm{S}$-integrable}) are equivalent.} Moreover, any $M\geq 0$ is a uniform bound of the sequence whenever $|f|\leq M$.
\end{lemma}

\begin{PROOF}{\ref{t65}}
     Assume that $f \colon X \to \bbR$ is a $\Xi$-integrable function, $M\geq 0$ and $|f|\leq M$. By \autoref{t29}, any $n < \omega$ there exists some $P_{n} \in \bfP^{\Xi}$ such that $\supsum(f, P_{n}) - \infsum(f, P_{n}) < \frac{1}{n+1}$. Without loss of generality, assume that $P_{n + 1} \ll P_{n}$ whenever $n < \omega$. For any $n < \omega$, consider $f_{n} \coloneqq \sum_{b \in P_{n}} c_{b} \chi_{b}$, where $c_{b} \coloneqq \sup(f[b])$ for every $b \in P_{n}$. Notice that $|f_n|\leq M$. 
     Define $g_{n} \coloneqq f_{n} - f$, which is clearly non-negative. 
     
     Let us prove that $\langle f_{n} \colon n < \omega \rangle$ $\Xi^*$-converges to $f$. For this, let $\varp > 0$ and $\varp'>0$. We can find $N < \omega$ such that $\frac{1}{N+ 1} < \varp\, \varp'$. Define $E_{n, \varp} \coloneqq \{ x \in X \colon \vert f_{n}(x) - f(x) \vert \geq \varp \}$. Notice that $\varp \chi_{E_{n, \varp}} \leq g_{n}$ and therefore, for $n \geq N$ we have,
    \begin{equation*}
        \begin{split}
            \varp \Xi^{\ast}(E_{n, \varp}) & \leq \overline{\int_{X}} g_{n} d \Xi \leq  \supsum(g_{n}, P_{n})  = \sum_{b \in P_{n}} \sup ( g_{n}[b] ) \Xi(b) = \sum_{b \in P_{n}} [c_{b} - \inf(f[b])] \Xi(b)\\
            & = \sum_{b \in P_{n}} c_{b} \Xi(b) - \sum_{b \in P_{n}} \inf(f[b]) \Xi(b)  = \supsum(f, P_{n}) - \infsum(f, P_{n})\\
            & \leq  \supsum(f, P_{N}) - \infsum(f, P_{N}) < \frac{1}{N+1} < \varp\, \varp'.
        \end{split}
    \end{equation*}
    As a consequence, for any $n \geq N$, $\Xi^{\ast}(E_{n, \varp}) \leq \varp'$. Thus, $f$ is the $\Xi^*$-limit of the sequence $\langle f_{n} \colon n < \omega \rangle$. 
\end{PROOF}

Before stating the limit theorem, we introduce the notion of \emph{almost everywhere}. 

\begin{definition}\label{t455}
    We say that a property $\varphi$ holds \emph{$\Xi$-almost everywhere in $X$} if $\Xi^*(\{ x \in X \colon \neg\varphi(x) \}) = 0$.

    A function $f\colon X\to \R$ is \emph{bounded $\Xi$-almost everywhere} if there is some $M\geq 0$ that bounds $f$ almost everywhere, i.e.\ $\Xi^*(\set{x\in X}{|f(x)|>M})=0$.
\end{definition}

%Notice that, by \autoref{t63}~\ref{t63.a}, $\varphi$ holds almost everywhere iff $\Xi^{\ast}(\{ x \in X \colon \varphi(x) \}) = 0.$

At this point, we can fulfill the promise made at the end of \autoref{3.5}: $\Xi$-integrability and Lebesgue integrability are equivalent in complete measure spaces (of finite measure) for bounded functions. 

\begin{theorem}\label{t46}
    Assume that $(X,\cB,\Xi)$ is a complete measurable space (i.e.\ $\cB$ is a $\sigma$-algebra, $\Xi$ is $\sigma$-additive and every $\Xi^*$-measure zero set is in $\cB$). Then, any bounded real valued function on $X$ is $\Xi$-integrable iff it is $\cB$-measurable (and hence Lebesgue integrable), and both integrals coincide.
\end{theorem}
\begin{PROOF}{\ref{t46}}
    One implication was discussed after \autoref{t38}, at the end of \autoref{3.5}.

    For the converse, 
    let $f\colon X\to \R$ be a bounded $\Xi$-integrable function. By \autoref{t65}, there is some uniformly bounded sequence $\Seq{f_n}{n<\omega}$ of simple functions on $\cB$ such that $f_n\to^{\Xi^*} f$. Since the measurable space is complete, we have that $\Seq{f_n}{n<\omega}$ converges in measure to $f$ and that $f$ is $\cB$-measurable (see e.g.~\cite[\S 22, Thm.~C and~E]{halmosmeasure}).
\end{PROOF}

Recall that, in measurable spaces, a sequence convergent in measure to a function $f$ contains a subsequence that converges pointwise to $f$ almost everywhere (see e.g.~\cite[\S 22, Thm.~D]{halmosmeasure}). However, this result does not generalize in the context of fams. If $\Xi$ is a free probability fam on a subalgebra of $\pts(\omega)$ and $\Seq{f_n\colon\omega\to \R}{n<\omega}$ satisfies that $f_n(i)=0$ for all $i>n$, then $f_n\rightarrow^{\Xi^*} 0$ but $\Seq{f_n(i)}{n\geq i}$ can be quite arbitrary for any $i<\omega$.

Inspired in results from~\cite[Ch.~4]{BhaskaraRa}, we provide the following limit theorem. 
%\DM{Quer\'ia verificar si esta referencia tiene el error que ten\'iamos al principio, pero ellos asumen que $f$ es acotado de entrada. Por otra parte, el teorema abajo si se sigue de la referencia?} \Andres{4.4.18 con la equivalencia de $T_{1}$-medible y $S$-integrable (4.5.7) se me parece mucho a lo que necesitábamos originalmente. De hecho 4.4.18 acota con una función integrable. Con la nueva versión que demostró el profe no veo tan claro que se siga.}, but we include our proof for completeness.

\begin{theorem}\label{t66}
    Assume that $f_n,f\colon X\to \bbR$ ($n<\omega$) such that each $f_n$ is $\Xi$-integrable. If $\Seq{f_n}{n<\omega}$ is uniformly bounded and $\Xi^*$-converges to $f$, then:
    \begin{enumerate}[label = \normalfont (\alph*)]
        \item\label{t66-a} $f$ is bounded $\Xi$-almost everywhere.
        \item\label{t66-b} $\Seq{f_n}{n<\omega}$ is Cauchy with respect to the semi-norm $\rho_\Xi(f,g)\coloneqq \int_X |f-g|d\Xi$, i.e.\ for any $\varp>0$ there is some $N<\omega$ such that $\int_X |f_m-f_n|<\varp$ for all $m,n\geq N$.

        This implies that the sequence $\lseq{\int_X f_n d\Xi}{n<\omega}$ is Cauchy.
        \item\label{t66-c} $f$ is $\hat \Xi$-integrable (in the complement of a measure zero set) and $\lim_{n\to\infty} \int_X f_n d\Xi =\int_X f d\hat\Xi$.
        \item\label{t66-d} In each of the following cases, $f$ is $\Xi$-integrable and $\lim_{n\to\infty} \int_X f_n d\Xi =\int_X f d\Xi$.
        \begin{enumerate}[label = \normalfont (\alph{enumi}-\roman*)]
            \item\label{t66-di} $f$ is bounded (everywhere).
            \item\label{t66-dii} $b\in\cB$ for all $b\subseteq X$ satisfying $\Xi^*(b)=0$.
        \end{enumerate}
    \end{enumerate}
\end{theorem}

\begin{PROOF}{\ref{t66}}
    Let $M>0$ be an uniform bound of $\Seq{f_n}{n<\omega}$. For $\varp>0$ and $n<\omega$, define the sets  $E_{\varp} \coloneqq \{ x \in X \colon M + \varp \leq \vert f(x) \vert\}$ and $E_{n, \varp} \coloneqq \{ x \in X  \colon \vert f_{n}(x) - f(x) \vert \geq \varp \}$. 
    
    \ref{t66-a}: Let $\varp > 0$ and $\varp' > 0$. By $\Xi^*$-convergence, we can find an $N < \omega$ such that $\Xi^*(E_{n,\varp})<\varp'$ for any $n \geq N$. Clearly $E_{\varp} \subseteq E_{N, \varp}$, so $\Xi^{\ast}(E_{\varp}) < \varp'$. Since $\varp'$ is arbitrary, it follows that $\Xi^{\ast}(E_{\varp}) = 0$, which proves the desired result. 

    \ref{t66-b}: For any $n < \omega$, let $L_n\coloneqq \int_X f_n d\Xi$. Let us show that $\Seq{L_n}{n<\omega}$ is a Cauchy sequence. Let $\varp>0$ and set $\varp' \coloneqq \frac{\varp}{18M}$ and $\varp'' \coloneqq \frac{\varp}{6 (\Xi(X)+1)}$. Since $f_n \rightarrow^{\Xi^*} f$, there is some $N<\omega$ such that $\Xi^*(E_{n,\varp''})<\varp'$ for all $n\geq N$, so there is some $E'_{n,\varp''}\in\cB$ of measure ${<}\varp'$ that contains $E_{n,\varp''}$. 
    Likewise, there is some $E'\in\cB$ of measure ${<}\varp'$ containing $E_{M}$. 
    Then, for $m,n\geq N$, 
    \begin{align*}
        |L_m - L_n|  & \leq \int_X |f_m-f_n| d\Xi \leq \int_{E'} |f_m-f_n| d\Xi + \int_{X\menos E'} |f_m-f| d\Xi + \int_{X\menos E'} |f_n-f| d\Xi.
    \end{align*}
    Notice that $\int_{E'} |f_m-f_n| d\Xi \leq 2M\Xi(E') < 2M\varp'<\frac{\varp}{3}$. 
    Now,
    \begin{equation*}
        \begin{split}
            \int_{X\menos E'} |f_m-f| d\Xi & = \int_{X\menos (E_{m,\varp''}'\cup E')} |f_m-f| d\Xi + \int_{E_{m,\varp''}'\menos E'} |f_m-f| d\Xi\\
            & \leq \int_{X \setminus (E_{m,\varp''}'\cup E')} \varp'' d \Xi + \int_{E_{m, \varp''}'\menos E'} 3M d \Xi < \varp'' \Xi(X) + 3M \varp' < \frac{\varp}{3}. 
         \end{split}
    \end{equation*}

    Analogously, $\int_{X} \vert f_{n} - f \vert d \Xi < \frac{\varp}{3}$, so $|L_m-L_n|<\varp$. Therefore $\langle L_{n} \colon n < \omega \rangle$ is a Cauchy.

    Let $L\coloneqq \lim_{n\to \infty}L_n$. 

    \ref{t66-d}: First assume that $f$ is bounded everywhere. Using \autoref{t29}~\ref{t29iv}, we aim to show that $f$ is $\Xi$-integrable and $\int_X f d\Xi = L$. Let $\varp>0$ and set $\varp'\coloneqq \frac{\varp}{18M}$. By $\Xi^*$-convergence, there is some $N<\omega$ such that $\Xi^*(E_{n,\varp/(6(\Xi(X)+1))}) < \varp'$ and $|L-L_n|<\frac{\varp}{3}$ for all $n\geq N$. For such an $n$, there is some $E''\in \cB$ containing $E'_{n,\varp/(6(\Xi(X)+1))}$ of measure ${<}\varp'$. 
    By \autoref{t29}~\ref{t29iii}, there is some $P\ll \{E'',X\menos E''\}$ such that $|S^\Xi(f_n,\tau)-L_n|<\frac{\varp}{3}$ for all $\tau\in\prod P$. Now,
    \begin{align*}
        |S^\Xi(f_n,\tau)-S^\Xi(f,\tau)| & \leq \sum_{\substack{b\in P\\ b\cap E'' = \emptyset}} |f_n(\tau(b))- f(\tau(b))|\,\Xi(b) + \sum_{\substack{b\in P\\ b\subseteq E''}} |f_n(\tau(b))- f(\tau(b))|\,\Xi(b)\\
        &\leq \frac{\varp}{6(\Xi(X)+1)}\Xi(X\menos E'') + 3M \varp'<\frac{\varp}{3}.
    \end{align*}
    Therefore, $|S^\Xi(f,\tau)-L|<\varp$.

    Now, in the case~\ref{t66-dii}, we get that $E_M\in\cB$ and $\Xi(E_M)=0$. We can apply~\ref{t66-di} after restricting everything to $X\menos E_M$. Then, $f$ is $\Xi$-integrable in the complement of a measure zero set and $\lim_{n\to\infty} \int_X f_n d\Xi =\int_X f d\Xi$.

    \ref{t66-c} Direct from~\ref{t66-di} applied to $\hat\Xi$ (\autoref{t93} is also used).
\end{PROOF}

\begin{remark}
    In \autoref{t66} it is required to assume that the sequence is uniformly bounded. A counter-example can be found with the Riemann integral. For $n<\omega$, let $I_n$ be the interval $\left[1-\frac{1}{2^n},1-\frac{1}{2^{n+1}}\right)$, which has length $\frac{1}{2^{n+1}}$. Set $c_n\coloneqq 2^{n+1}$ and $f_n\coloneqq \sum_{k\leq n}c_k \chi_{I_k}$, which is a Riemann integrable simple function on $[0,1]$. Then $f_n\rightarrow^{(\Xi^1)^*} f $, where $f=\sum_{k<\omega}c_k\chi_{I_k}$ (defined on $[0,1]$). Notice that $\int_0^1 f_nd\Xi^1 = n+1$ and $f$ is not even bounded, moreover, $\lim_{x\to 1^-} \int_0^x f d\Xi^1 = \infty$. 
\end{remark}

\begin{remark}
    One would expect in \autoref{t66} that, without the assumptions in~\ref{t66-d}, $f$ is $\Xi$-integrable in the sense of \autoref{t30} and $\lim_{n\to\infty} \int_X f_n d\Xi =\int_X f d\Xi$. However, there are counter-examples in the context of the Riemann integral: Let $C\subseteq [0,1]$ be the Cantor ternary set. Recall that the set of endpoints of the intervals used in the construction of $C$, excluding $0$ and $1$, is 
    \[D\coloneqq \largeset{\sum_{i<k}\frac{s(i)}{3^{i+1}}}{0<k<\omega,\ s\colon k\to 3,\ s(k-1)\neq 0,\ \forall i<k-1\ (s(i)\neq 1)},\]
    which is dense in $C$. Notice that any $q\in D$ has a unique representation as in the displayed sum. Define $f_n\colon[0,1]\to \R$ such that, 
    \[f(x)\coloneqq
      \begin{cases}
           \frac{1}{2^n} & x\in C,\\
          1 & x\notin C.
      \end{cases}\]    
    %for $x=\sum_{i<k}\frac{s(i)}{3^{i+1}}$ in $D$, $f_n(x)\coloneqq \frac{1}{n}$, and $f_n(x)\coloneqq 1$ when $x\in [0,1]\menos D$. 
    On the other hand, define $f\colon [0,1]\to \R$ by
    \[f(x) \coloneqq 
    \begin{cases}
        k & x=\sum_{i<k}\frac{s(i)}{3^{i+1}}\in D,\\
        1 & x\in [0,1]\menos D.
    \end{cases}\]
    It is clear that $\Seq{f_n}{n<\omega}$ is uniformly bounded and converges in outer measure to $f$, the latter because $C$ has Jordan-measure zero. Also, each $f_n$ is Riemann integrable by \autoref{t526} (because it is a simple function on the Jordan algebra) and $\int_0^1 f_n(x)dx =1$. However, although $f$ is bounded almost everywhere, it is not Riemann integrable, not even in the sense of \autoref{t29}~\ref{t29iii}. The reason is that, for every $b\in \cC^1|_{[0,1]}$ of $\Xi^1$-measure zero (only $\emptyset$ and $\{1\}$ satisfies this), $f\frestr([0,1]\menos b)$ is unbounded, so \autoref{t30} applies. 

    Still, $f$ is integrable with respect to the Jordan measure by \autoref{t66}~\ref{t66-c}. This also shows that \autoref{t560} cannot be extended to unbounded functions bounded outside a measure zero set.
\end{remark}

Recall that a sequence $\langle f_{n} \colon n < \omega \rangle$ of real functions on $X$ \emph{uniformly converges to $f$} if for any $\varp > 0$, there exists $N < \omega$ such that, for any $n\geq N$, $\vert f_{n} - f \vert < \varp.$ It is clear that any uniform convergent sequence is convergent in outer measure. As a consequence, we get a well-known result in an analytic context. 

\begin{theorem}\label{t67}
    Assume that $f_n,f\colon X\to \bbR$ ($n<\omega$) such that each $f_n$ is $\Xi$-integrable. If $\langle f_{n} \colon n < \omega \rangle$ uniformly converges to $f$, then $f$ is bounded, $\Xi$-integrable, and
    \[\lim_{n\to\infty} \int_X f_n d\Xi =\int_X f d\Xi.\]
\end{theorem}
\begin{PROOF}{\ref{t67}}
    Uniform convergence implies that $f$ is bounded, i.e.\ there is some $M>0$ such that $|f|\leq M$. Then, for any $\varp>0$, $|f_n|\leq M+\varp$ for all but finitely many $n$. Thus, \autoref{t66} applies.
\end{PROOF}

We are now ready to prove \autoref{t560}, the main result in this section. 

\begin{PROOF}[Proof of \autoref{t560}]{\ref{t560}}
    The implication ``$\Rightarrow$" follows easily by \autoref{t93}. 

    For the converse, 
    assume that $f \in \cI(\hat{\Xi}).$ By \autoref{t526} and \autoref{t65}, there exists a sequence $\langle f_{n} \colon n < \omega \rangle$ of uniformly bounded simple $\Xi$-integrable functions on $X$ such that $f_{n} \rightarrow^{\hat{\Xi}^*} f$. Using \autoref{t520}, it follows that $f_{n} \rightarrow^{\Xi^*} f$.
    %\[\lim_{n\to\infty}\Xi^*\left(\set{x\in X}{|f_n(x)-f(x)|\geq\varp}\right)=0 \text{ for all $\varp>0$.}\]
    Then, we can apply \autoref{t66} to conclude that $f \in \cI(\Xi)$ and 
    \[\int_{X} f d \Xi = \lim_{n \to \infty} \int_{X} f_{n} d \Xi = \lim_{n \to \infty} \int_{X} f_{n} d \hat{\Xi} = \int_{X} f d \hat{\Xi}. \qedhere\] 
\end{PROOF}

We conclude this section by proving that, as expected, $\Xi$-null sets do not affect the value of the integral. This follows directly from the following result.

\begin{lemma}\label{t71}
    If $E\in\cJ^\Xi$ with $\hat{\Xi}(E) = 0$ and $f\colon X\to\bbR$ is bounded on $E$, then $f$ is $\Xi$-integrable on $E$ and $\int_{E} f d \Xi = 0$. 
\end{lemma}

\begin{PROOF}{\ref{t71}}
    By \autoref{t560}, it is enough to show that $f$ is $\hat{\Xi}$-integrable and $\int_E fd \hat\Xi=0$.
    Let $M < \omega$ be such that $\vert f \vert \leq M$, $P \in \bfP^{\hat\Xi}$, and consider $P_{E}$ and $\hat{P}_{E}$ as in \autoref{t62}. Notice that $$ \supsum(\chi_{E} f, \hat{P}_{E}) = \sum_{b \in \hat{P}_{E}} \sup(\chi_{E}f[b])\hat\Xi(b) = \sum_{b \in P_{E}} \sup(f[b]) \hat\Xi(b) \leq \sum_{b \in P_{E}} M \hat\Xi(b) = M \hat\Xi(E) = 0.$$  
    Similarly, $\infsum(f \chi_{E}, \hat{P}_{E}) \geq 0$, 
    so the result follows.
    %we can conclude that $\chi_{E} f \in \cI(\hat\Xi)$ and \[\int_{E} f d \hat\Xi = \int_{X} \chi_{F} f d \hat\Xi = 0. \] \qedhere    
\end{PROOF}

%We also have additivity in subsets:

\begin{corollary}\label{t77}
     Let $E\in\cJ^\Xi$ and assume that $\hat\Xi(E) = 0$. If $f\colon X\to\bbR$ is $\Xi$-integrable, then it is $\Xi$-integrable on $X\menos E$ and
     $$\int_{X} f d \Xi = \int_{X \setminus E} f d \Xi.$$ 
\end{corollary}

\begin{PROOF}{\ref{t77}}
    Immediate from  \autoref{t73} and  \autoref{t71}.%, we have: \[\int_{X} f d \Xi = \int_{X \setminus E} f d \Xi  + \int_{E} f d \Xi = 0 + \int_{X \setminus E} f d \Xi = \int_{X \setminus E} f d \Xi.\] \qedhere 
\end{PROOF}

%%%%%%%%%%%%%%%%%%%%%%%%%%%%%%%%%%%%%%%%%%%%%%%%%%%%%%%%%%%%%%%%%%%%%%%%%%%%%%%%%%%%%%%%%%%%%%%%%%%%%%%%%%%%%%%%%%%%%%%%%%%%%%%%%%%%%%%%%%%%%%%%%%%%%%%%%%%%%%%%%%%%%%%%%%%%%%%%%%%%%%%%%%%%%%%%%%%%%%%%%%%%%%%%%%

\section{Integration and the Stone space}\label{stone}

Recall that the Stone space of any Boolean algebra is a compact zero-dimensional Hausdorff topological space.

Assume that $S$ is a compact zero-dimensional topological space, let $\cp(S)$ be the Boolean algebra of clopen subsets of $S$, and let $\mu_-$ be a fam on $\cp(S)$ with $\mu_-(S)=\delta\in(0,\infty)$. Since $S$ is compact, we have that $\mu_-$ is $\sigma$-additive on $\cp(S)$. Therefore, there is a unique ($\sigma$-additive) measure $\mu$ on the completion $\cM$ of the $\sigma$-algebra on $S$ generated by the clopen sets, such that $\mu$ extends $\mu_-$ (this classical argument appears in e.g.~\cite[Pg.~21]{kamburelis}). For all this section, we fix a compact zero dimensional space $S$, a finite fam $\mu_-$ on $\cp(S)$, and reserve the symbols $\cM$ and $\mu$ for the objects defined in this paragraph.

By the Riesz Representation Theorem, $\mu_-$ can be even extended to a unique \emph{regular} measure $\mu_+$ on the Borel $\sigma$-algebra $\calB(S)$ (see e.g.~\cite[18.7.1]{Semadeni1971}), i.e.\ for any $B\in\calB(S)$:
\[\mu_+(B) = \inf\set{\mu_+(U)}{U\subseteq S \text{ open, }B\subseteq U} = \sup\set{\mu_+(K)}{K\subseteq B \text{ compact}}.\]
For any compact $K\subseteq S$, if $U\subseteq S$ is open and $K\subseteq U$ then, by compactness, there is some clopen $C\in \cp(S)$ such that $K\subseteq C\subseteq U$. So $\mu_+(K)$ must be defined by $\inf\set{\mu_-(C)}{C\in\cp(S),\ K\subseteq C}$ (which will be  $\inf\set{\mu_+(U)}{U\subseteq S \text{ open, }K\subseteq U}$ after defining $\mu_+$ with the desired property). Afterwards, $\mu_+(U)$ must be defined by $\mu_+(S\menos U)$ for any open $U\subseteq S$, which shows the uniqueness of $\mu_+$. However, $\mu_+$ will not be used in the results of this section, so we do not further discuss about it.

%\DM{Ser\'ia divertido probar directamente la existencia de esta $\mu_+$, pero puede ser algo extensivo.} \Andres{Profe, quiere que lo intentemos hacer?} \DM{Si da tiempo, s\'i. The idea is to define the family $\cC$ of sets of the form $U\cap K$ where $U$ is open and $K$ is compact. This family forms a semi-group. We can extend the $\mu_+$ so far to $\cC$: note that $U\menos K$ is an open set, so we can define $\mu_+(U\cap K)\coloneqq \mu_+(U) - \mu_+(U\menos K)$. It is enough to show that (1) this definition does not depend on $U$ and $K$ (i.e.\ $U\cap K = U'\cap K'$ implies $\mu_+(U) - \mu_+(U\menos K) = \mu_+(U') - \mu_+(U'\menos K')$) and (2) $\mu_+$ is $\sigma$-additive on $\cC$.}

We also fix, for this whole section, a set $X$, a Boolean subalgebra $\cB$ of $\calP(X)$ and a finite fam $\Xi$ on $\cB$. Recall that $\cB$ is isomorphic with the Boolean algebra of clopen subsets of the Stone space $\stone(\cB)$ via the map $a\mapsto [a]\coloneqq \set{u\in \stone(\cB)}{a\in u}$. Then, the fam $\Xi$ can be transfered to the algebra of clopen subsets of $\stone(\cB)$ in the natural way, namely $\mu^\Xi_-([a])\coloneqq \Xi(a)$ for $a\in\cB$. Since $\stone(\cB)$ is compact, we have that $\mu^\Xi_-$ is $\sigma$-additive on the algebra of clopen sets. Therefore, by the previous discussion, there is a unique ($\sigma$-additive) measure $\mu^\Xi$ on the completion $\cM^\Xi\coloneqq \cM^{\mu^\Xi_-}$ of the $\sigma$-algebra on $\stone(\cB)$ generated by the clopen sets, such that $\mu^\Xi$ extends $\mu^\Xi_-$.

We aim to present a connection between $\Xi$-integration and the Lebesgue integral on the space $\la\stone(\cB),\cM^\Xi,\mu^\Xi\ra$. %\redq{completar carreta} \Andres{Esta carreta la echamos en la introducción}

Kamburelis~\cite[Lemma~2.3]{kamburelis} shows that any real-valued continuous function on the Stone space of a complete Boolean algebra is measurable with respect to the $\sigma$-algebra generated by the clopen sets. The following lemma generalizes this result.

\begin{lemma}\label{u1}
    Let $S$ be a compact zero dimensional space. Then
    any continuous function from $S$ into $\bbR$ is measurable with respect to the $\sigma$-algebra $\sigma(\cp(S))$ generated by the clopen sets.\footnote{A further generalization is stated in~\cite[Pg.~306]{Semadeni1971}: If $X$ is compact, then the $\sigma$-algebra generated by the open $F_\sigma$ sets is the smallest $\sigma$-algebra where all the real-valued continuous functions are measurable. In the case that $X$ is also zero-dimensional, the open $F_\sigma$ subsets of $X$ are precisely the countable unions of clopen sets, and hence this smallest $\sigma$-algebra is $\sigma(\cp(X))$.}
\end{lemma}
\begin{PROOF}{\ref{u1}}
    Assume that $f\colon S\to\bbR$ is continuous, and we show that $G_y\coloneqq f^{-1}[(-\infty,y)]\in \sigma(\cp(S))$ for any $y\in\bbR$. Choose an increasing sequence $\la q_n\colon n<\omega\ra$ of rational numbers converging to $y$. Then $G_{q_n}\subseteq F_n\coloneqq f^{-1}[(\infty,q_n]]\subseteq G_{q_{n+1}}$. Now, $F_n$ is a closed set contained in an open set so, by compactness, there is some clopen set $C_n$ such that $F_n\subseteq C_n\subseteq G_{y_{n+1}}$. Therefore, $G_y=\bigcup_{n<\omega}G_{q_n}=\bigcup_{n<\omega}C_n$, which is in $\sigma(\cp(S))$.
\end{PROOF}

Contrary to what one may expect, there could be open sets in $S$ that are not in $\cM$.

\begin{example}\label{u2}
    For this example, assume that $X$ is an infinite set, $\cB=\pts(X)$, $U\subseteq \calP(X)$ is a free ultrafilter and $\Xi=\Xi_U$. Since $\cB$ is a complete Boolean algebra, its Stone space is a extremely disconnected compact space.\footnote{Which means that the closure of open sets are open.} The ultrafilter $U$ itself is a member of the Stone space, but the closed set $\{U\}$ (and also its complement) is not in $\cM$. Indeed, $\{U\}$ has outer measure $1$ because $\Xi(b)=1$ for all $b\in\cB$ such that $U\in[b]$ (i.e.\ $b\in U$). But the complement of $\{U\}$ also has outer measure $1$: if $\la b_n\colon n<\omega\ra$ is a countable sequence of members of $\cB$ such that $\set{[b_n]}{n<\omega}$ covers $\stone(\cB)\menos\{U\}$, then $\bigcap_{n<\omega}[a_n]\subseteq\{U\}$ where $a_n=X\menos b_n$. If $\bigcap_{n<\omega}[a_n]=\{U\}$ then $U$ is the ultrafilter generated by $\{a_n\colon n<\omega\}$, but this contradicts the fact that no countable family can generate a free ultrafilter. Therefore, $\bigcap_{n<\omega}[a_n]=\emptyset$, which means that $\set{[b_n]}{n<\omega}$ covers the whole space $\stone(\cB)$ and, therefore, $\sum_{n<\omega}\mu^\Xi_-([b_n])\geq 1$.

    This implies that, for $A\in\calB(X)$, $\mu^\Xi_+(A)=1$ if $U\in A$, otherwise $\mu^\Xi_+(A)=0$. Then, $\mu^\Xi_+$ naturally extends to the measure on $\pts(\stone(\cB))$ determined by the principal ultrafilter generated by $\{U\}$.
\end{example}

Given a topological space $Z$, recall that $x,y\in Z$ are \emph{indistinguishable}, denoted by $x\equiv_Z y$, iff $x$ and $y$ have the same open neighborhoods. The \emph{Kolmogorov quotient of $Z$}
$Z/{\equiv_Z}$ is always a $T_0$ space (i.e.\ every pair of points are distinguishable), and its topology is essentially the same as $Z$, namely, the open sets in $Z/{\equiv_Z}$ are precisely the sets of the form $\set{\bar x_{\equiv_Z}}{x\in G}$ for some open $G$ in $Z$, where $\bar x_{\equiv_Z}$ denotes the $\equiv_Z$-equivalence class of $x$. Note that $Z$ is $T_0$ iff it is homeomorphic with its Kolmogorov quotient.

\begin{lemma}\label{u3}
    The space $\stone(\cp(S))$ is homeomorphic with $S/{\equiv_S}$. In particular, if additionally $S$ is $T_0$, then $S$ and $\stone(\cp(S))$ are homeomorphic via the map $x\mapsto u_x=\set{C\in\cp(S)}{x\in C}$.
\end{lemma}
\begin{PROOF}{\ref{u3}}
    We first show that $\stone(\cp(S))=\set{u_x}{x\in S}$. Let $u\in \stone(\cp(S))$ and assume, towards a contradiction, that for all $x\in S$ there is some clopen $C^x\in u$ not containing $x$. Let $C_x\coloneqq S\menos C^x\notin u$, which is a clopen set containing $x$. Then, $\set{C_x}{x\in S}$ is an open cover of $S$ so, by compactness, there is some finite cover of $S$ of clopen sets not in $u$, which contradicts that $u$ is an ultrafilter on $\cp(S)$.

    Define $F\colon S\to \stone(\cp(S))$ by $F(x)\coloneqq u_x$. It is clear that $F$ is a continuous open map (open because $F[c]=[c]$ for any clopen $c$ on $S$), and $F(x)= F(y)$ iff $x\equiv_S y$ for all $x,y\in S$, which indicates that $\stone(\cp(S))$ is homeomorphic with the Kolgomorov space of $S$.
\end{PROOF}

We now look at the integration theory of $(\cB,\Xi)$ in terms of the Stone space. Bounded real-valued functions on $X$ can be pulled up to the Stone space in the following way.

\begin{definition}\label{u4}
    Let $f\colon X\to \bbR$ be a bounded function. Define the function $\hat{f}$, with domain the set of ultrafilters $u$ on $\cB$ such that $f$ is $u$-integrable (see \autoref{t17} and \autoref{t25}), by
    $$\hat{f}(u)\coloneqq \int_X f d \Xi_u.$$
\end{definition}

\begin{lemma}\label{u6}
    Let $f\colon X\to\bbR$ be bounded. Then, $\hat{f}$ is continuous on its domain.
\end{lemma}
\begin{PROOF}{\ref{u6}}
    Let $u_0\in\dom\hat f$ and $\varp>0$. Since $f$ is $u_0$-integrable, there is some $a\in u_0$ such that $\sup(f[a])-\hat f(u_0)<\varp$ and $\hat f(u_0)-\inf(f[a])<\varp$. Then $u_0\in[a]$.

    Now, for all $u\in[a]\cap \dom \hat f$, $\hat f(u)-\hat f(u_0)\leq \sup(f[a])-\hat f(u_0)<\varp$, and $\hat f(u_0) - \hat f(u)\leq \hat f(u_0) - \inf(f[a])<\varp$, so $|\hat f(u) - \hat f(u_0)|<\varp$.
\end{PROOF}

Inspired by Lebesgue-Vitali Theorem, which states that the Riemann integrable bounded functions on a closed interval of $\R$ are precisely the almost everywhere continuous functions, we prove the following theorem.

\begin{theorem}\label{u8}
    Let $f\colon X\to \bbR$ be bounded. Then $f$ is $\Xi$-integrable iff the complement of $\dom \hat f$ has $\mu^\Xi$-measure zero (in other words, $\hat f$ is continuous $\mu^\Xi$-almost everywhere).
\end{theorem}
\begin{PROOF}{\ref{u8}}
    Recall from \autoref{t25} that $\dom\hat{f}=\set{u\in\stone(\cB)}{\osc_f(u)=0}$. For $\varp>0$, notice that $[\osc_f<\varp]\coloneqq \osc_f^{-1}[[0,\varp)]$ is an open subset of $\stone(\cB)$, hence its complement $[\osc_f\geq\varp]$ is closed. Indeed, if $u\in\stone(\cB)$ and $\osc_f(u)<\varp$, then there is some $a\in u$ such that $\sup(f[a])-\inf(f[a])<\varp$, so $u\in [a]\subseteq[\osc_f<\varp]$.

    We first assume that the complement of $\dom\hat{f}$ is not a $\mu^\Xi$-measure zero set, and show that $f$ is not $\Xi$-integrable. Since $\dom\hat{f}=\bigcap_{n<\omega}[\osc_f<2^{-n}]$, we can find some $N<\omega$ such that $C\coloneqq [\osc_f\geq 2^{-N}]$ has positive outer measure with value $z$.\footnote{Remember that some closed (and open) sets in $\stone(\cB)$ may not be in $\cM^\Xi$.} Let $P\in\bfP^\Xi$ be arbitrary, and set $P'\coloneqq \set{b\in P}{[b]\cap C\neq\emptyset}$. Then 
    \[\supsum(f,P)-\infsum(f,P)\geq \sum_{b\in P'}(\sup(f[b])-\inf(f[b]))\Xi(b).\]
    For $b\in P'$, $[b]$ intersects $C$, so $\sup(f[b])-\inf(f[b])\geq 2^{-N}$. On the other hand, $\set{[b]}{b\in P'}$ covers $C$, so $\sum_{b\in P'}\Xi(b)= \sum_{b\in P'}\mu^\Xi_-([b])\geq z$.
    Thus
    \[\supsum(f,P)-\infsum(f,P)\geq 2^{-N}\sum_{b\in P'}\Xi(b)\geq 2^{-N}z>0.\]
    Therefore $\supsum(f,P)-\infsum(f,P)\geq 2^{-N}z$ for any $P\in\bfP^\Xi$. In conclusion, $f$ is not $\Xi$-integrable by \autoref{t29}.

    We now assume that the complement of $\dom\hat f$ has $\mu^\Xi$-measure zero, which implies that, for any real $r>0$, $[\osc_f\geq r]$ has measure zero. Pick some real $M>0$ such that $\sup_{x\in X}|f(x)|\leq M$. Fix $\varp>0$. For some $\varp_1>0$ and $\varp_2>0$ depending only on $\varp$ (at the end of the proof we will say what they are), there is some countable $Z\subseteq\cB$ such that $\{[b]\colon b\in Z\}$ covers $[\osc_f\geq \varp_2]$ and $\sum_{b\in Z}\Xi(b)<\varp_1$. Since $[\osc_f\geq \varp_2]$ is closed, it is compact, so we can find some finite $F\subseteq Z$ such that $\set{[b]}{b\in F}$ still covers $[\osc_f\geq \varp_2]$. We can modify $F$ to make it pairwise disjoint.
    
    Let $b'\coloneqq X\menos \bigcup F$. Since $[b']\cap [\osc_f\geq\varp_2]=\emptyset$, for any $u\in[b']$, $\osc_f(u)<\varp_2$, so there is some $c_u\in u$, $c_u\subseteq b'$, such that $\sup(f[c_u])-\inf(f[c_u])<\varp_2$. Then $\set{[c_u]}{u\in [b']}$ covers $[b']$, so, by compactness, there is some finite subcovering of $[b']$, that is, we can find some finite $F'\subseteq \cB$ such that $\bigcup F' = b'$ and $\sup(f[c])-\inf(f[c])<\varp_2$ for all $c\in F'$. We can further modify $F'$ so that it is a partition of $b'$.
    
    Consider $P\coloneqq F\cup F'$, which is in $\bfP^\Xi$. Then,
    \[\sum_{b\in F}(\sup(f[b])-\inf(f[b]))\Xi(b) \leq 2M \sum_{b\in F}\Xi(b) < 2M\varp_1\]
    and
    \[\sum_{c\in F'}(\sup(f[c])-\inf(f[c]))\Xi(c) \leq \varp_2\Xi(b')\leq \varp_2\Xi(X).\]
    Therefore $\supsum(f,P)-\infsum(f,P)<2M\varp_1 + \Xi(X)\varp_2$. So, we could have defined $\varp_1\coloneqq \frac{\varp}{4M}$ and $\varp_2 \coloneqq \frac{\varp}{2(\Xi(X)+1)}$ to conclude that $\supsum(f,P)-\infsum(f,P)<\varp$. Therefore, $f$ is $\Xi$-integrable.
\end{PROOF}

% \begin{example}
%     For this example,
%     let $\cB$ be the $\sigma$-algebra on $[0,1]$ and $\Xi$ the well-known Lebesgue measure on $\cB$. By identifying $x$ and $u_x$ for $x\in[0,1]$, we assume that $[0,1]\subseteq \stone(\cB)$, but note that $[0,1]$ is discrete as a subspace of $\stone(\cB)$. Therefore, any closed set in $\stone(\cB)$ contained in $[0,1]$ must be a finite subset of $[0,1]$, hence of $\mu^\Xi_-$-measure zero. This implies that the $\mu^\Xi$-outer measure of $\stone(\cB)\menos [0,1]$ is $1$. But the $\mu^\Xi$-outer measure of $[0,1]$ is also $1$: if $\set{[b_n]}{n<\omega}$ (with each $b_n\in\cB$) is an open cover of $[0,1]$ (identified with $\{u_x\colon x\in[0,1]\}$), then $\set{b_n}{n<\omega}$ is an open cover of (the true) $[0,1]$, so we must have $\sum_{n<\omega}\mu^\Xi_-([b_n])=\sum_{n<\omega}\Xi(b_n) \geq 1$.

%     This shows that the Stone space of $\cB$ contains an open set that is not in $\Mwf^\Xi$.
% \end{example}

The previous result can be used to show that the $\Xi$-integral is determined by the Lebesgue integral in the Stone space.

\begin{theorem}\label{u10}
    Assume that $f\colon X\to\bbR$ is bounded and $\Xi$-integrable. Then $\hat f$ is a $\cM^\Xi$-measurable and Lebesgue integrable function on its domain, and its Lebesgue integral coincides with $\int_X f d\Xi$.
\end{theorem}
\begin{PROOF}{\ref{u10}}
    %Since $\hat f$ is continuous (by \autoref{u6}), \red{we obtain by \autoref{u1} that it is $\cM|_{\dom\hat f}$-measurable.} \DM{Esto no est\'a bien aplicado. Pareciera que toca extender la medida al espacio de Borel...}
    %By \autoref{u8}, the complement of $\dom\hat f$ has $\mu^\Xi$-measure zero, so $\hat f$ can be handled as a $\cM$-measurable function (by just sending the points outside its domain to $0$).
    %Moreover, $\hat f$ is Lebesgue integrable because it is bounded $\mu^\Xi$-almost everywhere (i.e.\ in its domain).

    %Now, we prove that the Lebesgue integral of $\hat f$ coincides with $\int_X f d \Xi$.
    Since $f\in\cI(\Xi)$, we can find a $\ll$-decreasing sequence $\la P_n\, \colon n<\omega\ra$ in $\Pbf^\Xi$ such that
    $$\int_X f d \Xi = \lim_{n\to\infty} \supsum(f,P_n) = \lim_{n\to\infty}\infsum(f,P_n).$$    
    For $n<\omega$, let $\hat P_n\coloneqq \set{[b]}{b\in P_n}$, which is a partition of $\stone(\cB)$ into clopen sets. Define $\bar f_n\colon X\to\bbR$ and $\unbar{f}_n\colon X\to\bbR$ by
    $$\bar f_n\coloneqq \sum_{b\in P_n}\sup(f[b])\chi_b \text{ and }\unbar{f}_n\coloneqq \sum_{b\in P_n}\inf(f[b])\chi_b.$$
    Straightforward calculations show that $\hat{\bar{f}}_n$ and $\hat{\unbar{f}}_n$ have domain $\stone(\cB)$ and that they are the simple functions
    \[\hat{\bar{f}}_n = \sum_{b\in P_n}\sup(f[b])\chi_{[b]} \text{ and } \hat{\unbar{f}}_n = \sum_{b\in P_n}\inf(f[b])\chi_{[b]}.\]
    Clearly, the Lebesgue integral of $\hat{\bar{f}}_n$ is equal to $\int_X \bar f_n d\Xi$, which is $\supsum(f,P_n)$, and the  Lebesgue integral of $\hat{\unbar{f}}_n$ is equal to $\int_X \unbar{f}_n d\Xi$, which is $\infsum(f,P_n)$. Even more, if $m<n<\omega$ then $P_n\ll P_m$, so $\unbar{f}_m \leq \unbar{f}_n \leq \bar{f}_n \leq \bar{f}_m$ and $\hat{\unbar{f}}_m \leq \hat{\unbar{f}}_n \leq \hat{\bar{f}}_n \leq \hat{\bar{f}}_m$ (everywhere in their domains). Therefore the functions $\hat{\bar f}\coloneqq \lim_{n\to\infty} \hat{\bar{f}}_n$ and $\hat{\unbar{f}}\coloneqq \lim_{n\to\infty} \hat{\unbar{f}}_n$ (pointwise limit) are $\cM^\Xi$-measurable and, by the dominated convergence theorem,
    \begin{equation}\label{u10eq}
    \begin{split} 
        \int_{\stone(\cB)}\hat{\unbar{f}} d \mu^\Xi & =\lim_{n\to \infty} \int_{\stone(\cB)}\hat{\unbar{f}}_n d \mu^\Xi = \lim_{n\to\infty}\infsum(f,P_n) = \int_{X} f d \Xi, \text{ and }\\[1ex]
        \int_{\stone(\cB)}\hat{\bar{f}} d \mu^\Xi & =\lim_{n\to \infty} \int_{\stone(\cB)}\hat{\bar{f}}_n d \mu^\Xi = \lim_{n\to\infty}\supsum(f,P_n) = \int_{X} f d \Xi.
    \end{split}\tag{$\oplus$}
    \end{equation}
    Then $\int_{\stone(\cB)}(\hat{\bar{f}}-\hat{\unbar{f}})d \mu^\Xi = 0$, so $\hat{\bar{f}}=\hat{\unbar{f}}$ $\mu^\Xi$-almost everywhere. 
    
    On the other hand, if $u\in\dom \hat f$, then $\hat{\unbar{f}}_n(u) \leq \hat{f}(u) \leq \hat{\bar{f}}_n(u)$ (because $\unbar{f}_n \leq f\leq \bar{f}_n$ everywhere in $X$, and the $u$-integral is monotone), so $\hat{\unbar{f}}_n = \hat{f} = \hat{\bar{f}}_n$ $\mu^\Xi$-almost everywhere. Therefore $\hat f$ is $\Mcal-\Xi$-measurable and, by \eqref{u10eq}, the Lebesgue integral of $\hat{f}$ is $\int_{X} f d \Xi$.
\end{PROOF}

\begin{corollary}\label{u12}
    Let $S$ be a compact zero-dimensional space and $f\colon S\to\bbR$ bounded. Then $f$ is $\mu_-$-integrable iff $f$ is continuous $\mu$-almost everywhere. Moreover, in such a case both the $\mu_-$-integral and the Lebesgue $\mu$-integral of $f$ coincide.
\end{corollary}
\begin{PROOF}{\ref{u12}}
    Recall from the proof of \autoref{u3} that all ultrafilters on $\cp(S)$ have the form $u_x$ for all $x\in X$. Hence, by \autoref{t25}, $f$ is $u_x$-integrable iff $f$ is continuous on $x$, and \autoref{t17} indicates that $\hat f(u_x)=f(x)$.    
    
    On the other hand, by \autoref{u3}, the Stone space of $\cp(S)$ is $S/{\equiv_S}$, which has essentially the same topology as $S$, and hence, the $\sigma$-algebra $\cM$, the fam $\mu_-$ and the measure $\mu$ can be transferred to this quotient in a very natural way. Here, $\hat{f}$ is interpreted as the function on $S/{\equiv_S}$ with 
    $$\dom\hat f = \set{e\in S/{\equiv_S}}{\forall \varp>0\ \exists c\in\cp(S)\ (e\subseteq c\text{ and }\sup f[c]-\inf f[c]<\varp)}$$ 
    and, for $e\in\dom\hat f$, $\hat f(e)= f(x)$ for some (hence all) $x\in e$. Recall that $\hat f$ is continuous in its domain, which implies that $f$ is continuous on $\bigcup\dom\hat f$.
    
    By \autoref{u8}, $f$ is $\mu_-$-integrable iff the complement of $\dom\hat{f}$ has $\mu$-measure zero. The latter is equivalent to say that the complement of $\bigcup\dom\hat f$ (in $S$) has $\mu$-measure zero, which means that $f$ is continuous $\mu$-almost everywhere.

    Now, the Lebesgue integral of $\hat f$ in $S/{\equiv_S}$ can be directly identified with the Lebesgue integral of $f$, so the last part of the corollary follows by \autoref{u10}. 
\end{PROOF}

\begin{corollary}\label{u13}
    Let $f\colon\stone(\cB)\to\R$ bounded. Then $f$ is $\mu^\Xi_-$-integrable iff $f$ is continuous $\mu^\Xi$-almost everywhere, in which case both the $\mu^\Xi_-$-integral and the Lebesgue $\mu$-integral of $f$ coincide. 
\end{corollary}

\begin{corollary}\label{u132}
    Let $f\colon X\to\bbR$ bounded. Then the following statements are equivalent.
    \begin{enumerate}[label = \normalfont (\roman*)]
        \item\label{u132i} $f$ is $\Xi$-integrable.
        \item\label{u132ii} $\hat f$ is $\mu^\Xi_-$-integrable.
        \item $\hat f$ is continuous $\mu^\Xi$-almost everywhere.
    \end{enumerate}
    Moreover, the integrals in~\ref{u132i} and~\ref{u132ii} coincide with the Lebesgue $\mu^\Xi$-integral of $\hat f$.
\end{corollary}

\begin{corollary}\label{u35}
    The sets in $\cJ^{\mu_-}$ are precisely the subsets of $S$ with $\mu$-measure zero boundary.
\end{corollary}
\begin{PROOF}{\ref{u35}}
    Let $E\subseteq S$ and denote by $\partial E$ the boundary of $E$ in $S$. Without loss of generality assume that $\emptyset \subsetneq E \subsetneq S$, hence $D \coloneqq \{ s \in S \colon \chi_{E}$ is discontinuous at $s\} = \partial E$, so the result follows easily applying \autoref{u12}.  %
    %Assume that $E\in\cJ^{\mu_-}$, which means that $\chi_E$ is $\mu_-$-integrable. By \autoref{u12}. $\chi_{E}$ is continuous $\mu$-almost everywhere, that is, $\mu(\partial E) = \mu(D) = 0$. For the converse, if $\mu(\partial E) = 0$, then $\chi_{E}$ is $\mu$.
\end{PROOF}

%DEBT
% \begin{corollary}\label{u25}
%     Assume that $g\colon \stone(\cB)\to\bbR$ is bounded and $\mu_-^\Xi$-integrable. Then there is some $\Xi$-integrable function $f\colon X\to \bbR$ such that $\hat f = g$ $\mu^\Xi$-almost everywhere.
% \end{corollary}
% \begin{PROOF}{\ref{u25}}
%     Let $Z\subseteq\stone(\cB)$ be the set of discontinuous points of $g$, which has $\mu^\Xi$-measure zero.
%     For each $b\in\cB$ such that $[b]\nsubseteq Z$, define $f^*(b)\coloneqq \sup g[[b]\menos Z]$ and $f_*(b)\coloneqq \inf g[[b]\menos Z]$. Define $f\colon X\to \bbR$ by $f(x)\coloneqq \inf\set{\bar f(b)}{b\in\cB,\ x\in b,\ [b]\nsubseteq Z}$. We show that $f$ is as required, concretely, that $\stone(\cB)\menos Z\subseteq\dom\hat f$. Fix $u\in\stone(\cB)\menos Z$. By continuity, for all $\varp>0$ there is some $b_\varp\in u$ such that $|g(u)-g(v)|<\varp$ for all $v\in [b_\varp]\menos Z$. Then, for $x,y\in b_\varp$,
%     \[|f(x)-f(y)|\]
    
%     \redq{To complete (not sure yet if this is correct)}

%     %%%%Failed attempt:
%     %Define $f(x)\coloneqq g(u_x)$. By \autoref{u12}, $g$ is continuous $\mu^\Xi$-almost everywhere, that is, outside some $\mu^\Xi$-measure zero set $Z\subseteq \stone(\cB)$. We show that $u\in \dom\hat f$ (i.e.\ $f$ is $u$-integrable) for any $u\in\stone(\cB)\menos Z$. By continuity, for all $\varp>0$ there is some $b_\varp\in u$ such that $|g(u)-g(v)|<\varp$ for all $v\in [b_\varp]\menos Z$.
% \end{PROOF}

The Lebesgue-Vitali Theorem can be proved as a corollary of our results. To prove this, it is enough to show that the integration theory of the Cantor space is equivalent to the integration theory of any closed interval in $\bbR$.

%\begin{example}\label{u30}
    Recall that the \emph{Cantor space} is ${}^\omega 2=\prod_{n<\omega}\{0,1\}$ endowed with the product topology, where $\{0,1\}$ has the discrete topology.
    This is a compact metrizable space, it is perfect and zero-dimensional. The \emph{basic clopen sets} have the form $[s]\coloneqq \set{x\in{}^\omega 2}{s\subseteq x}$ for $s\in{}^{<\omega}2$ (the set of finite binary sequences). Hence, the clopen subsets of ${}^\omega2$ are precisely the finite unions of basic clopen sets.

    There is a unique fam $\Lb_-\colon\cp(\cantor
    )\to[0,1]$ such that $\Lb_-([s])=2^{-n}$ for any $s\in{}^{<\omega}2$ of length $n$ (so it is a probability fam). This is actually a $\sigma$-additive measure, so it can be extended to a probability measure $\Lb$ on the completion of its Borel $\sigma$-algebra, which is known as the \emph{Lebesgue measure of the Cantor space}.
%\end{example}

The Lebesgue measure theory of the Cantor space is equivalent to the measure theory of the unit interval $[0,1]$, via the function $\iota_2\colon \cantor\to[0,1]$ defined by
\[\iota_2(x)=\sum_{n<\omega}\frac{x(n)}{2^{n+1}}.\]
This function is continuous, and homeomorphism modulo countable sets, namely, there are countable sets $Q_2\subseteq\cantor$ and $Q_2^*\subseteq[0,1]$ such that $\cantor\menos Q_2$ is homeomorphic with $[0,1]\menos Q^*_2$ via $\iota_2$. This leads to the equivalence of their measure theory, according to~\cite[Ch.~VII, Prop.~3.19]{Levy} and~\cite[Sec.~6]{MU25}. 

\begin{theorem}\label{u40}
    Let $A\subseteq\cantor$ and $B\subseteq[0,1]$.
    \begin{enumerate}[label=\rm (\alph*)]
        \item\label{u40a} $A$ is Lebesgue measurable in $\cantor$ iff $\iota_2[A]$ is Lebesgue measurable in $[0,1]$, in which case they have the same measure.
        
        \item\label{u40b} $B$ is Lebesgue measurable in $[0,1]$ iff $\iota_2^{-1}[B]$ is Lebesgue measurable in $\cantor$, in which case they have the same measure.
    \end{enumerate}
\end{theorem}

As a consequence, $\Lb_-$-integration in $\cantor$ is equivalent to Riemann integration in $[0,1]$, and the same for Lebesgue integration.

\begin{corollary}\label{u45}
    Let $f\colon\cantor\to\bbR$ and $g\colon[0,1]\to\bbR$ such that $g\circ\iota_2 =f$. Then:
    \begin{enumerate}[label=\normalfont(\alph*)]
        \item\label{u45a} $f$ is Lebesgue measurable iff $g$ is Lebesgue measurable.
        
        \item\label{u45b} $f$ is Lebesgue integrable iff $g$ is Lebesgue integrable, and their Lebesgue integrals coincide.

        \item\label{u45j} For $B\subseteq [0,1]$, $B\in \cJ^{\Xi^1|_{[0,1]}}$ iff $\iota_2^{-1}[B]\in \cJ^{\Lb_-}$, and their Jordan measure coincide.
        
        \item\label{u45c} When $g$ is bounded,
        $f$ is $\Lb_-$-integrable iff $g$ is Riemann integrable, and their integrals coincide.
        
        \item\label{u45d} For $A\subseteq\cantor$, if $A\in\cJ^{\Lb_-}$ then $\iota_2[A]\in\cJ^{\Xi^1|_{[0,1]}}$, and their Jordan measures coincide. %DM: se cumple el rec\'iproco? Veo la medida exterior, pero no la interior
    \end{enumerate}
\end{corollary}
\begin{PROOF}{\ref{u45}}
    \ref{u45a}: Immediate by \autoref{u40}. %~\ref{u40a} and~\ref{u40b}.

     \ref{u45b}: 
     Assume that $g$ is Lebesgue integrable and, without loss of generality, suppose that it is non-negative. Let $\sigma \colon \cantor \to \bbR$ be a simple function on the Lebesgue $\sigma$-algebra such that $\sigma \leq f$. Define $\tau \colon [0, 1] \to \bbR$ by
    \[
    \tau(y) \coloneqq
    \begin{cases}
    \sigma(\iota_2^{-1}(y)), & \text{if } y \in [0, 1] \setminus Q_{2}^{\ast}, \\
    0, & \text{if } y \in Q_{2}^{\ast}, \\
    \end{cases}
    \]
    which is a simple function on the Lebesgue $\sigma$-algebra. 
    If $y \in [0,1] \setminus Q_{2}^{\ast}$, then $ \tau(y) \leq f(\iota_2^{-1}(y)) = (g \circ \iota_2)(\iota_2^{-1}(y)) = g(y),$ hence $\tau \leq g$ and, by \autoref{u40}, $\int_{\cantor} \sigma d \Lb = \int_{[0, 1]} \tau d \lambda$. As a consequence,
    \begin{equation}\label{e.u45b}
        \begin{split}
            \int_{\cantor} f d \Lb & \leq % \sup \left\{ \int_{\cantor} \sigma d \Lb \colon \sigma  \leq f \text{ and } \sigma \text{ is simple}  \right \}\\
            %& \leq \sup  \left\{ \int_{[0, 1]} \tau d \lambda \colon \tau \leq g \text{ and } \tau \text{ is simple}  \right  \} 
             \int_{[0, 1]} g d \lambda < \infty.
        \end{split}
    \end{equation}
    This shows that $f$ is Lebesgue integrable. 

    Conversely, assume that $f$ is Lebesgue integrable, and without loss of generality, assume that it is non-negative.  Let $\tau \colon [0, 1] \to \bbR$ be a simple function on the Lebesgue $\sigma$-algebra such that $\tau \leq g$. Define $\sigma \colon \cantor \to \bbR$ by $\sigma\coloneqq \tau\circ \iota_2$. 
    % \[
    % \sigma(x) \coloneqq
    % \begin{cases}
    % \tau(\iota_{2}(x)), & \text{if } x \in \cantor \setminus Q_{2} \\
    % 0, & \text{if } x \in Q_{2}. \\
    % \end{cases}
    % \]
    It is not hard to check that $\sigma \leq f$ is a simple function on the Lebesgue $\sigma$-algebra and $\int_{\cantor} \sigma d \Lb = \int_{[0, 1]} \tau d \lambda$. As a consequence,  
        \begin{equation}\label{e.u45b1}
        \begin{split}
            \int_{[0, 1]} g d \lambda & \leq % \sup \left\{ \int_{\cantor} \sigma d \Lb \colon \sigma  \leq f \text{ and } \sigma \text{ is simple}  \right \}\\
            %& \leq \sup  \left\{ \int_{[0, 1]} \tau d \lambda \colon \tau \leq g \text{ and } \tau \text{ is simple}  \right  \} 
             \int_{\cantor} f d \Lb < \infty.
        \end{split}
    \end{equation}
    This shows that $g$ is Lebesgue integrable and, by~(\ref{e.u45b}) and~(\ref{e.u45b1}), the integrals coincide.  

    \ref{u45j} ($\Rightarrow$): If $B\in\cJ^{\Xi^1|_{[0,1]}}$ then $\chi_B$ is continuous except on a finite set, so the same is true for $\chi_B\circ \iota_2=\chi_{\iota^{-1}_2[B]}$. Hence, by \autoref{u12}, $\chi_{\iota^{-1}_2[B]}$ is $\Lb_-$-integrable, i.e.\ $\chi_{\iota^{-1}_2[B]} \in \cJ^{\Lb_-}$. The Jordan measures coincide by \autoref{u40}.
    
    \ref{u45c}: Assume that $g$ is a Riemann integrable bounded function. By \ref{u45j} ($\Rightarrow$), $\cC^1|_{[0,1]}\subseteq \iota^\to(\cJ^{\Lb_-})$ so, by \autoref{t93}, \autoref{s2},~\ref{t560} and~\ref{u40}, $f=\iota_2\circ g$ is $\Lb_-$ integrable and $\int_{\cantor} f d\Lb_- = \int_0^1 g(y) dy$, 
    
    %By ~\autoref{t65}, we can find a sequence $\langle \tau_{n} \colon [0, 1] \to \bbR  \colon n < \omega \rangle$ of simple functions on $\cC^1|_{[0,1]}$  that $(\Xi^1|_{[0,1]})^{\ast}$-converges to $g$. For any $n < \omega$, consider $\sigma_n \colon \cantor \to \bbR$ defined by $\sigma_{n} \coloneqq \tau_{n} \circ \iota_{2}$. \red{It is not hard to check} that $\langle \sigma_{n} \colon n < \omega \rangle$ $\Lb^{\ast}_-$-converges to $f$, \DM{se puede hacer un lema general de eso.} and that each $\sigma_{n}$ is continuous $\Lb$-almost everywhere (because $\tau_n$ is continuous $\lambda$-everywhere and by \autoref{u40}). Therefore, by~\autoref{u12}, $\sigma_n$ is $\Lb_{-}$-integrable. This implies, by~\autoref{t66}, that $f$ is $\Lb_{-}$-integrable.

    To prove the converse, assume that $f$ is $\Lb_{-}$-integrable and let $\langle \sigma_{n} \colon \cantor \to \bbR \colon n < \omega \rangle$ be a sequence of simple functions on $\cp(\cantor)$ $\Lb^{\ast}_-$-converging to $f$. Fix $n < \omega$. We can express $\sigma_n=\sum_{i<k_n}c_{n,i}\chi_{[s_{n,i}]}$ where each $s_{n,i}\in {}^{<\omega}2$ and $H\coloneqq\set{[s_{n,i}]}{i<k_n}$ covers $2^\omega$. For each $s\in{}^{<\omega}2$, $\iota_2[[s]]$ is a closed interval contained in $[0,1]$, so $\set{\iota_2[[s_{n,i}]]}{i<k_n}$ covers $[0,1]$ and two different of these interval intersects in at most one endpoint. Define $\tau_n\colon [0,1]\to\R$ by
    \[\tau_n(y) \coloneqq
      \begin{cases}
          \sigma_n(x) & \text{if $y=\iota_2(x)$ is in the interior of some $\iota_2[[s_{n,i}]]$,}\\
          0 & \text{otherwise.}
      \end{cases}\]
    This map is well-defined and a simple function on the Jordan algebra, hence Riemann integrable by \autoref{t526}. 
    
    We show that $\langle \tau_{n} \colon n < \omega \rangle$ $(\Xi^1|_{[0,1]})^{\ast}$-converges to $g$. Let $\varp>0$ and $\varp'>0$, so there is some $N<\omega$ such that $\Lb_-^*(\set{x\in\cantor}{|f(x)-\sigma_n(x)|\geq\varp})<\varp'$ for all $n\geq N$. For such an $n$, there is some $C_n\in\cp(\cantor)$ containing $\set{x\in\cantor}{|f(x)-\sigma_n(x)|\geq\varp}$ with $\Lb_-(C_n)<\varp'$. 
    For $y\in [0,1]$, if $y$ is in the interior of some $\iota_2[[s_{n,i}]]$ then $\tau_n(y)=\sigma_n(x)$ for some $x\in\cantor$ such that $y=\iota_2(x)$, so $|g(y) - \tau_n(y)| = |f(x)-\sigma_n(x)|$; on the other hand, $y$ is the set of endpoints of $\iota_2[[s_{n,i}]]$ for $i<k_n$, which we denote by $F_n$. Hence, $\set{y\in [0,1]}{|g(y)-\tau_n(y)|\geq\varp}$ is contained in $\iota_2(C_n)\cup F_n$. Since $\iota_2(C_n)$ is a union of closed intervals, $\iota_2(C_n)$ is in the Jordan algebra, as well as $F_n$ (because it is finite), with Jordan measure $\Lb_-(C_n)$ (by \autoref{u40}) and $0$, respectively. Thus, 
    $(\Xi^1|_{[0,1]})^*(\set{y\in [0,1]}{|g(y)-\tau_n(y)|\geq\varp})<\varp'$.
    
    Therefore, by \autoref{t66}, it follows that  $g$ is Riemann integrable.

    \ref{u45j} ($\Leftarrow$): Direct by~\ref{u45c} because $\chi_B\circ \iota_2=\chi_{\iota^{-1}_2[B]}$.

    \ref{u45d}:  
    %Let $A \subseteq \cantor$. 
    Assume $A\in\cJ^{\Lb_-}$. For any $\varp>0$, there are $C,D\in\cp(\cantor)$ such that $C\subseteq A\subseteq D$ and $\Lb_-(D\menos C)<\varp$. Since the image of a clopen subset of $\cantor$ is a finite union of closed intervals, $\iota_2[C],\iota_2[D]\in \cJ^{\Xi^1|_{[0,1]}}$. Then, by \autoref{u40}, and using that $\iota_2$ is a homomorphism except on a countable set, $\hat \Xi^1(\iota_2(D\menos C)) = \hat\Xi^1(\iota_2[D])-\hat\Xi^1(\iota_2[C]) = \Lb_-(D)-\Lb_-(C) <\varp$. Hence, we can find $C',D'\in\cC^1|_{[0,1]}$ such that $C'\subseteq \iota_2[C] \subseteq \iota_2[A] \subseteq \iota_2[D] \subseteq D'$ and $\Xi^1(D'\menos C')<\varp$. This shows that $\iota_2[A]\in\cJ^{\Xi^1|_{[0,1]}}$. By \autoref{u45}, the Jordan measures of $A$ and $\iota_2[A]$ coincide.
    %
    % Consider $A^{\ast} \coloneqq A \setminus Q_{2}$ and $B^{\ast} \coloneqq \iota_{2}[A] \setminus Q_{2}^{\ast}$. It is clear that $\chi_{B^{\ast}} \circ \iota_{2} = \chi_{A^{\ast}}$. 
    % \begin{equation*}
    %     \begin{split}
    %         A \in \cJ^{\Lb_{-}} & \Leftrightarrow A^{\ast} \in \cJ^{\Lb_{-}} \Leftrightarrow \chi_{A^{\ast}} \in \cI(\Lb_{-}) \Leftrightarrow \chi_{B^{\ast}} \in \cI(\Xi^{1} |_{[0, 1]})\\
    %         & \Leftrightarrow B^{\ast} \in \cJ^{\Xi^{1} |_{[0, 1]}} \Leftrightarrow \iota_{2}[A] \in \cJ^{\Xi^{1} |_{[0, 1]}}. 
    %     \end{split}
    % \end{equation*}
\end{PROOF}

\begin{remark}\label{da601}
    In contrast to the situation in \autoref{u45}~\ref{u45j}, the converse implication in \autoref{u45}~\ref{u45d} does not hold in general. To see this, consider the set  
    $ A \coloneqq \cantor \setminus A^{\ast},$ where $A^{\ast}$ consists of sequences $s \in \cantor$ for which there exists some $N < \omega$ such that $s(N) = 1$ and $s(n) = 0$ for all $n > N$. Notice that $A \notin \cJ^{\Lb_{-}}$, as its outer measure is $1$ while its inner measure is $0$. However, $\iota_{2}[A] = [0,1] \in \cJ^{\Xi^{1} |_{[0,1]}}$.  
\end{remark}

So Riemann-Vitali Theorem follows.

\begin{corollary}\label{u50}
    Let $g\colon[a,b]\to \bbR$ be a bounded function. Then $g$ is Riemann integrable iff $g$ is continuous almost everywhere, in which case the Riemann and the Lebesgue integral coincide.
\end{corollary}
\begin{PROOF}{\ref{u50}}
    Without loss of generality, assume that $a = 0$ and $b = 1$. Assume that $g$ is Riemann integrable and consider $f \coloneqq g \circ \iota_2$. By \autoref{u45}~\ref{u45c}, we have that $f$ is $\Lb_{-}$-integrable, by \autoref{u12}, it follows that $f$ is continuous $\Lb$-almost everywhere, and therefore $g$ is continuous $\lambda$-almost everywhere. To prove the converse, assume that $g$ is continuous almost everywhere. Then $f$ is continuous $\Lb$-almost everywhere, and by \autoref{u12}, it follows that $f$ is $\Lb_{-}$-integrable. The conclusion follows from \autoref{u45}~\ref{u45c}.
\end{PROOF}

%\DM{Habr\'ia alguna forma (sencilla) de sacar a \autoref{u50} para $\R^n$ como corolario?} \Andres{Yo esto no lo veo de forma sencilla. Estaba intentando pensar si es posible hacer algo parecido a lo que hicimos entre $[0, 1]$ y $\cantor$, pero no puede existir una biyección continua entre $[0, 1]$ y $[0, 1] \times [0, 1]$, por ejemplo. Me hubiera gustado tener más tiempo para pensar en esto!} \DM{Basta con mirar a $(\cantor)^n$ y hacer el $n$-producto de $\iota_2$. Se debe seguir f\'acil as\'i.} \Andres{Esto implicaría hacer \autoref{u45} para el caso de $\bbR^{n}$, no?} \DM{De cierta forma, pero gracias a \autoref{u45} se ahorra mucho trabajo... creo. Bueno, no, hay que meterse con productos}
 
As a consequence of our results, we can easily find a counter-example of the converse of \autoref{t38}.

\begin{example}\label{u20}
    Let $\la x_n\colon n<\omega\ra$ be an increasing sequence of positive real numbers converging to $1$. Consider $S\coloneqq \set{x_n}{n<\omega}\cup\{1\}$ as a subspace of $\R$. Then $S$ is compact, Hausdorff and zero dimensional. Note that the clopen subsets of $S$ are precisely the finite subsets of $\set{x_n}{n<\omega}$ and their complements in $S$. Let $\mu_-$ be a fam on $\cp(S)$ such that all finite subsets of $\set{x_n}{n<\omega}$ have measure $0$, and their complements have measure $1$.

    The identity function $\id_S$ on $S$ is continuous, and hence, $\mu_-$integrable by \autoref{u12}. However, $\set{z\in S}{\id_S(z)=z<1}=\set{x_n}{n<\omega}$ is not in $\cJ^{\mu_-}$ because it has $\mu_-$-inner measure $0$ and $\mu_-$-outer measure $1$.
\end{example}

\bibliographystyle{alpha}

%\bibliography{appl}
\end{document}